\catcode`\@=\active
\magnification=1200

\tolerance10000
\hfuzz500pt
\def\cmp #1{{Commun.  Math.  Phys.} {\bf #1}}

\def\jmp #1{{J. Math. Phys.} {\bf #1}}
\def\pl #1{{Phys. Lett.} {\bf #1}}
\def\np #1{{Nucl. Phys.} {\bf #1}}

\def\pr #1{{Phys. Rev.} {\bf #1}}

\newcount\sectnum
\global\sectnum=0
\newcount\subnum
\newcount\convnum
\convnum=1
\newcount\noreffe

\def\reffe{\advance\noreffe by 1}
\newcount\nota

\def\cita{\global\advance\nota by 1}
\def\cito{{[\the\nota]}} 
 \newcount\thmnum  
\global\thmnum=0
\newcount\subsec
\global\subsec=0
\newcount\eqnum
\global\eqnum=0
\newcount\norefe

\def\refe{\advance\norefe by 1}
 \newcount\formu
\def\clearformu{\eqnum=1}

%FORMULE

\def\ref#1{\num  \xdef#1{(\the\sectnum.\the\eqnum)}
  }
\def\eqalignref#1{\eqalignnum  \xdef#1{( \the\sectnum.\the\eqnum)}}

\def\today{\ifcase\month\or January\or February\or March\or
        April\or May\or June\or July\or August\or September\or
        October\or November\or December\fi\space\number\day,
        \number\year}

%FILETEMPORANEI

\newwrite\fileack
\immediate\openout\fileack=ref.tmp2
\immediate\write\fileack{\parindent 30pt}

\newwrite\fileref
\immediate\openout\fileref=ref.tmp1
\immediate\write\fileref{\parindent 30pt}

\def\ack#1{ \hskip-11pt
\immediate\write\fileack
           {\par\hangindent\parindent}
\ignorespaces
              \immediate\write\fileack{ {#1}}
\ignorespaces}
 
\def\Ascr{{\Cal A}}

\def\Gscr{{\Cal G}}
\def\Hscr{{\Cal H}}

\def\Lscr{{\Cal L}}

\def\Oscr{{\Cal O}}
\def\Pscr{{\Cal P}}

\def\tr{{\text{\rm Tr}\;}}
\def\sqr#1#2{{\vcenter {\hrule height.#2pt
       \hbox {\vrule width.#2pt height#1pt \kern#1pt
       \vrule width.#2pt}
       \hrule height.#2pt}}}

%FONTI

\font\subsecfnt=cmbx10 at 11.20 truept 
\font\titfnt=cmbx10 at 14.40 truept
 at 13.00 truept
\font\tenmsx=msam10
\font\sevenmsx=msam7
\font\fivemsx=msam5
\font\tenmsy=msbm10
\font\sevenmsy=msbm7
\font\fivemsy=msbm5
\newfam\msxfam
\newfam\msyfam
\textfont\msxfam=\tenmsx  \scriptfont\msxfam=\sevenmsx
  \scriptscriptfont\msxfam=\fivemsx
\textfont\msyfam=\tenmsy  \scriptfont\msyfam=\sevenmsy
  \scriptscriptfont\msyfam=\fivemsy
 
%SPAZIATURE

\def\blank{\vskip 12pt}
\def\blankii{\blank\blank}
\def\blankm{\vskip 6pt} 
\def\rt#1{\hfill{#1}}
\def\lora{\longrightarrow}
\def\rsa{\rightsquigarrow}

\def\Lra{\Leftrightarrow}

\def\endproof{\rt{$\squaro$}\par\blankii}
%SUDDIVISIONE IN CAPITOLI, SEZIONI E SOTTOSEZIONI

\def\sezione#1{\subnum=0\def\sectname{#1}\global\advance\sectnum by 1\thmnum=1
\goodbreak\vskip 28pt plus 30 pt\noindent
               { \bf\the\sectnum .$\quad$ {#1}}
  \vskip15pt \clearformu}

\def\prop#1#2{\global\advance\thmnum by 1
	\xdef#1{Proposition \the\sectnum.\the\thmnum}
	\bigbreak\noindent{\bf Proposition \the\sectnum.\the\thmnum.}
	{\it#2} }

\def\definizione#1#2{
	\xdef#1{Definition  \the\sectnum.\the\thmnum}
	\bigbreak\noindent{\bf Definition   \the\sectnum.\the\thmnum.} 
	{\it#2}\global\advance\thmnum by 1 }

\def\notazione#1#2{
	\xdef#1{Notations  \the\sectnum.\the\thmnum}
	\bigbreak\noindent{\bf Notations  \the\sectnum.\the\thmnum.} 
	{\it#2}\global\advance\thmnum by 1 }

 %TEOREMI ETC

\def\convention#1#2{
        \xdef#1{Convention  {\bf\romannumeral\convnum}}
        \bigbreak\noindent{\bf Convention \romannumeral\convnum.}
        {\it#2}\global\advance\convnum by 1 }

\def\remark#1#2{
        \xdef#1{Remark  \the\sectnum.\the\thmnum}
        \bigbreak\noindent{\bf Remark  \the\sectnum.\the\thmnum.}
        {\it#2}\global\advance\thmnum by 1 }

 \def\claim#1#2{
        \xdef#1{Claim  \the\sectnum.\the\thmnum }
        \bigbreak\noindent{\bf Claim  \the\sectnum.\the\thmnum.}
        {\it#2}\global\advance\thmnum by 1 }
\def\theorem#1#2{
        \xdef#1{Theorem  \the\sectnum.\the\thmnum }
        \bigbreak\noindent{\bf Theorem  \the\sectnum.\the\thmnum.}
        {\it#2}\global\advance\thmnum by 1 }

\def\proof#1{\vskip10pt{\it Proof of
#1}}

\def\lemma#1#2{
	\xdef#1{Lemma  \the\sectnum.\the\thmnum}
	\bigbreak\noindent{\bf Lemma  \the\sectnum.\the\thmnum.} 
	{\it#2}\global\advance\thmnum by 1}

\def\cor#1#2{\global\advance\thmnum by 1
	\xdef#1{Corollary \the\sectnum.\the\thmnum}
	\bigbreak\noindent{\bf Corollary \the\sectnum.\the\thmnum.} 
	{\it#2} }
\def\proof#1{\vskip10pt{\it Proof of 
#1}}

\def\num{\global\advance\eqnum by 1
	\eqno({ \the\sectnum.\the\eqnum)}}

\def\eqalignnum{\global\advance\eqnum by 1
	({\rm\the\sectnum}.\the\eqnum)}

\def\ref#1{\num  \xdef#1{(\the\sectnum.\the\eqnum)}}
\def\eqalignref#1{\eqalignnum  \xdef#1{( \the\sectnum.\the\eqnum)}}

\def\title#1{\centerline{\bf\titfnt#1}}

\def\Alph#1{\ifcase#1\or A\or B\or C\or D\or E\or F\or G\or H\fi}
\def\sottosezione#1{\global\advance\subnum by 1 
	\vskip12pt plus4pt \goodbreak \noindent
	{ \bf\the\sectnum.{\the\subnum}.}\quad{\bf \subsecfnt #1 }\vskip10pt  }
\def\newsubsec#1{\global\subnum=1 \vskip6pt\noindent  
	{ \the\sectnum.\Alph\subnum.}{\bf #1}  }
\def\today{\ifcase\month\or January\or February\or March\or 
	April\or May\or June\or July\or August\or September\or 
	October\or November\or December\fi\space\number\day, 
	\number\year}
 
%CITAZIONI
\newwrite\fileref
\immediate\openout\fileref=ref.tmp1
\immediate\write\fileref{\parindent 30pt}

%SIMBOLI MATEMATICI

\def\foot#1{\cita\plainfootnote{$^\cito$}{#1}}

\def\mapdown#1{\Big\downarrow\rlap{$\vcenter{\hbox{$\scriptstyle#1$}}$}}

\def\mapdownn#1
{\bigg\downarrow\rlap{$\vcenter{\hbox{$\scriptstyle#1$}}$}}

\def\riga#1#2{\mathop{ \buildrel #1
 \over{\hbox to #2pt{\rightarrowfill}}}}
\def\liga#1#2{\mathop{ \buildrel #1
 \over{\hbox to #2pt{\leftarrowfill}}}}

 \def\incul{\hookrightarrow}
\def\rsa{\rightsquigarrow}

\def\dota{ {\buildrel \circ\over{ \!  A}}}

\def\dy{\displaystyle}
\def\ad{\mathop{{\text{\rm {ad}}}}}
\def\Ad{\mathop{{\text{\rm Ad}}}}

\def\adP{\mathop{{\text{\rm ad}} \, P}}

\def\form#1{\Omega^{#1}(M, \text {ad} P)}

\def\hexnumber@#1{\ifcase#1 0\or1\or2\or3\or4\or5\or6\or7\or8\or9\or
	A\or B\or C\or D\or E\or F\fi }

\edef\msx@{\hexnumber@\msxfam}
\edef\msy@{\hexnumber@\msyfam}

 \mathchardef\boxdot="2\msx@00
\mathchardef\boxplus="2\msx@01
\mathchardef\boxtimes="2\msx@02
\mathchardef\square="0\msx@03
\mathchardef\blacksquare="0\msx@04
\mathchardef\centerdot="2\msx@05
\mathchardef\lozenge="0\msx@06
\mathchardef\blacklozenge="0\msx@07
\mathchardef\circlearrowright="3\msx@08
\mathchardef\circlearrowleft="3\msx@09
\mathchardef\rightleftharpoons="3\msx@0A
\mathchardef\leftrightharpoons="3\msx@0B
\mathchardef\boxminus="2\msx@0C
\mathchardef\Vdash="3\msx@0D
\mathchardef\Vvdash="3\msx@0E
\mathchardef\vDash="3\msx@0F
\mathchardef\twoheadrightarrow="3\msx@10
\mathchardef\twoheadleftarrow="3\msx@11
\mathchardef\leftleftarrows="3\msx@12
\mathchardef\rightrightarrows="3\msx@13
\mathchardef\upuparrows="3\msx@14
\mathchardef\downdownarrows="3\msx@15
\mathchardef\upharpoonright="3\msx@16

\mathchardef\downharpoonright="3\msx@17
\mathchardef\upharpoonleft="3\msx@18
\mathchardef\downharpoonleft="3\msx@19
\mathchardef\rightarrowtail="3\msx@1A
\mathchardef\leftarrowtail="3\msx@1B
\mathchardef\leftrightarrows="3\msx@1C
\mathchardef\rightleftarrows="3\msx@1D
\mathchardef\Lsh="3\msx@1E
\mathchardef\Rsh="3\msx@1F
\mathchardef\rightsquigarrow="3\msx@20
\mathchardef\leftrightsquigarrow="3\msx@21
\mathchardef\looparrowleft="3\msx@22
\mathchardef\looparrowright="3\msx@23
\mathchardef\circeq="3\msx@24
\mathchardef\succsim="3\msx@25
\mathchardef\gtrsim="3\msx@26
\mathchardef\gtrapprox="3\msx@27
\mathchardef\multimap="3\msx@28
\mathchardef\therefore="3\msx@29
\mathchardef\because="3\msx@2A
\mathchardef\doteqdot="3\msx@2B

\mathchardef\triangleq="3\msx@2C
\mathchardef\precsim="3\msx@2D
\mathchardef\lesssim="3\msx@2E
\mathchardef\lessapprox="3\msx@2F
\mathchardef\eqslantless="3\msx@30
\mathchardef\eqslantgtr="3\msx@31
\mathchardef\curlyeqprec="3\msx@32
\mathchardef\curlyeqsucc="3\msx@33
\mathchardef\preccurlyeq="3\msx@34
\mathchardef\leqq="3\msx@35
\mathchardef\leqslant="3\msx@36
\mathchardef\lessgtr="3\msx@37
\mathchardef\backprime="0\msx@38
\mathchardef\risingdotseq="3\msx@3A
\mathchardef\fallingdotseq="3\msx@3B
\mathchardef\succcurlyeq="3\msx@3C
\mathchardef\geqq="3\msx@3D
\mathchardef\geqslant="3\msx@3E
\mathchardef\gtrless="3\msx@3F
\mathchardef\sqsubset="3\msx@40
\mathchardef\sqsupset="3\msx@41
\mathchardef\vartriangleright="3\msx@42
\mathchardef\vartriangleleft="3\msx@43
\mathchardef\trianglerighteq="3\msx@44
\mathchardef\trianglelefteq="3\msx@45
\mathchardef\bigstar="0\msx@46
\mathchardef\between="3\msx@47
\mathchardef\blacktriangledown="0\msx@48
\mathchardef\blacktriangleright="3\msx@49
\mathchardef\blacktriangleleft="3\msx@4A
\mathchardef\vartriangle="0\msx@4D
\mathchardef\blacktriangle="0\msx@4E
\mathchardef\triangledown="0\msx@4F
\mathchardef\eqcirc="3\msx@50
\mathchardef\lesseqgtr="3\msx@51
\mathchardef\gtreqless="3\msx@52
\mathchardef\lesseqqgtr="3\msx@53
\mathchardef\gtreqqless="3\msx@54
\mathchardef\Rrightarrow="3\msx@56
\mathchardef\Lleftarrow="3\msx@57
\mathchardef\veebar="2\msx@59
\mathchardef\barwedge="2\msx@5A
\mathchardef\doublebarwedge="2\msx@5B
\mathchardef\angle="0\msx@5C
\mathchardef\measuredangle="0\msx@5D
\mathchardef\sphericalangle="0\msx@5E
\mathchardef\varpropto="3\msx@5F
\mathchardef\smallsmile="3\msx@60
\mathchardef\smallfrown="3\msx@61
\mathchardef\Subset="3\msx@62
\mathchardef\Supset="3\msx@63
\mathchardef\Cup="2\msx@64

\mathchardef\Cap="2\msx@65

\mathchardef\curlywedge="2\msx@66
\mathchardef\curlyvee="2\msx@67
\mathchardef\leftthreetimes="2\msx@68
\mathchardef\rightthreetimes="2\msx@69
\mathchardef\subseteqq="3\msx@6A
\mathchardef\supseteqq="3\msx@6B
\mathchardef\bumpeq="3\msx@6C
\mathchardef\Bumpeq="3\msx@6D
\mathchardef\lll="3\msx@6E

\mathchardef\ggg="3\msx@6F

\mathchardef\circledS="0\msx@73
\mathchardef\pitchfork="3\msx@74
\mathchardef\dotplus="2\msx@75
\mathchardef\backsim="3\msx@76
\mathchardef\backsimeq="3\msx@77
\mathchardef\complement="0\msx@7B
\mathchardef\intercal="2\msx@7C
\mathchardef\circledcirc="2\msx@7D
\mathchardef\circledast="2\msx@7E
\mathchardef\circleddash="2\msx@7F
\def\ulcorner{\delimiter"4\msx@70\msx@70 }
\def\urcorner{\delimiter"5\msx@71\msx@71 }
\def\llcorner{\delimiter"4\msx@78\msx@78 }
\def\lrcorner{\delimiter"5\msx@79\msx@79 }
\def\yen{\mathhexbox\msx@55 }
\def\checkmark{\mathhexbox\msx@58 }
\def\circledR{\mathhexbox\msx@72 }
\def\maltese{\mathhexbox\msx@7A }
\mathchardef\lvertneqq="3\msy@00
\mathchardef\gvertneqq="3\msy@01
\mathchardef\nleq="3\msy@02
\mathchardef\ngeq="3\msy@03
\mathchardef\nless="3\msy@04
\mathchardef\ngtr="3\msy@05
\mathchardef\nprec="3\msy@06
\mathchardef\nsucc="3\msy@07
\mathchardef\lneqq="3\msy@08
\mathchardef\gneqq="3\msy@09
\mathchardef\nleqslant="3\msy@0A
\mathchardef\ngeqslant="3\msy@0B
\mathchardef\lneq="3\msy@0C
\mathchardef\gneq="3\msy@0D
\mathchardef\npreceq="3\msy@0E
\mathchardef\nsucceq="3\msy@0F
\mathchardef\precnsim="3\msy@10
\mathchardef\succnsim="3\msy@11
\mathchardef\lnsim="3\msy@12
\mathchardef\gnsim="3\msy@13
\mathchardef\nleqq="3\msy@14
\mathchardef\ngeqq="3\msy@15
\mathchardef\precneqq="3\msy@16
\mathchardef\succneqq="3\msy@17
\mathchardef\precnapprox="3\msy@18
\mathchardef\succnapprox="3\msy@19
\mathchardef\lnapprox="3\msy@1A
\mathchardef\gnapprox="3\msy@1B
\mathchardef\nsim="3\msy@1C
\mathchardef\napprox="3\msy@1D
\mathchardef\ncong="3\msy@1D
\def\napprox{\not\approx}
\mathchardef\varsubsetneq="3\msy@20
\mathchardef\varsupsetneq="3\msy@21
\mathchardef\nsubseteqq="3\msy@22
\mathchardef\nsupseteqq="3\msy@23
\mathchardef\subsetneqq="3\msy@24
\mathchardef\supsetneqq="3\msy@25
\mathchardef\varsubsetneqq="3\msy@26
\mathchardef\varsupsetneqq="3\msy@27
\mathchardef\subsetneq="3\msy@28
\mathchardef\supsetneq="3\msy@29
\mathchardef\nsubseteq="3\msy@2A
\mathchardef\nsupseteq="3\msy@2B
\mathchardef\nparallel="3\msy@2C
\mathchardef\nmid="3\msy@2D
\mathchardef\nshortmid="3\msy@2E
\mathchardef\nshortparallel="3\msy@2F
\mathchardef\nvdash="3\msy@30
\mathchardef\nVdash="3\msy@31
\mathchardef\nvDash="3\msy@32
\mathchardef\nVDash="3\msy@33
\mathchardef\ntrianglerighteq="3\msy@34
\mathchardef\ntrianglelefteq="3\msy@35
\mathchardef\ntriangleleft="3\msy@36
\mathchardef\ntriangleright="3\msy@37
\mathchardef\nleftarrow="3\msy@38
\mathchardef\nrightarrow="3\msy@39
\mathchardef\nLeftarrow="3\msy@3A
\mathchardef\nRightarrow="3\msy@3B
\mathchardef\nLeftrightarrow="3\msy@3C
\mathchardef\nleftrightarrow="3\msy@3D
\mathchardef\divideontimes="2\msy@3E
\mathchardef\varnothing="0\msy@3F
\mathchardef\nexists="0\msy@40
\mathchardef\mho="0\msy@66
\mathchardef\eth="0\msy@67
\mathchardef\eqsim="3\msy@68
\mathchardef\beth="0\msy@69
\mathchardef\gimel="0\msy@6A
\mathchardef\daleth="0\msy@6B
\mathchardef\lessdot="3\msy@6C
\mathchardef\gtrdot="3\msy@6D
\mathchardef\ltimes="2\msy@6E
\mathchardef\rtimes="2\msy@6F
\mathchardef\shortmid="3\msy@70
\mathchardef\shortparallel="3\msy@71
\mathchardef\smallsetminus="2\msy@72
\mathchardef\thicksim="3\msy@73
\mathchardef\thickapprox="3\msy@74
\mathchardef\approxeq="3\msy@75
\mathchardef\succapprox="3\msy@76
\mathchardef\precapprox="3\msy@77
\mathchardef\curvearrowleft="3\msy@78
\mathchardef\curvearrowright="3\msy@79
\mathchardef\digamma="0\msy@7A
\mathchardef\varkappa="0\msy@7B
\mathchardef\hslash="0\msy@7D
\mathchardef\hbar="0\msy@7E
\mathchardef\backepsilon="3\msy@7F

\def\Bbb{\ifmmode\let\next\Bbb@\else
\def\next{\errmessage{Use \string\Bbb\space only in math mode}}\fi\next}
\def\Bbb@#1{{\Bbb@@{#1}}}
\def\Bbb@@#1{\fam\msyfam#1}

\def\squaro{\mathchoice\sqr35\sqr35\sqr{2.1}5\sqr{1.5}5}

\catcode`\@=11
\ifx\amstexloaded@\relax\catcode`\@=\active
  \endinput\else\let\amstexloaded@\relax\fi
\newlinechar=`\^^J
\def\W@{\immediate\write\sixt@@n}

\toksdef\toks@@=2
\long\def\rightappend@#1\to#2{\toks@{\\{#1}}\toks@@
 =\expandafter{#2}\xdef#2{\the\toks@@\the\toks@}\toks@{}\toks@@{}}
\def\alloclist@{}
\newif\ifalloc@
\def\showallocations{{\def\\{\immediate\write\m@ne}\alloclist@}\alloc@true}
\def\alloc@#1#2#3#4#5{\global\advance\count1#1by\@ne
 \ch@ck#1#4#2\allocationnumber=\count1#1
 \global#3#5=\allocationnumber
 \edef\next@{\string#5=\string#2\the\allocationnumber}
 \expandafter\rightappend@\next@\to\alloclist@}
\newcount\count@@
\newcount\count@@@
\def\FN@{\futurelet\next}
\def\DN@{\def\next@}
\def\DNii@{\def\nextii@}
\def\RIfM@{\relax\ifmmode}
\def\RIfMIfI@{\relax\ifmmode\ifinner}
\def\setboxz@h{\setbox\z@\hbox}
\def\wdz@{\wd\z@}
\def\boxz@{\box\z@}
\def\setbox@ne{\setbox\@ne}
\def\wd@ne{\wd\@ne}
\def\iterate{\body\expandafter\iterate\else\fi}
\def\err@#1{\errmessage{AmS-TeX error: #1}}
\newhelp\defaulthelp@{Sorry, I already gave what help I could...^^J
Maybe you should try asking a human?^^J
An error might have occurred before I noticed any problems.^^J
``If all else fails, read the instructions.''}
\def\Err@{\errhelp\defaulthelp@\err@}
\def\eat@#1{}
\def\in@#1#2{\def\in@@##1#1##2##3\in@@{\ifx\in@##2\in@false\else\in@true\fi}
 \in@@#2#1\in@\in@@}
\newif\ifin@
\def\space@.{\futurelet\space@\relax}
\space@. 
\newhelp\athelp@
{Only certain combinations beginning with @ make sense to me.^^J
Perhaps you wanted \string\@\space for a printed @?^^J
I've ignored the character or group after @.}
{\catcode`\~=\active  \lccode`\~=`\@ \lowercase{\gdef~{\FN@\at@}}}
\def\at@{\let\next@\at@@
 \ifcat\noexpand\next a\else\ifcat\noexpand\next0\else
 \ifcat\noexpand\next\relax\else
   \let\next\at@@@\fi\fi\fi
 \next@}
\def\at@@#1{\expandafter
 \ifx\csname\space @\string#1\endcsname\relax
  \expandafter\at@@@ \else
  \csname\space @\string#1\expandafter\endcsname\fi}
\def\at@@@#1{\errhelp\athelp@ \err@{\Invalid@@ @}}
\def\atdef@#1{\expandafter\def\csname\space @\string#1\endcsname}
\newhelp\defahelp@{If you typed \string\define\space cs instead of
\string\define\string\cs\space^^J
I've substituted an inaccessible control sequence so that your^^J
definition will be completed without mixing me up too badly.^^J
If you typed \string\define{\string\cs} the inaccessible control sequence^^J
was defined to be \string\cs, and the rest of your^^J
definition appears as input.}
\newhelp\defbhelp@{I've ignored your definition, because it might^^J
conflict with other uses that are important to me.}
\def\define{\FN@\define@}
\def\define@{\ifcat\noexpand\next\relax
 \expandafter\define@@\else\errhelp\defahelp@                               %1
 \err@{\string\define\space must be followed by a control
 sequence}\expandafter\def\expandafter\nextii@\fi}                          %2
\def\undefined@@@@@@@@@@{}
\def\preloaded@@@@@@@@@@{}
\def\next@@@@@@@@@@{}
\def\define@@#1{\ifx#1\relax\errhelp\defbhelp@                              %1
 \err@{\string#1\space is already defined}\DN@{\DNii@}\else
 \expandafter\ifx\csname\expandafter\eat@\string                            %2
 #1@@@@@@@@@@\endcsname\undefined@@@@@@@@@@\errhelp\defbhelp@
 \err@{\string#1\space can't be defined}\DN@{\DNii@}\else
 \expandafter\ifx\csname\expandafter\eat@\string#1\endcsname\relax          %3
 \global\let#1\undefined\DN@{\def#1}\else\errhelp\defbhelp@
 \err@{\string#1\space is already defined}\DN@{\DNii@}\fi
 \fi\fi\next@}

\def\predefine#1#2{\let#1#2}
\def\undefine#1{\let#1\undefined}
\message{page layout,}
\newdimen\captionwidth@
\captionwidth@\hsize
\advance\captionwidth@-1.5in
\def\pagewidth#1{\hsize#1\relax
 \captionwidth@\hsize\advance\captionwidth@-1.5in}

\def\hcorrection#1{\advance\hoffset#1\relax}
\def\vcorrection#1{\advance\voffset#1\relax}
\message{accents/punctuation,}

\let\graveaccent\`
\let\acuteaccent\'
\let\tildeaccent\~
\let\hataccent\^
\let\underscore\_
\let\B\=
\let\D\.
\let\ic@\/
\def\/{\unskip\ic@}
\def\textfonti{\the\textfont\@ne}
\def\t#1#2{{\edef\next@{\the\font}\textfonti\accent"7F \next@#1#2}}
\def~{\unskip\nobreak\ \ignorespaces}
\def\.{.\spacefactor\@m}
\atdef@;{\leavevmode\null;}
\atdef@:{\leavevmode\null:}
\atdef@?{\leavevmode\null?}
\edef\@{\string @}
\def\relaxnext@{\let\next\relax}
\atdef@-{\relaxnext@\leavevmode
 \DN@{\ifx\next-\DN@-{\FN@\nextii@}\else
  \DN@{\leavevmode\hbox{-}}\fi\next@}%
 \DNii@{\ifx\next-\DN@-{\leavevmode\hbox{---}}\else
  \DN@{\leavevmode\hbox{--}}\fi\next@}%
 \FN@\next@}
\def\srdr@{\kern.16667em}
\def\drsr@{\kern.02778em}
\def\sldl@{\drsr@}
\def\dlsl@{\srdr@}
\atdef@"{\unskip\relaxnext@
 \DN@{\ifx\next\space@\DN@. {\FN@\nextii@}\else
  \DN@.{\FN@\nextii@}\fi\next@.}%
 \DNii@{\ifx\next`\DN@`{\FN@\nextiii@}\else
  \ifx\next\lq\DN@\lq{\FN@\nextiii@}\else
  \DN@####1{\FN@\nextiv@}\fi\fi\next@}%
 \def\nextiii@{\ifx\next`\DN@`{\sldl@``}\else\ifx\next\lq
  \DN@\lq{\sldl@``}\else\DN@{\dlsl@`}\fi\fi\next@}%
 \def\nextiv@{\ifx\next'\DN@'{\srdr@''}\else
  \ifx\next\rq\DN@\rq{\srdr@''}\else\DN@{\drsr@'}\fi\fi\next@}%
 \FN@\next@}

\def\textfontii{\the\textfont\tw@}
\def\lbrace@{\delimiter"4266308 }
\def\rbrace@{\delimiter"5267309 }
\def\{{\RIfM@\lbrace@\else{\textfontii f}\spacefactor\@m\fi}
\def\}{\RIfM@\rbrace@\else
 \let\@sf\empty\ifhmode\edef\@sf{\spacefactor\the\spacefactor}\fi
 {\textfontii g}\@sf\relax\fi}
\let\lbrace\{
\let\rbrace\}
\def\AmSTeX{{\textfontii A\kern-.1667em%
  \lower.5ex\hbox{M}\kern-.125emS}-\TeX}
\message{line and page breaks,}
\def\vmodeerr@#1{\Err@{\string#1\space not allowed between paragraphs}}
\def\mathmodeerr@#1{\Err@{\string#1\space not allowed in math mode}}
\def\linebreak{\RIfM@\mathmodeerr@\linebreak\else
 \ifhmode\unskip\unkern\break\else\vmodeerr@\linebreak\fi\fi}

\newskip\saveskip@
\def\allowlinebreak{\RIfM@\mathmodeerr@\allowlinebreak\else
 \ifhmode\saveskip@\lastskip\unskip
 \allowbreak\ifdim\saveskip@>\z@\hskip\saveskip@\fi
 \else\vmodeerr@\allowlinebreak\fi\fi}
\def\nolinebreak{\RIfM@\mathmodeerr@\nolinebreak\else
 \ifhmode\saveskip@\lastskip\unskip
 \nobreak\ifdim\saveskip@>\z@\hskip\saveskip@\fi
 \else\vmodeerr@\nolinebreak\fi\fi}
\def\newline{\relaxnext@
 \DN@{\RIfM@\expandafter\mathmodeerr@\expandafter\newline\else
  \ifhmode\ifx\next\par\else
  \expandafter\unskip\expandafter\null\expandafter\hfill\expandafter\break\fi
  \else
  \expandafter\vmodeerr@\expandafter\newline\fi\fi}%
 \FN@\next@}
\def\dmatherr@#1{\Err@{\string#1\space not allowed in display math mode}}
\def\nondmatherr@#1{\Err@{\string#1\space not allowed in non-display math
 mode}}
\def\onlydmatherr@#1{\Err@{\string#1\space allowed only in display math mode}}
\def\nonmatherr@#1{\Err@{\string#1\space allowed only in math mode}}
\def\mathbreak{\RIfMIfI@\break\else
 \dmatherr@\mathbreak\fi\else\nonmatherr@\mathbreak\fi}
\def\nomathbreak{\RIfMIfI@\nobreak\else
 \dmatherr@\nomathbreak\fi\else\nonmatherr@\nomathbreak\fi}
\def\allowmathbreak{\RIfMIfI@\allowbreak\else
 \dmatherr@\allowmathbreak\fi\else\nonmatherr@\allowmathbreak\fi}
\def\pagebreak{\RIfM@
 \ifinner\nondmatherr@\pagebreak\else\postdisplaypenalty-\@M\fi
 \else\ifvmode\removelastskip\break\else\vadjust{\break}\fi\fi}
\def\nopagebreak{\RIfM@
 \ifinner\nondmatherr@\nopagebreak\else\postdisplaypenalty\@M\fi
 \else\ifvmode\nobreak\else\vadjust{\nobreak}\fi\fi}
\def\nonvmodeerr@#1{\Err@{\string#1\space not allowed within a paragraph
 or in math}}
\def\vnonvmode@#1#2{\relaxnext@\DNii@{\ifx\next\par\DN@{#1}\else
 \DN@{#2}\fi\next@}%
 \ifvmode\DN@{#1}\else
 \DN@{\FN@\nextii@}\fi\next@}
\def\newpage{\vnonvmode@{\vfill\break}{\nonvmodeerr@\newpage}}
\def\smallpagebreak{\vnonvmode@\smallbreak{\nonvmodeerr@\smallpagebreak}}
\def\medpagebreak{\vnonvmode@\medbreak{\nonvmodeerr@\medpagebreak}}
\def\bigpagebreak{\vnonvmode@\bigbreak{\nonvmodeerr@\bigpagebreak}}
\def\NoBlackBoxes{\global\overfullrule\z@}
\def\BlackBoxes{\global\overfullrule5\p@}
\def\Invalid@#1{\def#1{\Err@{\Invalid@@\string#1}}}
\def\Invalid@@{Invalid use of }
\message{figures,}
\Invalid@\caption
\Invalid@\captionwidth
\newdimen\smallcaptionwidth@
\def\topspace{\mid@false\ins@}
\def\midspace{\mid@true\ins@}
\newif\ifmid@
\def\captionfont@{}
\def\ins@#1{\relaxnext@\allowbreak
 \smallcaptionwidth@\captionwidth@\gdef\thespace@{#1}%
 \DN@{\ifx\next\space@\DN@. {\FN@\nextii@}\else
  \DN@.{\FN@\nextii@}\fi\next@.}%
 \DNii@{\ifx\next\caption\DN@\caption{\FN@\nextiii@}%
  \else\let\next@\nextiv@\fi\next@}%
 \def\nextiv@{\vnonvmode@
  {\ifmid@\expandafter\midinsert\else\expandafter\topinsert\fi
   \vbox to\thespace@{}\endinsert}
  {\ifmid@\nonvmodeerr@\midspace\else\nonvmodeerr@\topspace\fi}}%
 \def\nextiii@{\ifx\next\captionwidth\expandafter\nextv@
  \else\expandafter\nextvi@\fi}%
 \def\nextv@\captionwidth##1##2{\smallcaptionwidth@##1\relax\nextvi@{##2}}%
 \def\nextvi@##1{\def\thecaption@{\captionfont@##1}%
  \DN@{\ifx\next\space@\DN@. {\FN@\nextvii@}\else
   \DN@.{\FN@\nextvii@}\fi\next@.}%
  \FN@\next@}%
 \def\nextvii@{\vnonvmode@
  {\ifmid@\expandafter\midinsert\else
  \expandafter\topinsert\fi\vbox to\thespace@{}\nobreak\smallskip
  \setboxz@h{\noindent\ignorespaces\thecaption@\unskip}%
  \ifdim\wdz@>\smallcaptionwidth@\centerline{\vbox{\hsize\smallcaptionwidth@
   \noindent\ignorespaces\thecaption@\unskip}}%
  \else\centerline{\boxz@}\fi\endinsert}
  {\ifmid@\nonvmodeerr@\midspace
  \else\nonvmodeerr@\topspace\fi}}%
 \FN@\next@}
\message{comments,}
\def\newcodes@{\catcode`\\12\catcode`\{12\catcode`\}12\catcode`\#12%
 \catcode`\%12\relax}
\def\oldcodes@{\catcode`\\0\catcode`\{1\catcode`\}2\catcode`\#6%
 \catcode`\%14\relax}
\def\comment{\newcodes@\endlinechar=10 \comment@}
{\lccode`\0=`\\
\lowercase{\gdef\comment@#1^^J{\comment@@#10endcomment\comment@@@}%
\gdef\comment@@#10endcomment{\FN@\comment@@@}%
\gdef\comment@@@#1\comment@@@{\ifx\next\comment@@@\let\next\comment@
 \else\def\next{\oldcodes@\endlinechar=`\^^M\relax}%
 \fi\next}}}
\def\pr@m@s{\ifx'\next\DN@##1{\prim@s}\else\let\next@\egroup\fi\next@}
\def\prime{{\null\prime@\null}}
\mathchardef\prime@="0230
\let\dsize\displaystyle

\let\ssize\scriptstyle

\message{math spacing,}
\def\,{\RIfM@\mskip\thinmuskip\relax\else\kern.16667em\fi}
\def\!{\RIfM@\mskip-\thinmuskip\relax\else\kern-.16667em\fi}
\let\thinspace\,
\let\negthinspace\!
\def\medspace{\RIfM@\mskip\medmuskip\relax\else\kern.222222em\fi}
\def\negmedspace{\RIfM@\mskip-\medmuskip\relax\else\kern-.222222em\fi}
\def\thickspace{\RIfM@\mskip\thickmuskip\relax\else\kern.27777em\fi}
\let\;\thickspace
\def\negthickspace{\RIfM@\mskip-\thickmuskip\relax\else
 \kern-.27777em\fi}
\atdef@,{\RIfM@\mskip.1\thinmuskip\else\leavevmode\null,\fi}
\atdef@!{\RIfM@\mskip-.1\thinmuskip\else\leavevmode\null!\fi}
\atdef@.{\RIfM@&&\else\leavevmode.\spacefactor3000 \fi}
\def\and{\DOTSB\;\mathchar"3026 \;}

\message{fractions,}
\def\frac#1#2{{#1\over#2}}

\newdimen\ex@
\ex@.2326ex
\Invalid@\thickness
\def\thickfrac{\relaxnext@
 \DN@{\ifx\next\thickness\let\next@\nextii@\else
 \DN@{\nextii@\thickness1}\fi\next@}%
 \DNii@\thickness##1##2##3{{##2\above##1\ex@##3}}%
 \FN@\next@}

\def\thickfracwithdelims#1#2{\relaxnext@\def\ldelim@{#1}\def\rdelim@{#2}%
 \DN@{\ifx\next\thickness\let\next@\nextii@\else
 \DN@{\nextii@\thickness1}\fi\next@}%
 \DNii@\thickness##1##2##3{{##2\abovewithdelims
 \ldelim@\rdelim@##1\ex@##3}}%
 \FN@\next@}

\def\:{\nobreak\hskip.1111em\mathpunct{}\nonscript\mkern-\thinmuskip{:}\hskip
 .3333emplus.0555em\relax}
\def\snug{\unskip\kern-\mathsurround}
\message{smash commands,}
\def\topsmash{\top@true\bot@false\smash@}
\def\botsmash{\top@false\bot@true\smash@}
\newif\iftop@
\newif\ifbot@
\def\smash{\top@true\bot@true\smash@}
\def\smash@{\RIfM@\expandafter\mathpalette\expandafter\mathsm@sh\else
 \expandafter\makesm@sh\fi}
\def\finsm@sh{\iftop@\ht\z@\z@\fi\ifbot@\dp\z@\z@\fi\leavevmode\boxz@}
\message{large operator symbols,}
\def\LimitsOnSums{\global\let\slimits@\displaylimits}
\def\NoLimitsOnSums{\global\let\slimits@\nolimits}
\LimitsOnSums
\mathchardef\coprod@="1360       \def\coprod{\DOTSB\coprod@\slimits@}
\mathchardef\bigvee@="1357       \def\bigvee{\DOTSB\bigvee@\slimits@}
\mathchardef\bigwedge@="1356     \def\bigwedge{\DOTSB\bigwedge@\slimits@}
\mathchardef\biguplus@="1355     \def\biguplus{\DOTSB\biguplus@\slimits@}
\mathchardef\bigcap@="1354       \def\bigcap{\DOTSB\bigcap@\slimits@}
\mathchardef\bigcup@="1353       \def\bigcup{\DOTSB\bigcup@\slimits@}
\mathchardef\prod@="1351         \def\prod{\DOTSB\prod@\slimits@}
\mathchardef\sum@="1350          \def\sum{\DOTSB\sum@\slimits@}
\mathchardef\bigotimes@="134E    \def\bigotimes{\DOTSB\bigotimes@\slimits@}
\mathchardef\bigoplus@="134C     \def\bigoplus{\DOTSB\bigoplus@\slimits@}
\mathchardef\bigodot@="134A      \def\bigodot{\DOTSB\bigodot@\slimits@}
\mathchardef\bigsqcup@="1346     \def\bigsqcup{\DOTSB\bigsqcup@\slimits@}
\message{integrals,}
\def\LimitsOnInts{\global\let\ilimits@\displaylimits}
\def\NoLimitsOnInts{\global\let\ilimits@\nolimits}
\NoLimitsOnInts
\def\int{\DOTSI\intop\ilimits@}
\def\oint{\DOTSI\ointop\ilimits@}
\def\intic@{\mathchoice{\hskip.5em}{\hskip.4em}{\hskip.4em}{\hskip.4em}}
\def\negintic@{\mathchoice
 {\hskip-.5em}{\hskip-.4em}{\hskip-.4em}{\hskip-.4em}}
\def\intkern@{\mathchoice{\!\!\!}{\!\!}{\!\!}{\!\!}}
\def\intdots@{\mathchoice{\plaincdots@}
 {{\cdotp}\mkern1.5mu{\cdotp}\mkern1.5mu{\cdotp}}
 {{\cdotp}\mkern1mu{\cdotp}\mkern1mu{\cdotp}}
 {{\cdotp}\mkern1mu{\cdotp}\mkern1mu{\cdotp}}}
\newcount\intno@
\def\iint{\DOTSI\intno@\tw@\FN@\ints@}
\def\iiint{\DOTSI\intno@\thr@@\FN@\ints@}
\def\iiiint{\DOTSI\intno@4 \FN@\ints@}
\def\idotsint{\DOTSI\intno@\z@\FN@\ints@}
\def\ints@{\findlimits@\ints@@}
\newif\iflimtoken@
\newif\iflimits@
\def\findlimits@{\limtoken@true\ifx\next\limits\limits@true
 \else\ifx\next\nolimits\limits@false\else
 \limtoken@false\ifx\ilimits@\nolimits\limits@false\else
 \ifinner\limits@false\else\limits@true\fi\fi\fi\fi}
\def\multint@{\int\ifnum\intno@=\z@\intdots@                                %1
 \else\intkern@\fi                                                          %2
 \ifnum\intno@>\tw@\int\intkern@\fi                                         %3
 \ifnum\intno@>\thr@@\int\intkern@\fi                                       %4
 \int}                                                                      %5
\def\multintlimits@{\intop\ifnum\intno@=\z@\intdots@\else\intkern@\fi
 \ifnum\intno@>\tw@\intop\intkern@\fi
 \ifnum\intno@>\thr@@\intop\intkern@\fi\intop}
\def\ints@@{\iflimtoken@                                                    %1
 \def\ints@@@{\iflimits@\negintic@\mathop{\intic@\multintlimits@}\limits    %2
  \else\multint@\nolimits\fi                                                %3
  \eat@}                                                                    %4
 \else                                                                      %5
 \def\ints@@@{\iflimits@\negintic@
  \mathop{\intic@\multintlimits@}\limits\else
  \multint@\nolimits\fi}\fi\ints@@@}
\def\LimitsOnNames{\global\let\nlimits@\displaylimits}
\def\NoLimitsOnNames{\global\let\nlimits@\nolimits@}
\LimitsOnNames
\def\nolimits@{\relaxnext@
 \DN@{\ifx\next\limits\DN@\limits{\nolimits}\else
  \let\next@\nolimits\fi\next@}%
 \FN@\next@}
\message{operator names,}
\def\newmcodes@{\mathcode`\'"27\mathcode`\*"2A\mathcode`\."613A%
 \mathcode`\-"2D\mathcode`\/"2F\mathcode`\:"603A }
\def\operatorname#1{\mathop{\newmcodes@\kern\z@\fam\z@#1}\nolimits@}
\def\operatornamewithlimits#1{\mathop{\newmcodes@\kern\z@\fam\z@#1}\nlimits@}
\def\qopname@#1{\mathop{\fam\z@#1}\nolimits@}
\def\qopnamewl@#1{\mathop{\fam\z@#1}\nlimits@}
\def\arccos{\qopname@{arccos}}
\def\arcsin{\qopname@{arcsin}}
\def\arctan{\qopname@{arctan}}
\def\arg{\qopname@{arg}}
\def\cos{\qopname@{cos}}
\def\cosh{\qopname@{cosh}}
\def\cot{\qopname@{cot}}
\def\coth{\qopname@{coth}}
\def\csc{\qopname@{csc}}
\def\deg{\qopname@{deg}}
\def\det{\qopnamewl@{det}}
\def\dim{\qopname@{dim}}
\def\exp{\qopname@{exp}}
\def\gcd{\qopnamewl@{gcd}}
\def\hom{\qopname@{hom}}
\def\inf{\qopnamewl@{inf}}
\def\injlim{\qopnamewl@{inj\,lim}}
\def\ker{\qopname@{ker}}
\def\lg{\qopname@{lg}}
\def\lim{\qopnamewl@{lim}}
\def\liminf{\qopnamewl@{lim\,inf}}
\def\limsup{\qopnamewl@{lim\,sup}}
\def\ln{\qopname@{ln}}
\def\log{\qopname@{log}}
\def\max{\qopnamewl@{max}}
\def\min{\qopnamewl@{min}}
\def\Pr{\qopnamewl@{Pr}}
\def\projlim{\qopnamewl@{proj\,lim}}
\def\sec{\qopname@{sec}}
\def\sin{\qopname@{sin}}
\def\sinh{\qopname@{sinh}}
\def\sup{\qopnamewl@{sup}}
\def\tan{\qopname@{tan}}
\def\tanh{\qopname@{tanh}}
\def\varinjlim{\mathop{\vtop{\ialign{##\crcr
 \hfil\rm lim\hfil\crcr\noalign{\nointerlineskip}\rightarrowfill\crcr
 \noalign{\nointerlineskip\kern-\ex@}\crcr}}}}
\def\varprojlim{\mathop{\vtop{\ialign{##\crcr
 \hfil\rm lim\hfil\crcr\noalign{\nointerlineskip}\leftarrowfill\crcr
 \noalign{\nointerlineskip\kern-\ex@}\crcr}}}}
\def\varliminf{\mathop{\underline{\vrule height\z@ depth.2exwidth\z@
 \hbox{\rm lim}}}}

\newdimen\buffer@
\buffer@\fontdimen13 \tenex
\newdimen\buffer
\buffer\buffer@

\def\ResetBuffer{\fontdimen13 \tenex\buffer@\global\buffer\buffer@}
\def\shave#1{\mathop{\hbox{$\m@th\fontdimen13 \tenex\z@                     %1
 \displaystyle{#1}$}}\fontdimen13 \tenex\buffer}

\message{multilevel sub/superscripts,}
\Invalid@\\
\def\Let@{\relax\iffalse{\fi\let\\=\cr\iffalse}\fi}
\Invalid@\vspace
\def\vspace@{\def\vspace##1{\crcr\noalign{\vskip##1\relax}}}
\def\multilimits@{\bgroup\vspace@\Let@
 \baselineskip\fontdimen10 \scriptfont\tw@
 \advance\baselineskip\fontdimen12 \scriptfont\tw@
 \lineskip\thr@@\fontdimen8 \scriptfont\thr@@
 \lineskiplimit\lineskip
 \vbox\bgroup\ialign\bgroup\hfil$\m@th\scriptstyle{##}$\hfil\crcr}
\def\Sb{_\multilimits@}
\def\endSb{\crcr\egroup\egroup\egroup}
\def\Sp{^\multilimits@}

\def\spreadlines#1{\RIfMIfI@\onlydmatherr@\spreadlines\else
 \openup#1\relax\fi\else\onlydmatherr@\spreadlines\fi}
\def\Mathstrut@{\copy\Mathstrutbox@}
\newbox\Mathstrutbox@
\setbox\Mathstrutbox@\null
\setboxz@h{$\m@th($}
\ht\Mathstrutbox@\ht\z@
\dp\Mathstrutbox@\dp\z@
\message{matrices,}
\newdimen\spreadmlines@
\def\spreadmatrixlines#1{\RIfMIfI@
 \onlydmatherr@\spreadmatrixlines\else
 \spreadmlines@#1\relax\fi\else\onlydmatherr@\spreadmatrixlines\fi}
\def\matrix{\null\,\vcenter\bgroup\Let@\vspace@
 \normalbaselines\openup\spreadmlines@\ialign
 \bgroup\hfil$\m@th##$\hfil&&\quad\hfil$\m@th##$\hfil\crcr
 \Mathstrut@\crcr\noalign{\kern-\baselineskip}}
\def\endmatrix{\crcr\Mathstrut@\crcr\noalign{\kern-\baselineskip}\egroup
 \egroup\,}
\def\format{\crcr\egroup\iffalse{\fi\ifnum`}=0 \fi\format@}
\newtoks\hashtoks@
\hashtoks@{#}
\def\format@#1\\{\def\preamble@{#1}%
 \def\l{$\m@th\the\hashtoks@$\hfil}%
 \def\c{\hfil$\m@th\the\hashtoks@$\hfil}%
 \def\r{\hfil$\m@th\the\hashtoks@$}%
 \edef\preamble@@{\preamble@}\ifnum`{=0 \fi\iffalse}\fi
 \ialign\bgroup\span\preamble@@\crcr}
\def\smallmatrix{\null\,\vcenter\bgroup\vspace@\Let@
 \baselineskip9\ex@\lineskip\ex@
 \ialign\bgroup\hfil$\m@th\scriptstyle{##}$\hfil&&\thickspace\hfil
 $\m@th\scriptstyle{##}$\hfil\crcr}
\def\endsmallmatrix{\crcr\egroup\egroup\,}

\newmuskip\dotsspace@
\dotsspace@1.5mu
\def\strip@#1 {#1}
\def\spacehdots#1\for#2{\multispan{#2}\xleaders
 \hbox{$\m@th\mkern\strip@#1 \dotsspace@.\mkern\strip@#1 \dotsspace@$}\hfill}
\def\hdotsfor#1{\spacehdots\@ne\for{#1}}
\def\multispan@#1{\omit\mscount#1\unskip\loop\ifnum\mscount>\@ne\sp@n\repeat}
\def\spaceinnerhdots#1\for#2\after#3{\multispan@{\strip@#2 }#3\xleaders
 \hbox{$\m@th\mkern\strip@#1 \dotsspace@.\mkern\strip@#1 \dotsspace@$}\hfill}
\def\innerhdotsfor#1\after#2{\spaceinnerhdots\@ne\for#1\after{#2}}
\def\cases{\bgroup\spreadmlines@\jot\left\{\,\matrix\format\l&\quad\l\\}
\def\endcases{\endmatrix\right.\egroup}
\message{multiline displays,}
\newif\ifinany@
\newif\ifinalign@
\newif\ifingather@
\def\strut@{\copy\strutbox@}
\newbox\strutbox@
\setbox\strutbox@\hbox{\vrule height8\p@ depth3\p@ width\z@}
\def\topaligned{\null\,\vtop\aligned@}
\def\botaligned{\null\,\vbox\aligned@}
\def\aligned{\null\,\vcenter\aligned@}
\def\aligned@{\bgroup\vspace@\Let@
 \ifinany@\else\openup\jot\fi\ialign
 \bgroup\hfil\strut@$\m@th\displaystyle{##}$&
 $\m@th\displaystyle{{}##}$\hfil\crcr}
\def\endaligned{\crcr\egroup\egroup}

\def\alignedat#1{\null\,\vcenter\bgroup\doat@{#1}\vspace@\Let@
 \ifinany@\else\openup\jot\fi\ialign\bgroup\span\preamble@@\crcr}
\newcount\atcount@
\def\doat@#1{\toks@{\hfil\strut@$\m@th
 \displaystyle{\the\hashtoks@}$&$\m@th\displaystyle
 {{}\the\hashtoks@}$\hfil}%                                                 %1
 \atcount@#1\relax\advance\atcount@\m@ne                                    %2
 \loop\ifnum\atcount@>\z@\toks@=\expandafter{\the\toks@&\hfil$\m@th
 \displaystyle{\the\hashtoks@}$&$\m@th
 \displaystyle{{}\the\hashtoks@}$\hfil}\advance
  \atcount@\m@ne\repeat                                                     %3
 \xdef\preamble@{\the\toks@}\xdef\preamble@@{\preamble@}}

\def\gathered{\null\,\vcenter\bgroup\vspace@\Let@
 \ifinany@\else\openup\jot\fi\ialign
 \bgroup\hfil\strut@$\m@th\displaystyle{##}$\hfil\crcr}
\def\endgathered{\crcr\egroup\egroup}
\newif\iftagsleft@
\def\TagsOnLeft{\global\tagsleft@true}
\def\TagsOnRight{\global\tagsleft@false}
\TagsOnLeft
\newif\ifmathtags@
\def\TagsAsMath{\global\mathtags@true}
\def\TagsAsText{\global\mathtags@false}
\TagsAsText
\def\tagform@#1{\hbox{\rm(\ignorespaces#1\unskip)}}
\def\thetag{\leavevmode\tagform@}
\def\tag#1$${\iftagsleft@\leqno\else\eqno\fi                                %1
 \maketag@#1\maketag@                                                       %2
 $$}                                                                        %3
\def\maketag@{\FN@\maketag@@}
\def\maketag@@{\ifx\next"\expandafter\maketag@@@\else\expandafter\maketag@@@@
 \fi}
\def\maketag@@@"#1"#2\maketag@{\hbox{\rm#1}}                                %1
\def\maketag@@@@#1\maketag@{\ifmathtags@\tagform@{$\m@th#1$}\else
 \tagform@{#1}\fi}
\interdisplaylinepenalty\@M
\def\allowdisplaybreaks{\RIfMIfI@
 \onlydmatherr@\allowdisplaybreaks\else
 \interdisplaylinepenalty\z@\fi\else\onlydmatherr@\allowdisplaybreaks\fi}
\Invalid@\allowdisplaybreak
\Invalid@\displaybreak
\Invalid@\intertext
\def\allowdisplaybreak@{\def\allowdisplaybreak{\crcr\noalign{\allowbreak}}}
\def\displaybreak@{\def\displaybreak{\crcr\noalign{\break}}}
\def\intertext@{\def\intertext##1{\crcr\noalign{%
 \penalty\postdisplaypenalty \vskip\belowdisplayskip
 \vbox{\normalbaselines\noindent##1}%
 \penalty\predisplaypenalty \vskip\abovedisplayskip}}}
\newskip\centering@
\centering@\z@ plus\@m\p@
\def\align{\relax\ifingather@\DN@{\csname align (in
  \string\gather)\endcsname}\else
 \ifmmode\ifinner\DN@{\onlydmatherr@\align}\else
  \let\next@\align@\fi
 \else\DN@{\onlydmatherr@\align}\fi\fi\next@}
\newhelp\andhelp@
{An extra & here is so disastrous that you should probably exit^^J
and fix things up.}
\newif\iftag@
\newcount\and@
\def\align@{\inalign@true\inany@true
 \vspace@\allowdisplaybreak@\displaybreak@\intertext@
 \def\tag{\global\tag@true\ifnum\and@=\z@\DN@{&&}\else
          \DN@{&}\fi\next@}%
 \iftagsleft@\DN@{\csname align \endcsname}\else
  \DN@{\csname align \space\endcsname}\fi\next@}
\def\Tag@{\iftag@\else\errhelp\andhelp@\err@{Extra & on this line}\fi}
\newdimen\lwidth@
\newdimen\rwidth@
\newdimen\maxlwidth@
\newdimen\maxrwidth@
\newdimen\totwidth@
\def\measure@#1\endalign{\lwidth@\z@\rwidth@\z@\maxlwidth@\z@\maxrwidth@\z@
 \global\and@\z@                                                            %1
 \setbox@ne\vbox                                                            %2
  {\everycr{\noalign{\global\tag@false\global\and@\z@}}\Let@                %3
  \halign{\setboxz@h{$\m@th\displaystyle{\@lign##}$}%                       %4
   \global\lwidth@\wdz@                                                     %5
   \ifdim\lwidth@>\maxlwidth@\global\maxlwidth@\lwidth@\fi                  %6
   \global\advance\and@\@ne                                                 %7
   &\setboxz@h{$\m@th\displaystyle{{}\@lign##}$}\global\rwidth@\wdz@        %8
   \ifdim\rwidth@>\maxrwidth@\global\maxrwidth@\rwidth@\fi                  %9
   \global\advance\and@\@ne                                                %10
   &\Tag@
   \eat@{##}\crcr#1\crcr}}%                                                %11
 \totwidth@\maxlwidth@\advance\totwidth@\maxrwidth@}                       %12
\def\displ@y@{\global\dt@ptrue\openup\jot
 \everycr{\noalign{\global\tag@false\global\and@\z@\ifdt@p\global\dt@pfalse
 \vskip-\lineskiplimit\vskip\normallineskiplimit\else
 \penalty\interdisplaylinepenalty\fi}}}
\def\black@#1{\noalign{\ifdim#1>\displaywidth
 \dimen@\prevdepth\nointerlineskip                                          %1
 \vskip-\ht\strutbox@\vskip-\dp\strutbox@                                   %2
 \vbox{\noindent\hbox to#1{\strut@\hfill}}%                                 %3
 \prevdepth\dimen@                                                          %4
 \fi}}
\expandafter\def\csname align \space\endcsname#1\endalign
 {\measure@#1\endalign\global\and@\z@                                       %1
 \ifingather@\everycr{\noalign{\global\and@\z@}}\else\displ@y@\fi           %2
 \Let@\tabskip\centering@                                                   %3
 \halign to\displaywidth
  {\hfil\strut@\setboxz@h{$\m@th\displaystyle{\@lign##}$}%                  %4
  \global\lwidth@\wdz@\boxz@\global\advance\and@\@ne                        %5
  \tabskip\z@skip                                                           %6
  &\setboxz@h{$\m@th\displaystyle{{}\@lign##}$}%                            %7
  \global\rwidth@\wdz@\boxz@\hfill\global\advance\and@\@ne                  %8
  \tabskip\centering@                                                       %9
  &\setboxz@h{\@lign\strut@\maketag@##\maketag@}%                          %10
  \dimen@\displaywidth\advance\dimen@-\totwidth@
  \divide\dimen@\tw@\advance\dimen@\maxrwidth@\advance\dimen@-\rwidth@     %11
  \ifdim\dimen@<\tw@\wdz@\llap{\vtop{\normalbaselines\null\boxz@}}%        %12
  \else\llap{\boxz@}\fi                                                    %13
  \tabskip\z@skip                                                          %14
  \crcr#1\crcr                                                             %15
  \black@\totwidth@}}                                                      %16
\newdimen\lineht@
\expandafter\def\csname align \endcsname#1\endalign{\measure@#1\endalign
 \global\and@\z@
 \ifdim\totwidth@>\displaywidth\let\displaywidth@\totwidth@\else
  \let\displaywidth@\displaywidth\fi                                        %1
 \ifingather@\everycr{\noalign{\global\and@\z@}}\else\displ@y@\fi
 \Let@\tabskip\centering@\halign to\displaywidth
  {\hfil\strut@\setboxz@h{$\m@th\displaystyle{\@lign##}$}%
  \global\lwidth@\wdz@\global\lineht@\ht\z@                                 %2
  \boxz@\global\advance\and@\@ne
  \tabskip\z@skip&\setboxz@h{$\m@th\displaystyle{{}\@lign##}$}%
  \global\rwidth@\wdz@\ifdim\ht\z@>\lineht@\global\lineht@\ht\z@\fi         %3
  \boxz@\hfil\global\advance\and@\@ne
  \tabskip\centering@&\kern-\displaywidth@                                  %4
  \setboxz@h{\@lign\strut@\maketag@##\maketag@}%
  \dimen@\displaywidth\advance\dimen@-\totwidth@
  \divide\dimen@\tw@\advance\dimen@\maxlwidth@\advance\dimen@-\lwidth@
  \ifdim\dimen@<\tw@\wdz@
   \rlap{\vbox{\normalbaselines\boxz@\vbox to\lineht@{}}}\else
   \rlap{\boxz@}\fi
  \tabskip\displaywidth@\crcr#1\crcr\black@\totwidth@}}
\expandafter\def\csname align (in \string\gather)\endcsname
  #1\endalign{\vcenter{\align@#1\endalign}}
\Invalid@\endalign
\newif\ifxat@
\def\alignat{\RIfMIfI@\DN@{\onlydmatherr@\alignat}\else
 \DN@{\csname alignat \endcsname}\fi\else
 \DN@{\onlydmatherr@\alignat}\fi\next@}
\newif\ifmeasuring@
\newbox\savealignat@
\expandafter\def\csname alignat \endcsname#1#2\endalignat                   %1
 {\inany@true\xat@false
 \def\tag{\global\tag@true\count@#1\relax\multiply\count@\tw@
  \xdef\tag@{}\loop\ifnum\count@>\and@\xdef\tag@{&\tag@}\advance\count@\m@ne
  \repeat\tag@}%
 \vspace@\allowdisplaybreak@\displaybreak@\intertext@
 \displ@y@\measuring@true                                                   %2
 \setbox\savealignat@\hbox{$\m@th\displaystyle\Let@
  \attag@{#1}%                                                              %3
  \vbox{\halign{\span\preamble@@\crcr#2\crcr}}$}%
 \measuring@false                                                           %4
 \Let@\attag@{#1}%                                                          %5
 \tabskip\centering@\halign to\displaywidth
  {\span\preamble@@\crcr#2\crcr                                             %6
  \black@{\wd\savealignat@}}}                                               %7
\Invalid@\endalignat
\def\xalignat{\RIfMIfI@
 \DN@{\onlydmatherr@\xalignat}\else
 \DN@{\csname xalignat \endcsname}\fi\else
 \DN@{\onlydmatherr@\xalignat}\fi\next@}
\expandafter\def\csname xalignat \endcsname#1#2\endxalignat
 {\inany@true\xat@true
 \def\tag{\global\tag@true\def\tag@{}\count@#1\relax\multiply\count@\tw@
  \loop\ifnum\count@>\and@\xdef\tag@{&\tag@}\advance\count@\m@ne\repeat\tag@}%
 \vspace@\allowdisplaybreak@\displaybreak@\intertext@
 \displ@y@\measuring@true\setbox\savealignat@\hbox{$\m@th\displaystyle\Let@
 \attag@{#1}\vbox{\halign{\span\preamble@@\crcr#2\crcr}}$}%
 \measuring@false\Let@
 \attag@{#1}\tabskip\centering@\halign to\displaywidth
 {\span\preamble@@\crcr#2\crcr\black@{\wd\savealignat@}}}
\def\attag@#1{\let\Maketag@\maketag@\let\TAG@\Tag@                          %1
 \let\Tag@=0\let\maketag@=0%                                                %2
 \ifmeasuring@\def\llap@##1{\setboxz@h{##1}\hbox to\tw@\wdz@{}}%
  \def\rlap@##1{\setboxz@h{##1}\hbox to\tw@\wdz@{}}\else
  \let\llap@\llap\let\rlap@\rlap\fi                                         %3
 \toks@{\hfil\strut@$\m@th\displaystyle{\@lign\the\hashtoks@}$\tabskip\z@skip
  \global\advance\and@\@ne&$\m@th\displaystyle{{}\@lign\the\hashtoks@}$\hfil
  \ifxat@\tabskip\centering@\fi\global\advance\and@\@ne}%                   %4
 \iftagsleft@
  \toks@@{\tabskip\centering@&\Tag@\kern-\displaywidth
   \rlap@{\@lign\maketag@\the\hashtoks@\maketag@}%
   \global\advance\and@\@ne\tabskip\displaywidth}\else
  \toks@@{\tabskip\centering@&\Tag@\llap@{\@lign\maketag@
   \the\hashtoks@\maketag@}\global\advance\and@\@ne\tabskip\z@skip}\fi      %5
 \atcount@#1\relax\advance\atcount@\m@ne
 \loop\ifnum\atcount@>\z@
 \toks@=\expandafter{\the\toks@&\hfil$\m@th\displaystyle{\@lign
  \the\hashtoks@}$\global\advance\and@\@ne
  \tabskip\z@skip&$\m@th\displaystyle{{}\@lign\the\hashtoks@}$\hfil\ifxat@
  \tabskip\centering@\fi\global\advance\and@\@ne}\advance\atcount@\m@ne
 \repeat                                                                    %6
 \xdef\preamble@{\the\toks@\the\toks@@}%                                    %7
 \xdef\preamble@@{\preamble@}%                                              %8
 \let\maketag@\Maketag@\let\Tag@\TAG@}                                      %9
\Invalid@\endxalignat
\def\xxalignat{\RIfMIfI@
 \DN@{\onlydmatherr@\xxalignat}\else\DN@{\csname xxalignat
  \endcsname}\fi\else
 \DN@{\onlydmatherr@\xxalignat}\fi\next@}
\expandafter\def\csname xxalignat \endcsname#1#2\endxxalignat{\inany@true
 \vspace@\allowdisplaybreak@\displaybreak@\intertext@
 \displ@y\setbox\savealignat@\hbox{$\m@th\displaystyle\Let@
 \xxattag@{#1}\vbox{\halign{\span\preamble@@\crcr#2\crcr}}$}%
 \Let@\xxattag@{#1}\tabskip\z@skip\halign to\displaywidth
 {\span\preamble@@\crcr#2\crcr\black@{\wd\savealignat@}}}
\def\xxattag@#1{\toks@{\tabskip\z@skip\hfil\strut@
 $\m@th\displaystyle{\the\hashtoks@}$&%
 $\m@th\displaystyle{{}\the\hashtoks@}$\hfil\tabskip\centering@&}%
 \atcount@#1\relax\advance\atcount@\m@ne\loop\ifnum\atcount@>\z@
 \toks@=\expandafter{\the\toks@&\hfil$\m@th\displaystyle{\the\hashtoks@}$%
  \tabskip\z@skip&$\m@th\displaystyle{{}\the\hashtoks@}$\hfil
  \tabskip\centering@}\advance\atcount@\m@ne\repeat
 \xdef\preamble@{\the\toks@\tabskip\z@skip}\xdef\preamble@@{\preamble@}}
\Invalid@\endxxalignat
\newdimen\gwidth@
\newdimen\gmaxwidth@
\def\gmeasure@#1\endgather{\gwidth@\z@\gmaxwidth@\z@\setbox@ne\vbox{\Let@
 \halign{\setboxz@h{$\m@th\displaystyle{##}$}\global\gwidth@\wdz@
 \ifdim\gwidth@>\gmaxwidth@\global\gmaxwidth@\gwidth@\fi
 &\eat@{##}\crcr#1\crcr}}}
\def\gather{\RIfMIfI@\DN@{\onlydmatherr@\gather}\else
 \ingather@true\inany@true\def\tag{&}%
 \vspace@\allowdisplaybreak@\displaybreak@\intertext@
 \displ@y\Let@
 \iftagsleft@\DN@{\csname gather \endcsname}\else
  \DN@{\csname gather \space\endcsname}\fi\fi
 \else\DN@{\onlydmatherr@\gather}\fi\next@}
\expandafter\def\csname gather \space\endcsname#1\endgather
 {\gmeasure@#1\endgather\tabskip\centering@
 \halign to\displaywidth{\hfil\strut@\setboxz@h{$\m@th\displaystyle{##}$}%
 \global\gwidth@\wdz@\boxz@\hfil&
 \setboxz@h{\strut@{\maketag@##\maketag@}}%
 \dimen@\displaywidth\advance\dimen@-\gwidth@
 \ifdim\dimen@>\tw@\wdz@\llap{\boxz@}\else
 \llap{\vtop{\normalbaselines\null\boxz@}}\fi
 \tabskip\z@skip\crcr#1\crcr\black@\gmaxwidth@}}
\newdimen\glineht@
\expandafter\def\csname gather \endcsname#1\endgather{\gmeasure@#1\endgather
 \ifdim\gmaxwidth@>\displaywidth\let\gdisplaywidth@\gmaxwidth@\else
 \let\gdisplaywidth@\displaywidth\fi\tabskip\centering@\halign to\displaywidth
 {\hfil\strut@\setboxz@h{$\m@th\displaystyle{##}$}%
 \global\gwidth@\wdz@\global\glineht@\ht\z@\boxz@\hfil&\kern-\gdisplaywidth@
 \setboxz@h{\strut@{\maketag@##\maketag@}}%
 \dimen@\displaywidth\advance\dimen@-\gwidth@
 \ifdim\dimen@>\tw@\wdz@\rlap{\boxz@}\else
 \rlap{\vbox{\normalbaselines\boxz@\vbox to\glineht@{}}}\fi
 \tabskip\gdisplaywidth@\crcr#1\crcr\black@\gmaxwidth@}}
\newif\ifctagsplit@
\def\CenteredTagsOnSplits{\global\ctagsplit@true}
\def\TopOrBottomTagsOnSplits{\global\ctagsplit@false}
\TopOrBottomTagsOnSplits
\def\split{\relax\ifinany@\let\next@\insplit@\else
 \ifmmode\ifinner\def\next@{\onlydmatherr@\split}\else
 \let\next@\outsplit@\fi\else
 \def\next@{\onlydmatherr@\split}\fi\fi\next@}
\def\insplit@{\global\setbox\z@\vbox\bgroup\vspace@\Let@\ialign\bgroup
 \hfil\strut@$\m@th\displaystyle{##}$&$\m@th\displaystyle{{}##}$\hfill\crcr}
\def\endsplit{\crcr\egroup\egroup\iftagsleft@\expandafter\lendsplit@\else
 \expandafter\rendsplit@\fi}
\def\rendsplit@{\global\setbox9 \vbox
 {\unvcopy\z@\global\setbox8 \lastbox\unskip}%                              %1
 \setbox@ne\hbox{\unhcopy8 \unskip\global\setbox\tw@\lastbox
 \unskip\global\setbox\thr@@\lastbox}%                                      %2
 \global\setbox7 \hbox{\unhbox\tw@\unskip}%                                 %3
 \ifinalign@\ifctagsplit@                                                   %4
  \gdef\split@{\hbox to\wd\thr@@{}&
   \vcenter{\vbox{\moveleft\wd\thr@@\boxz@}}}%                              %5
 \else\gdef\split@{&\vbox{\moveleft\wd\thr@@\box9}\crcr
  \box\thr@@&\box7}\fi                                                      %6
 \else                                                                      %7
  \ifctagsplit@\gdef\split@{\vcenter{\boxz@}}\else
  \gdef\split@{\box9\crcr\hbox{\box\thr@@\box7}}\fi
 \fi
 \split@}                                                                   %8
\def\lendsplit@{\global\setbox9\vtop{\unvcopy\z@}%                          %1
 \setbox@ne\vbox{\unvcopy\z@\global\setbox8\lastbox}%                       %2
 \setbox@ne\hbox{\unhcopy8\unskip\setbox\tw@\lastbox
  \unskip\global\setbox\thr@@\lastbox}%                                     %3
 \ifinalign@\ifctagsplit@                                                   %4
  \gdef\split@{\hbox to\wd\thr@@{}&
  \vcenter{\vbox{\moveleft\wd\thr@@\box9}}}%                                %5
  \else                                                                     %6
  \gdef\split@{\hbox to\wd\thr@@{}&\vbox{\moveleft\wd\thr@@\box9}}\fi
 \else
  \ifctagsplit@\gdef\split@{\vcenter{\box9}}\else
  \gdef\split@{\box9}\fi
 \fi\split@}
\def\outsplit@#1$${\align\insplit@#1\endalign$$}
\newdimen\multlinegap@
\multlinegap@1em
\newdimen\multlinetaggap@
\multlinetaggap@1em
\def\MultlineGap#1{\global\multlinegap@#1\relax}
\def\multlinegap#1{\RIfMIfI@\onlydmatherr@\multlinegap\else
 \multlinegap@#1\relax\fi\else\onlydmatherr@\multlinegap\fi}
\def\nomultlinegap{\multlinegap{\z@}}
\def\multline{\RIfMIfI@
 \DN@{\onlydmatherr@\multline}\else
 \DN@{\multline@}\fi\else
 \DN@{\onlydmatherr@\multline}\fi\next@}
\newif\iftagin@
\def\tagin@#1{\tagin@false\in@\tag{#1}\ifin@\tagin@true\fi}
\def\multline@#1$${\inany@true\vspace@\allowdisplaybreak@\displaybreak@
 \tagin@{#1}\iftagsleft@\DN@{\multline@l#1$$}\else
 \DN@{\multline@r#1$$}\fi\next@}
\newdimen\mwidth@
\def\rmmeasure@#1\endmultline{%
 \def\shoveleft##1{##1}\def\shoveright##1{##1}%                             %1
 \setbox@ne\vbox{\Let@\halign{\setboxz@h
  {$\m@th\@lign\displaystyle{}##$}\global\mwidth@\wdz@
  \crcr#1\crcr}}}
\newdimen\mlineht@
\newif\ifzerocr@
\newif\ifonecr@
\def\lmmeasure@#1\endmultline{\global\zerocr@true\global\onecr@false
 \everycr{\noalign{\ifonecr@\global\onecr@false\fi
  \ifzerocr@\global\zerocr@false\global\onecr@true\fi}}%                    %1
  \def\shoveleft##1{##1}\def\shoveright##1{##1}%
 \setbox@ne\vbox{\Let@\halign{\setboxz@h
  {$\m@th\@lign\displaystyle{}##$}\ifonecr@\global\mwidth@\wdz@
  \global\mlineht@\ht\z@\fi\crcr#1\crcr}}}
\newbox\mtagbox@
\newdimen\ltwidth@
\newdimen\rtwidth@
\def\multline@l#1$${\iftagin@\DN@{\lmultline@@#1$$}\else
 \DN@{\setbox\mtagbox@\null\ltwidth@\z@\rtwidth@\z@
  \lmultline@@@#1$$}\fi\next@}
\def\lmultline@@#1\endmultline\tag#2$${%
 \setbox\mtagbox@\hbox{\maketag@#2\maketag@}%                               %1
 \lmmeasure@#1\endmultline\dimen@\mwidth@\advance\dimen@\wd\mtagbox@
 \advance\dimen@\multlinetaggap@                                            %2
 \ifdim\dimen@>\displaywidth\ltwidth@\z@\else\ltwidth@\wd\mtagbox@\fi       %3
 \lmultline@@@#1\endmultline$$}
\def\lmultline@@@{\displ@y
 \def\shoveright##1{##1\hfilneg\hskip\multlinegap@}%
 \def\shoveleft##1{\setboxz@h{$\m@th\displaystyle{}##1$}%
  \setbox@ne\hbox{$\m@th\displaystyle##1$}%
  \hfilneg
  \iftagin@
   \ifdim\ltwidth@>\z@\hskip\ltwidth@\hskip\multlinetaggap@\fi
  \else\hskip\multlinegap@\fi\hskip.5\wd@ne\hskip-.5\wdz@##1}%              %1
  \halign\bgroup\Let@\hbox to\displaywidth
   {\strut@$\m@th\displaystyle\hfil{}##\hfil$}\crcr
   \hfilneg                                                                 %2
   \iftagin@                                                                %3
    \ifdim\ltwidth@>\z@                                                     %4
     \box\mtagbox@\hskip\multlinetaggap@                                    %5
    \else
     \rlap{\vbox{\normalbaselines\hbox{\strut@\box\mtagbox@}%
     \vbox to\mlineht@{}}}\fi                                               %6
   \else\hskip\multlinegap@\fi}                                             %7
\def\multline@r#1$${\iftagin@\DN@{\rmultline@@#1$$}\else
 \DN@{\setbox\mtagbox@\null\ltwidth@\z@\rtwidth@\z@
  \rmultline@@@#1$$}\fi\next@}
\def\rmultline@@#1\endmultline\tag#2$${\ltwidth@\z@
 \setbox\mtagbox@\hbox{\maketag@#2\maketag@}%
 \rmmeasure@#1\endmultline\dimen@\mwidth@\advance\dimen@\wd\mtagbox@
 \advance\dimen@\multlinetaggap@
 \ifdim\dimen@>\displaywidth\rtwidth@\z@\else\rtwidth@\wd\mtagbox@\fi
 \rmultline@@@#1\endmultline$$}
\def\rmultline@@@{\displ@y
 \def\shoveright##1{##1\hfilneg\iftagin@\ifdim\rtwidth@>\z@
  \hskip\rtwidth@\hskip\multlinetaggap@\fi\else\hskip\multlinegap@\fi}%
 \def\shoveleft##1{\setboxz@h{$\m@th\displaystyle{}##1$}%
  \setbox@ne\hbox{$\m@th\displaystyle##1$}%
  \hfilneg\hskip\multlinegap@\hskip.5\wd@ne\hskip-.5\wdz@##1}%
 \halign\bgroup\Let@\hbox to\displaywidth
  {\strut@$\m@th\displaystyle\hfil{}##\hfil$}\crcr
 \hfilneg\hskip\multlinegap@}
\def\endmultline{\iftagsleft@\expandafter\lendmultline@\else
 \expandafter\rendmultline@\fi}
\def\lendmultline@{\hfilneg\hskip\multlinegap@\crcr\egroup}
\def\rendmultline@{\iftagin@                                                %1
 \ifdim\rtwidth@>\z@                                                        %2
  \hskip\multlinetaggap@\box\mtagbox@                                       %3
 \else\llap{\vtop{\normalbaselines\null\hbox{\strut@\box\mtagbox@}}}\fi     %4
 \else\hskip\multlinegap@\fi                                                %5
 \hfilneg\crcr\egroup}
\def\bmod{\mskip-\medmuskip\mkern5mu\mathbin{\fam\z@ mod}\penalty900
 \mkern5mu\mskip-\medmuskip}
\def\pmod#1{\allowbreak\ifinner\mkern8mu\else\mkern18mu\fi
 ({\fam\z@ mod}\,\,#1)}
\def\pod#1{\allowbreak\ifinner\mkern8mu\else\mkern18mu\fi(#1)}
\def\mod#1{\allowbreak\ifinner\mkern12mu\else\mkern18mu\fi{\fam\z@ mod}\,\,#1}
\message{continued fractions,}
\newcount\cfraccount@
\def\cfrac{\bgroup\bgroup\advance\cfraccount@\@ne\strut
 \iffalse{\fi\def\\{\over\displaystyle}\iffalse}\fi}
\def\lcfrac{\bgroup\bgroup\advance\cfraccount@\@ne\strut
 \iffalse{\fi\def\\{\hfill\over\displaystyle}\iffalse}\fi}
\def\rcfrac{\bgroup\bgroup\advance\cfraccount@\@ne\strut\hfill
 \iffalse{\fi\def\\{\over\displaystyle}\iffalse}\fi}
\def\gloop@#1\repeat{\gdef\body{#1}\iterate}
\def\endcfrac{\gloop@\ifnum\cfraccount@>\z@\global\advance\cfraccount@\m@ne
 \egroup\hskip-\nulldelimiterspace\egroup\repeat}
\message{compound symbols,}
\def\binrel@#1{\setboxz@h{\thinmuskip0mu
  \medmuskip\m@ne mu\thickmuskip\@ne mu$#1\m@th$}%
 \setbox@ne\hbox{\thinmuskip0mu\medmuskip\m@ne mu\thickmuskip
  \@ne mu${}#1{}\m@th$}%
 \setbox\tw@\hbox{\hskip\wd@ne\hskip-\wdz@}}
\def\overset#1\to#2{\binrel@{#2}\ifdim\wd\tw@<\z@
 \mathbin{\mathop{\kern\z@#2}\limits^{#1}}\else\ifdim\wd\tw@>\z@
 \mathrel{\mathop{\kern\z@#2}\limits^{#1}}\else
 {\mathop{\kern\z@#2}\limits^{#1}}{}\fi\fi}
\def\underset#1\to#2{\binrel@{#2}\ifdim\wd\tw@<\z@
 \mathbin{\mathop{\kern\z@#2}\limits_{#1}}\else\ifdim\wd\tw@>\z@
 \mathrel{\mathop{\kern\z@#2}\limits_{#1}}\else
 {\mathop{\kern\z@#2}\limits_{#1}}{}\fi\fi}
\def\oversetbrace#1\to#2{\overbrace{#2}^{#1}}
\def\undersetbrace#1\to#2{\underbrace{#2}_{#1}}
\def\sideset#1\and#2\to#3{%
 \setbox@ne\hbox{$\dsize{\vphantom{#3}}#1{#3}\m@th$}%
 \setbox\tw@\hbox{$\dsize{#3}#2\m@th$}%
 \hskip\wd@ne\hskip-\wd\tw@\mathop{\hskip\wd\tw@\hskip-\wd@ne
  {\vphantom{#3}}#1{#3}#2}}
\def\rightarrowfill@#1{\setboxz@h{$#1-\m@th$}\ht\z@\z@
  $#1\m@th\copy\z@\mkern-6mu\cleaders
  \hbox{$#1\mkern-2mu\box\z@\mkern-2mu$}\hfill
  \mkern-6mu\mathord\rightarrow$}
\def\leftarrowfill@#1{\setboxz@h{$#1-\m@th$}\ht\z@\z@
  $#1\m@th\mathord\leftarrow\mkern-6mu\cleaders
  \hbox{$#1\mkern-2mu\copy\z@\mkern-2mu$}\hfill
  \mkern-6mu\box\z@$}
\def\leftrightarrowfill@#1{\setboxz@h{$#1-\m@th$}\ht\z@\z@
  $#1\m@th\mathord\leftarrow\mkern-6mu\cleaders
  \hbox{$#1\mkern-2mu\box\z@\mkern-2mu$}\hfill
  \mkern-6mu\mathord\rightarrow$}
\def\overrightarrow{\mathpalette\overrightarrow@}
\def\overrightarrow@#1#2{\vbox{\ialign{##\crcr\rightarrowfill@#1\crcr
 \noalign{\kern-\ex@\nointerlineskip}$\m@th\hfil#1#2\hfil$\crcr}}}

\def\overleftarrow{\mathpalette\overleftarrow@}
\def\overleftarrow@#1#2{\vbox{\ialign{##\crcr\leftarrowfill@#1\crcr
 \noalign{\kern-\ex@\nointerlineskip}$\m@th\hfil#1#2\hfil$\crcr}}}
\def\overleftrightarrow{\mathpalette\overleftrightarrow@}
\def\overleftrightarrow@#1#2{\vbox{\ialign{##\crcr\leftrightarrowfill@#1\crcr
 \noalign{\kern-\ex@\nointerlineskip}$\m@th\hfil#1#2\hfil$\crcr}}}
\def\underrightarrow{\mathpalette\underrightarrow@}
\def\underrightarrow@#1#2{\vtop{\ialign{##\crcr$\m@th\hfil#1#2\hfil$\crcr
 \noalign{\nointerlineskip}\rightarrowfill@#1\crcr}}}

\def\underleftarrow{\mathpalette\underleftarrow@}
\def\underleftarrow@#1#2{\vtop{\ialign{##\crcr$\m@th\hfil#1#2\hfil$\crcr
 \noalign{\nointerlineskip}\leftarrowfill@#1\crcr}}}
\def\underleftrightarrow{\mathpalette\underleftrightarrow@}
\def\underleftrightarrow@#1#2{\vtop{\ialign{##\crcr$\m@th\hfil#1#2\hfil$\crcr
 \noalign{\nointerlineskip}\leftrightarrowfill@#1\crcr}}}
\message{various kinds of dots,}
\let\DOTSI\relax
\let\DOTSB\relax

\newif\ifmath@
{\uccode`7=`\\ \uccode`8=`m \uccode`9=`a \uccode`0=`t \uccode`!=`h
 \uppercase{\gdef\math@#1#2#3#4#5#6\math@{\global\math@false\ifx 7#1\ifx 8#2%
 \ifx 9#3\ifx 0#4\ifx !#5\xdef\meaning@{#6}\global\math@true\fi\fi\fi\fi\fi}}}
\newif\ifmathch@
{\uccode`7=`c \uccode`8=`h \uccode`9=`\"
 \uppercase{\gdef\mathch@#1#2#3#4#5#6\mathch@{\global\mathch@false
  \ifx 7#1\ifx 8#2\ifx 9#5\global\mathch@true\xdef\meaning@{9#6}\fi\fi\fi}}}
\newcount\classnum@
\def\getmathch@#1.#2\getmathch@{\classnum@#1 \divide\classnum@4096
 \ifcase\number\classnum@\or\or\gdef\thedots@{\dotsb@}\or
 \gdef\thedots@{\dotsb@}\fi}
\newif\ifmathbin@
{\uccode`4=`b \uccode`5=`i \uccode`6=`n
 \uppercase{\gdef\mathbin@#1#2#3{\relaxnext@
  \DNii@##1\mathbin@{\ifx\space@\next\global\mathbin@true\fi}%
 \global\mathbin@false\DN@##1\mathbin@{}%
 \ifx 4#1\ifx 5#2\ifx 6#3\DN@{\FN@\nextii@}\fi\fi\fi\next@}}}
\newif\ifmathrel@
{\uccode`4=`r \uccode`5=`e \uccode`6=`l
 \uppercase{\gdef\mathrel@#1#2#3{\relaxnext@
  \DNii@##1\mathrel@{\ifx\space@\next\global\mathrel@true\fi}%
 \global\mathrel@false\DN@##1\mathrel@{}%
 \ifx 4#1\ifx 5#2\ifx 6#3\DN@{\FN@\nextii@}\fi\fi\fi\next@}}}
\newif\ifmacro@
{\uccode`5=`m \uccode`6=`a \uccode`7=`c
 \uppercase{\gdef\macro@#1#2#3#4\macro@{\global\macro@false
  \ifx 5#1\ifx 6#2\ifx 7#3\global\macro@true
  \xdef\meaning@{\macro@@#4\macro@@}\fi\fi\fi}}}
\def\macro@@#1->#2\macro@@{#2}
\newif\ifDOTS@
\newcount\DOTSCASE@
{\uccode`6=`\\ \uccode`7=`D \uccode`8=`O \uccode`9=`T \uccode`0=`S
 \uppercase{\gdef\DOTS@#1#2#3#4#5{\global\DOTS@false\DN@##1\DOTS@{}%
  \ifx 6#1\ifx 7#2\ifx 8#3\ifx 9#4\ifx 0#5\let\next@\DOTS@@\fi\fi\fi\fi\fi
  \next@}}}
{\uccode`3=`B \uccode`4=`I \uccode`5=`X
 \uppercase{\gdef\DOTS@@#1{\relaxnext@
  \DNii@##1\DOTS@{\ifx\space@\next\global\DOTS@true\fi}%
  \DN@{\FN@\nextii@}%
  \ifx 3#1\global\DOTSCASE@\z@\else
  \ifx 4#1\global\DOTSCASE@\@ne\else
  \ifx 5#1\global\DOTSCASE@\tw@\else\DN@##1\DOTS@{}%
  \fi\fi\fi\next@}}}
\newif\ifnot@
{\uccode`5=`\\ \uccode`6=`n \uccode`7=`o \uccode`8=`t
 \uppercase{\gdef\not@#1#2#3#4{\relaxnext@
  \DNii@##1\not@{\ifx\space@\next\global\not@true\fi}%
 \global\not@false\DN@##1\not@{}%
 \ifx 5#1\ifx 6#2\ifx 7#3\ifx 8#4\DN@{\FN@\nextii@}\fi\fi\fi
 \fi\next@}}}
\newif\ifkeybin@
\def\keybin@{\keybin@true
 \ifx\next+\else\ifx\next=\else\ifx\next<\else\ifx\next>\else\ifx\next-\else
 \ifx\next*\else\ifx\next:\else\keybin@false\fi\fi\fi\fi\fi\fi\fi}
\def\dots{\RIfM@\expandafter\mdots@\else\expandafter\tdots@\fi}
\def\tdots@{\unskip\relaxnext@
 \DN@{$\m@th\mathinner{\ldotp\ldotp\ldotp}\,
   \ifx\next,\,$\else\ifx\next.\,$\else\ifx\next;\,$\else\ifx\next:\,$\else
   \ifx\next?\,$\else\ifx\next!\,$\else$ \fi\fi\fi\fi\fi\fi}%
 \ \FN@\next@}
\def\mdots@{\FN@\mdots@@}
\def\mdots@@{\gdef\thedots@{\dotso@}%                                       %1
 \ifx\next\boldkey\gdef\thedots@\boldkey{\boldkeydots@}\else                %2
 \ifx\next\boldsymbol\gdef\thedots@\boldsymbol{\boldsymboldots@}\else       %3
 \ifx,\next\gdef\thedots@{\dotsc}%                                          %4
 \else\ifx\not\next\gdef\thedots@{\dotsb@}%                                 %5
 \else\keybin@
 \ifkeybin@\gdef\thedots@{\dotsb@}%                                         %6
 \else\xdef\meaning@{\meaning\next..........}\xdef\meaning@@{\meaning@}%    %7
  \expandafter\math@\meaning@\math@
  \ifmath@
   \expandafter\mathch@\meaning@\mathch@
   \ifmathch@\expandafter\getmathch@\meaning@\getmathch@\fi                 %8
  \else\expandafter\macro@\meaning@@\macro@                                 %9
  \ifmacro@                                                                %10
   \expandafter\not@\meaning@\not@\ifnot@\gdef\thedots@{\dotsb@}%          %11
  \else\expandafter\DOTS@\meaning@\DOTS@
  \ifDOTS@
   \ifcase\number\DOTSCASE@\gdef\thedots@{\dotsb@}%
    \or\gdef\thedots@{\dotsi}\else\fi                                      %12
  \else\expandafter\math@\meaning@\math@                                   %13
  \ifmath@\expandafter\mathbin@\meaning@\mathbin@
  \ifmathbin@\gdef\thedots@{\dotsb@}%                                      %14
  \else\expandafter\mathrel@\meaning@\mathrel@
  \ifmathrel@\gdef\thedots@{\dotsb@}%                                      %15
  \fi\fi\fi\fi\fi\fi\fi\fi\fi\fi\fi\fi
 \thedots@}
\def\plainldots@{\mathinner{\ldotp\ldotp\ldotp}}
\def\plaincdots@{\mathinner{\cdotp\cdotp\cdotp}}
\def\dotsi{\!\plaincdots@}
\let\dotsb@\plaincdots@
\newif\ifextra@
\newif\ifrightdelim@
\def\rightdelim@{\global\rightdelim@true                                    %1
 \ifx\next)\else                                                            %2
 \ifx\next]\else
 \ifx\next\rbrack\else
 \ifx\next\}\else
 \ifx\next\rbrace\else
 \ifx\next\rangle\else
 \ifx\next\rceil\else
 \ifx\next\rfloor\else
 \ifx\next\rgroup\else
 \ifx\next\rmoustache\else
 \ifx\next\right\else
 \ifx\next\bigr\else
 \ifx\next\biggr\else
 \ifx\next\Bigr\else                                                        %3
 \ifx\next\Biggr\else\global\rightdelim@false
 \fi\fi\fi\fi\fi\fi\fi\fi\fi\fi\fi\fi\fi\fi\fi}
\def\extra@{%
 \global\extra@false\rightdelim@\ifrightdelim@\global\extra@true            %1
 \else\ifx\next$\global\extra@true                                          %2
 \else\xdef\meaning@{\meaning\next..........}%                              %3
 \expandafter\macro@\meaning@\macro@\ifmacro@                               %4
 \expandafter\DOTS@\meaning@\DOTS@
 \ifDOTS@
 \ifnum\DOTSCASE@=\tw@\global\extra@true                                    %5
 \fi\fi\fi\fi\fi}
\newif\ifbold@
\def\dotso@{\relaxnext@
 \ifbold@
  \let\next\delayed@
  \DNii@{\extra@\plainldots@\ifextra@\,\fi}%
 \else
  \DNii@{\DN@{\extra@\plainldots@\ifextra@\,\fi}\FN@\next@}%
 \fi
 \nextii@}
\def\extrap@#1{%
 \ifx\next,\DN@{#1\,}\else
 \ifx\next;\DN@{#1\,}\else
 \ifx\next.\DN@{#1\,}\else\extra@
 \ifextra@\DN@{#1\,}\else
 \let\next@#1\fi\fi\fi\fi\next@}
\def\ldots{\DN@{\extrap@\plainldots@}%
 \FN@\next@}
\def\cdots{\DN@{\extrap@\plaincdots@}%
 \FN@\next@}

\def\dotsc{\relaxnext@
 \DN@{\ifx\next;\plainldots@\,\else
  \ifx\next.\plainldots@\,\else\extra@\plainldots@
  \ifextra@\,\fi\fi\fi}%
 \FN@\next@}
\def\cdot{\mathchar"2201 }
\def\longrightarrow{\DOTSB\relbar\joinrel\rightarrow}
\def\Longrightarrow{\DOTSB\Relbar\joinrel\Rightarrow}

\def\hookrightarrow{\DOTSB\lhook\joinrel\rightarrow}

\message{special superscripts,}
\def\dddot#1{{\mathop{#1}\limits^{\vbox to-1.4\ex@{\kern-\tw@\ex@
 \hbox{\rm...}\vss}}}}
\def\ddddot#1{{\mathop{#1}\limits^{\vbox to-1.4\ex@{\kern-\tw@\ex@
 \hbox{\rm....}\vss}}}}
\def\sphat{^{\mathchoice{}{}%
 {\,\,\botsmash{\hbox{\lower4\ex@\hbox{$\m@th\widehat{\null}$}}}}%
 {\,\botsmash{\hbox{\lower3\ex@\hbox{$\m@th\hat{\null}$}}}}}}

\def\spacute{^{\!\botsmash{\hbox{\lower\@ne ex\hbox{\'{}}}}}}
\def\spgrave{^{\mathchoice{}{}{}{\!}%
 \botsmash{\hbox{\lower\@ne ex\hbox{\`{}}}}}}
\def\spdot{^{\hbox{\raise\ex@\hbox{\rm.}}}}
\def\spddot{^{\hbox{\raise\ex@\hbox{\rm..}}}}
\def\spdddot{^{\hbox{\raise\ex@\hbox{\rm...}}}}
\def\spddddot{^{\hbox{\raise\ex@\hbox{\rm....}}}}
\def\spbreve{^{\!\botsmash{\hbox{\lower4\ex@\hbox{\u{}}}}}}

\message{\string\text,}
\def\textonlyfont@#1#2{\def#1{\RIfM@
 \Err@{Use \string#1\space only in text}\else#2\fi}}
\textonlyfont@\rm\tenrm
\textonlyfont@\it\tenit
\textonlyfont@\sl\tensl
\textonlyfont@\bf\tenbf
\def\oldnos#1{\RIfM@{\mathcode`\,="013B \fam\@ne#1}\else
 \leavevmode\hbox{$\m@th\mathcode`\,="013B \fam\@ne#1$}\fi}
\def\text{\RIfM@\expandafter\text@\else\expandafter\text@@\fi}
\def\text@@#1{\leavevmode\hbox{#1}}
\def\mathhexbox@#1#2#3{\text{$\m@th\mathchar"#1#2#3$}}
\def\dag{{\mathhexbox@279}}
\def\ddag{{\mathhexbox@27A}}
\def\S{{\mathhexbox@278}}
\def\P{{\mathhexbox@27B}}
\newif\iffirstchoice@
\firstchoice@true
\def\text@#1{\mathchoice
 {\hbox{\everymath{\displaystyle}\def\textfonti{\the\textfont\@ne}%
  \def\textfontii{\the\textfont\tw@}\textdef@@ T#1}}
 {\hbox{\firstchoice@false
  \everymath{\textstyle}\def\textfonti{\the\textfont\@ne}%
  \def\textfontii{\the\textfont\tw@}\textdef@@ T#1}}
 {\hbox{\firstchoice@false
  \everymath{\scriptstyle}\def\textfonti{\the\scriptfont\@ne}%
  \def\textfontii{\the\scriptfont\tw@}\textdef@@ S\rm#1}}
 {\hbox{\firstchoice@false
  \everymath{\scriptscriptstyle}\def\textfonti
  {\the\scriptscriptfont\@ne}%
  \def\textfontii{\the\scriptscriptfont\tw@}\textdef@@ s\rm#1}}}
\def\textdef@@#1{\textdef@#1\rm\textdef@#1\bf\textdef@#1\sl\textdef@#1\it}
\def\rmfam{0}
\def\textdef@#1#2{%
 \DN@{\csname\expandafter\eat@\string#2fam\endcsname}%
 \if S#1\edef#2{\the\scriptfont\next@\relax}%
 \else\if s#1\edef#2{\the\scriptscriptfont\next@\relax}%
 \else\edef#2{\the\textfont\next@\relax}\fi\fi}
\scriptfont\itfam\tenit \scriptscriptfont\itfam\tenit
\scriptfont\slfam\tensl \scriptscriptfont\slfam\tensl
\newif\iftopfolded@
\newif\ifbotfolded@
\def\topfoldedtext{\topfolded@true\botfolded@false\foldedtext@}
\def\botfoldedtext{\botfolded@true\topfolded@false\foldedtext@}
\def\foldedtext{\topfolded@false\botfolded@false\foldedtext@}
\Invalid@\foldedwidth
\def\foldedtext@{\relaxnext@
 \DN@{\ifx\next\foldedwidth\let\next@\nextii@\else
  \DN@{\nextii@\foldedwidth{.3\hsize}}\fi\next@}%
 \DNii@\foldedwidth##1##2{\setbox\z@\vbox
  {\normalbaselines\hsize##1\relax
  \tolerance1600 \noindent\ignorespaces##2}\ifbotfolded@\boxz@\else
  \iftopfolded@\vtop{\unvbox\z@}\else\vcenter{\boxz@}\fi\fi}%
 \FN@\next@}
\message{math font commands,}
\def\bold{\RIfM@\expandafter\bold@\else
 \expandafter\nonmatherr@\expandafter\bold\fi}
\def\bold@#1{{\bold@@{#1}}}
\def\bold@@#1{\fam\bffam\relax#1}
\def\slanted{\RIfM@\expandafter\slanted@\else
 \expandafter\nonmatherr@\expandafter\slanted\fi}
\def\slanted@#1{{\slanted@@{#1}}}
\def\slanted@@#1{\fam\slfam\relax#1}
\def\roman{\RIfM@\expandafter\roman@\else
 \expandafter\nonmatherr@\expandafter\roman\fi}
\def\roman@#1{{\roman@@{#1}}}
\def\roman@@#1{\fam\rmfam\relax#1}
\def\italic{\RIfM@\expandafter\italic@\else
 \expandafter\nonmatherr@\expandafter\italic\fi}
\def\italic@#1{{\italic@@{#1}}}
\def\italic@@#1{\fam\itfam\relax#1}
\def\Cal{\RIfM@\expandafter\Cal@\else
 \expandafter\nonmatherr@\expandafter\Cal\fi}
\def\Cal@#1{{\Cal@@{#1}}}
\def\Cal@@#1{\noaccents@\fam\tw@#1}
\mathchardef\Gamma="0000
\mathchardef\Delta="0001
\mathchardef\Theta="0002
\mathchardef\Lambda="0003
\mathchardef\Xi="0004
\mathchardef\Pi="0005
\mathchardef\Sigma="0006
\mathchardef\Upsilon="0007
\mathchardef\Phi="0008
\mathchardef\Psi="0009
\mathchardef\Omega="000A
\mathchardef\varGamma="0100
\mathchardef\varDelta="0101
\mathchardef\varTheta="0102
\mathchardef\varLambda="0103
\mathchardef\varXi="0104
\mathchardef\varPi="0105
\mathchardef\varSigma="0106
\mathchardef\varUpsilon="0107
\mathchardef\varPhi="0108
\mathchardef\varPsi="0109
\mathchardef\varOmega="010A
\let\alloc@@\alloc@
\def\hexnumber@#1{\ifcase#1 0\or 1\or 2\or 3\or 4\or 5\or 6\or 7\or 8\or
 9\or A\or B\or C\or D\or E\or F\fi}
\def\loadmsam{%
 \font@\tenmsa=msam10
 \font@\sevenmsa=msam7
 \font@\fivemsa=msam5
 \alloc@@8\fam\chardef\sixt@@n\msafam
 \textfont\msafam=\tenmsa
 \scriptfont\msafam=\sevenmsa
 \scriptscriptfont\msafam=\fivemsa
 \edef\next{\hexnumber@\msafam}%
 \mathchardef\dabar@"0\next39
 \edef\dashrightarrow{\mathrel{\dabar@\dabar@\mathchar"0\next4B}}%
 \edef\dashleftarrow{\mathrel{\mathchar"0\next4C\dabar@\dabar@}}%
 \let\dasharrow\dashrightarrow
 \edef\ulcorner{\delimiter"4\next70\next70 }%
 \edef\urcorner{\delimiter"5\next71\next71 }%
 \edef\llcorner{\delimiter"4\next78\next78 }%
 \edef\lrcorner{\delimiter"5\next79\next79 }%
 \edef\yen{{\noexpand\mathhexbox@\next55}}%
 \edef\checkmark{{\noexpand\mathhexbox@\next58}}%
 \edef\circledR{{\noexpand\mathhexbox@\next72}}%
 \edef\maltese{{\noexpand\mathhexbox@\next7A}}%
 \global\let\loadmsam\empty}%
\def\loadmsbm{%
 \font@\tenmsb=msbm10 \font@\sevenmsb=msbm7 \font@\fivemsb=msbm5
 \alloc@@8\fam\chardef\sixt@@n\msbfam
 \textfont\msbfam=\tenmsb
 \scriptfont\msbfam=\sevenmsb \scriptscriptfont\msbfam=\fivemsb
 \global\let\loadmsbm\empty
 }
\def\widehat#1{\ifx\undefined\msbfam \DN@{362}%
  \else \setboxz@h{$\m@th#1$}%
    \edef\next@{\ifdim\wdz@>\tw@ em%
        \hexnumber@\msbfam 5B%
      \else 362\fi}\fi
  \mathaccent"0\next@{#1}}
\def\widetilde#1{\ifx\undefined\msbfam \DN@{365}%
  \else \setboxz@h{$\m@th#1$}%
    \edef\next@{\ifdim\wdz@>\tw@ em%
        \hexnumber@\msbfam 5D%
      \else 365\fi}\fi
  \mathaccent"0\next@{#1}}
\message{\string\newsymbol,}
\def\newsymbol#1#2#3#4#5{\define#1{}%
  \count@#2\relax \advance\count@\m@ne % to push case 0 to the \else clause
 \ifcase\count@
   \ifx\undefined\msafam\loadmsam\fi \let\next@\msafam
 \or \ifx\undefined\msbfam\loadmsbm\fi \let\next@\msbfam
 \else  \Err@{\Invalid@@\string\newsymbol}\let\next@\tw@\fi
 \mathchardef#1="#3\hexnumber@\next@#4#5\space}

\def\Bbb{\RIfM@\expandafter\Bbb@\else
 \expandafter\nonmatherr@\expandafter\Bbb\fi}
\def\Bbb@#1{{\Bbb@@{#1}}}
\def\Bbb@@#1{\noaccents@\fam\msbfam\relax#1}
\message{bold Greek and bold symbols,}
\def\loadbold{%
 \font@\tencmmib=cmmib10 \font@\sevencmmib=cmmib7 \font@\fivecmmib=cmmib5
 \skewchar\tencmmib'177 \skewchar\sevencmmib'177 \skewchar\fivecmmib'177
 \alloc@@8\fam\chardef\sixt@@n\cmmibfam
 \textfont\cmmibfam\tencmmib
 \scriptfont\cmmibfam\sevencmmib \scriptscriptfont\cmmibfam\fivecmmib
 \font@\tencmbsy=cmbsy10 \font@\sevencmbsy=cmbsy7 \font@\fivecmbsy=cmbsy5
 \skewchar\tencmbsy'60 \skewchar\sevencmbsy'60 \skewchar\fivecmbsy'60
 \alloc@@8\fam\chardef\sixt@@n\cmbsyfam
 \textfont\cmbsyfam\tencmbsy
 \scriptfont\cmbsyfam\sevencmbsy \scriptscriptfont\cmbsyfam\fivecmbsy
 \let\loadbold\empty
}
\def\boldnotloaded#1{\Err@{\ifcase#1\or First\else Second\fi
       bold symbol font not loaded}}
\def\mathchari@#1#2#3{\ifx\undefined\cmmibfam
    \boldnotloaded@\@ne
  \else\mathchar"#1\hexnumber@\cmmibfam#2#3\space \fi}
\def\mathcharii@#1#2#3{\ifx\undefined\cmbsyfam
    \boldnotloaded\tw@
  \else \mathchar"#1\hexnumber@\cmbsyfam#2#3\space\fi}
\edef\bffam@{\hexnumber@\bffam}
\def\boldkey#1{\ifcat\noexpand#1A%
  \ifx\undefined\cmmibfam \boldnotloaded\@ne
  \else {\fam\cmmibfam#1}\fi
 \else
 \ifx#1!\mathchar"5\bffam@21 \else
 \ifx#1(\mathchar"4\bffam@28 \else\ifx#1)\mathchar"5\bffam@29 \else
 \ifx#1+\mathchar"2\bffam@2B \else\ifx#1:\mathchar"3\bffam@3A \else
 \ifx#1;\mathchar"6\bffam@3B \else\ifx#1=\mathchar"3\bffam@3D \else
 \ifx#1?\mathchar"5\bffam@3F \else\ifx#1[\mathchar"4\bffam@5B \else
 \ifx#1]\mathchar"5\bffam@5D \else
 \ifx#1,\mathchari@63B \else
 \ifx#1-\mathcharii@200 \else
 \ifx#1.\mathchari@03A \else
 \ifx#1/\mathchari@03D \else
 \ifx#1<\mathchari@33C \else
 \ifx#1>\mathchari@33E \else
 \ifx#1*\mathcharii@203 \else
 \ifx#1|\mathcharii@06A \else
 \ifx#10\bold0\else\ifx#11\bold1\else\ifx#12\bold2\else\ifx#13\bold3\else
 \ifx#14\bold4\else\ifx#15\bold5\else\ifx#16\bold6\else\ifx#17\bold7\else
 \ifx#18\bold8\else\ifx#19\bold9\else
  \Err@{\string\boldkey\space can't be used with #1}%
 \fi\fi\fi\fi\fi\fi\fi\fi\fi\fi\fi\fi\fi\fi\fi
 \fi\fi\fi\fi\fi\fi\fi\fi\fi\fi\fi\fi\fi\fi}
\def\boldsymbol#1{%
 \DN@{\Err@{You can't use \string\boldsymbol\space with \string#1}#1}%
 \ifcat\noexpand#1A%
   \let\next@\relax
   \ifx\undefined\cmmibfam \boldnotloaded\@ne
   \else {\fam\cmmibfam#1}\fi
 \else
  \xdef\meaning@{\meaning#1.........}%
  \expandafter\math@\meaning@\math@
  \ifmath@
   \expandafter\mathch@\meaning@\mathch@
   \ifmathch@
    \expandafter\boldsymbol@@\meaning@\boldsymbol@@
   \fi
  \else
   \expandafter\macro@\meaning@\macro@
   \expandafter\delim@\meaning@\delim@
   \ifdelim@
    \expandafter\delim@@\meaning@\delim@@
   \else
    \boldsymbol@{#1}%
   \fi
  \fi
 \fi
 \next@}
\def\mathhexboxii@#1#2{\ifx\undefined\cmbsyfam
    \boldnotloaded\tw@
  \else \mathhexbox@{\hexnumber@\cmbsyfam}{#1}{#2}\fi}
\def\boldsymbol@#1{\let\next@\relax\let\next#1%
 \ifx\next\cdot\mathcharii@201 \else
 \ifx\next\prime{{\null\mathcharii@030 \null}}\else
 \ifx\next\lbrack\mathchar"4\bffam@5B \else
 \ifx\next\rbrack\mathchar"5\bffam@5D \else
 \ifx\next\{\mathcharii@466 \else
 \ifx\next\lbrace\mathcharii@466 \else
 \ifx\next\}\mathcharii@567 \else
 \ifx\next\rbrace\mathcharii@567 \else
 \ifx\next\surd{{\mathcharii@170}}\else
 \ifx\next\S{{\mathhexboxii@78}}\else
 \ifx\next\P{{\mathhexboxii@7B}}\else
 \ifx\next\dag{{\mathhexboxii@79}}\else
 \ifx\next\ddag{{\mathhexboxii@7A}}\else
 \DN@{\Err@{You can't use \string\boldsymbol\space with \string#1}#1}%
 \fi\fi\fi\fi\fi\fi\fi\fi\fi\fi\fi\fi\fi}
\def\boldsymbol@@#1.#2\boldsymbol@@{\classnum@#1 \count@@@\classnum@        %1
 \divide\classnum@4096 \count@\classnum@                                    %2
 \multiply\count@4096 \advance\count@@@-\count@ \count@@\count@@@           %3
 \divide\count@@@\@cclvi \count@\count@@                                    %4
 \multiply\count@@@\@cclvi \advance\count@@-\count@@@                       %5
 \divide\count@@@\@cclvi                                                    %6
 \multiply\classnum@4096 \advance\classnum@\count@@                         %7
 \ifnum\count@@@=\z@                                                        %8
  \count@"\bffam@ \multiply\count@\@cclvi
  \advance\classnum@\count@
  \DN@{\mathchar\number\classnum@}%
 \else
  \ifnum\count@@@=\@ne                                                      %9
   \ifx\undefined\cmmibfam \DN@{\boldnotloaded\@ne}%
   \else \count@\cmmibfam \multiply\count@\@cclvi
     \advance\classnum@\count@
     \DN@{\mathchar\number\classnum@}\fi
  \else
   \ifnum\count@@@=\tw@                                                    %10
     \ifx\undefined\cmbsyfam
       \DN@{\boldnotloaded\tw@}%
     \else
       \count@\cmbsyfam \multiply\count@\@cclvi
       \advance\classnum@\count@
       \DN@{\mathchar\number\classnum@}%
     \fi
  \fi
 \fi
\fi}
\newif\ifdelim@
\newcount\delimcount@
{\uccode`6=`\\ \uccode`7=`d \uccode`8=`e \uccode`9=`l
 \uppercase{\gdef\delim@#1#2#3#4#5\delim@
  {\delim@false\ifx 6#1\ifx 7#2\ifx 8#3\ifx 9#4\delim@true
   \xdef\meaning@{#5}\fi\fi\fi\fi}}}
\def\delim@@#1"#2#3#4#5#6\delim@@{\if#32%
\let\next@\relax
 \ifx\undefined\cmbsyfam \boldnotloaded\@ne
 \else \mathcharii@#2#4#5\space \fi\fi}
\def\vert{\delimiter"026A30C }
\def\Vert{\delimiter"026B30D }
\let\|\Vert

\def\boldkeydots@#1{\bold@true\let\next=#1\let\delayed@=#1\mdots@@
 \boldkey#1\bold@false}  % = required!
\def\boldsymboldots@#1{\bold@true\let\next#1\let\delayed@#1\mdots@@
 \boldsymbol#1\bold@false}
\message{Euler fonts,}
\def\loadeufm{\loadmathfont{eufm}}

\def\frak{\mathfont@\frak}

\def\loadmathfont#1{% 
   \expandafter\font@\csname ten#1\endcsname=#110
   \expandafter\font@\csname seven#1\endcsname=#17
   \expandafter\font@\csname five#1\endcsname=#15
   \edef\next{\noexpand\alloc@@8\fam\chardef\sixt@@n
     \expandafter\noexpand\csname#1fam\endcsname}%
   \next
   \textfont\csname#1fam\endcsname \csname ten#1\endcsname
   \scriptfont\csname#1fam\endcsname \csname seven#1\endcsname
   \scriptscriptfont\csname#1fam\endcsname \csname five#1\endcsname
   \expandafter\def\csname #1\expandafter\endcsname\expandafter{%
      \expandafter\mathfont@\csname#1\endcsname}%
 \expandafter\gdef\csname load#1\endcsname{}%
}
\def\mathfont@#1{\RIfM@\expandafter\mathfont@@\expandafter#1\else
  \expandafter\nonmatherr@\expandafter#1\fi}
\def\mathfont@@#1#2{{\mathfont@@@#1{#2}}}
\def\mathfont@@@#1#2{\noaccents@
   \fam\csname\expandafter\eat@\string#1fam\endcsname
   \relax#2}
\message{math accents,}
\def\accentclass@{7}
\def\noaccents@{\def\accentclass@{0}}
\def\makeacc@#1#2{\def#1{\mathaccent"\accentclass@#2 }}
\makeacc@\hat{05E}
\makeacc@\check{014}
\makeacc@\tilde{07E}
\makeacc@\acute{013}
\makeacc@\grave{012}
\makeacc@\dot{05F}
\makeacc@\ddot{07F}
\makeacc@\breve{015}
\makeacc@\bar{016}

\newcount\skewcharcount@
\newcount\familycount@
\def\theskewchar@{\familycount@\@ne
 \global\skewcharcount@\the\skewchar\textfont\@ne                           %1
 \ifnum\fam>\m@ne\ifnum\fam<16
  \global\familycount@\the\fam\relax
  \global\skewcharcount@\the\skewchar\textfont\the\fam\relax\fi\fi          %2
 \ifnum\skewcharcount@>\m@ne
  \ifnum\skewcharcount@<128
  \multiply\familycount@256
  \global\advance\skewcharcount@\familycount@
  \global\advance\skewcharcount@28672
  \mathchar\skewcharcount@\else
  \global\skewcharcount@\m@ne\fi\else
 \global\skewcharcount@\m@ne\fi}                                            %3
\newcount\pointcount@
\def\getpoints@#1.#2\getpoints@{\pointcount@#1 }
\newdimen\accentdimen@
\newcount\accentmu@
\def\dimentomu@{\multiply\accentdimen@ 100
 \expandafter\getpoints@\the\accentdimen@\getpoints@
 \multiply\pointcount@18
 \divide\pointcount@\@m
 \global\accentmu@\pointcount@}
\def\Makeacc@#1#2{\def#1{\RIfM@\DN@{\mathaccent@
 {"\accentclass@#2 }}\else\DN@{\nonmatherr@{#1}}\fi\next@}}
\def\unbracefonts@{\let\Cal@\Cal@@\let\roman@\roman@@\let\bold@\bold@@
 \let\slanted@\slanted@@}
\def\mathaccent@#1#2{\ifnum\fam=\m@ne\xdef\thefam@{1}\else
 \xdef\thefam@{\the\fam}\fi                                                 %1
 \accentdimen@\z@                                                           %2
 \setboxz@h{\unbracefonts@$\m@th\fam\thefam@\relax#2$}%                     %3
 \ifdim\accentdimen@=\z@\DN@{\mathaccent#1{#2}}%                            %4
  \setbox@ne\hbox{\unbracefonts@$\m@th\fam\thefam@\relax#2\theskewchar@$}% %5a
  \setbox\tw@\hbox{$\m@th\ifnum\skewcharcount@=\m@ne\else
   \mathchar\skewcharcount@\fi$}%                                          %5b
  \global\accentdimen@\wd@ne\global\advance\accentdimen@-\wdz@
  \global\advance\accentdimen@-\wd\tw@                                     %5c
  \global\multiply\accentdimen@\tw@
  \dimentomu@\global\advance\accentmu@\@ne                                 %5d
 \else\DN@{{\mathaccent#1{#2\mkern\accentmu@ mu}%
    \mkern-\accentmu@ mu}{}}\fi                                             %6
 \next@}\Makeacc@\Hat{05E}
\Makeacc@\Check{014}
\Makeacc@\Tilde{07E}
\Makeacc@\Acute{013}
\Makeacc@\Grave{012}
\Makeacc@\Dot{05F}
\Makeacc@\Ddot{07F}
\Makeacc@\Breve{015}
\Makeacc@\Bar{016}
\def\Vec{\RIfM@\DN@{\mathaccent@{"017E }}\else
 \DN@{\nonmatherr@\Vec}\fi\next@}
\def\accentedsymbol#1#2{\csname newbox\expandafter\endcsname
  \csname\expandafter\eat@\string#1@box\endcsname
 \expandafter\setbox\csname\expandafter\eat@
  \string#1@box\endcsname\hbox{$\m@th#2$}\define
  #1{\copy\csname\expandafter\eat@\string#1@box\endcsname{}}}
\message{roots,}
\def\sqrt#1{\radical"270370 {#1}}
\let\underline@\underline
\let\overline@\overline
\def\underline#1{\underline@{#1}}
\def\overline#1{\overline@{#1}}
\Invalid@\leftroot
\Invalid@\uproot
\newcount\uproot@
\newcount\leftroot@
\def\root{\relaxnext@
  \DN@{\ifx\next\uproot\let\next@\nextii@\else
   \ifx\next\leftroot\let\next@\nextiii@\else
   \let\next@\plainroot@\fi\fi\next@}%
  \DNii@\uproot##1{\uproot@##1\relax\FN@\nextiv@}%
  \def\nextiv@{\ifx\next\space@\DN@. {\FN@\nextv@}\else
   \DN@.{\FN@\nextv@}\fi\next@.}%
  \def\nextv@{\ifx\next\leftroot\let\next@\nextvi@\else
   \let\next@\plainroot@\fi\next@}%
  \def\nextvi@\leftroot##1{\leftroot@##1\relax\plainroot@}%
   \def\nextiii@\leftroot##1{\leftroot@##1\relax\FN@\nextvii@}%
  \def\nextvii@{\ifx\next\space@
   \DN@. {\FN@\nextviii@}\else
   \DN@.{\FN@\nextviii@}\fi\next@.}%
  \def\nextviii@{\ifx\next\uproot\let\next@\nextix@\else
   \let\next@\plainroot@\fi\next@}%
  \def\nextix@\uproot##1{\uproot@##1\relax\plainroot@}%
  \bgroup\uproot@\z@\leftroot@\z@\FN@\next@}
\def\plainroot@#1\of#2{\setbox\rootbox\hbox{$\m@th\scriptscriptstyle{#1}$}%
 \mathchoice{\r@@t\displaystyle{#2}}{\r@@t\textstyle{#2}}
 {\r@@t\scriptstyle{#2}}{\r@@t\scriptscriptstyle{#2}}\egroup}
\def\r@@t#1#2{\setboxz@h{$\m@th#1\sqrt{#2}$}%
 \dimen@\ht\z@\advance\dimen@-\dp\z@
 \setbox@ne\hbox{$\m@th#1\mskip\uproot@ mu$}\advance\dimen@ 1.667\wd@ne
 \mkern-\leftroot@ mu\mkern5mu\raise.6\dimen@\copy\rootbox
 \mkern-10mu\mkern\leftroot@ mu\boxz@}
\def\boxed#1{\setboxz@h{$\m@th\displaystyle{#1}$}\dimen@.4\ex@
 \advance\dimen@3\ex@\advance\dimen@\dp\z@
 \hbox{\lower\dimen@\hbox{%
 \vbox{\hrule height.4\ex@
 \hbox{\vrule width.4\ex@\hskip3\ex@\vbox{\vskip3\ex@\boxz@\vskip3\ex@}%
 \hskip3\ex@\vrule width.4\ex@}\hrule height.4\ex@}%
 }}}
\message{commutative diagrams,}
\let\ampersand@\relax
\newdimen\minaw@
\minaw@11.11128\ex@
\newdimen\minCDaw@
\minCDaw@2.5pc
\def\minCDarrowwidth#1{\RIfMIfI@\onlydmatherr@\minCDarrowwidth
 \else\minCDaw@#1\relax\fi\else\onlydmatherr@\minCDarrowwidth\fi}
\newif\ifCD@
\def\CD{\bgroup\vspace@\relax\iffalse{\fi\let\ampersand@&\iffalse}\fi
 \CD@true\vcenter\bgroup\Let@\tabskip\z@skip\baselineskip20\ex@
 \lineskip3\ex@\lineskiplimit3\ex@\halign\bgroup
 &\hfill$\m@th##$\hfill\crcr}
\def\endCD{\crcr\egroup\egroup\egroup}
\newdimen\bigaw@
\atdef@>#1>#2>{\ampersand@                                                  %1
 \setboxz@h{$\m@th\ssize\;{#1}\;\;$}%                                       %2
 \setbox@ne\hbox{$\m@th\ssize\;{#2}\;\;$}%                                  %3
 \setbox\tw@\hbox{$\m@th#2$}%                                               %4
 \ifCD@\global\bigaw@\minCDaw@\else\global\bigaw@\minaw@\fi                 %5
 \ifdim\wdz@>\bigaw@\global\bigaw@\wdz@\fi
 \ifdim\wd@ne>\bigaw@\global\bigaw@\wd@ne\fi                                %6
 \ifCD@\enskip\fi                                                           %7
 \ifdim\wd\tw@>\z@
  \mathrel{\mathop{\hbox to\bigaw@{\rightarrowfill@\displaystyle}}%
    \limits^{#1}_{#2}}%                                                     %8
 \else\mathrel{\mathop{\hbox to\bigaw@{\rightarrowfill@\displaystyle}}%
    \limits^{#1}}\fi                                                        %9
 \ifCD@\enskip\fi                                                          %10
 \ampersand@}                                                              %11
\atdef@<#1<#2<{\ampersand@\setboxz@h{$\m@th\ssize\;\;{#1}\;$}%
 \setbox@ne\hbox{$\m@th\ssize\;\;{#2}\;$}\setbox\tw@\hbox{$\m@th#2$}%
 \ifCD@\global\bigaw@\minCDaw@\else\global\bigaw@\minaw@\fi
 \ifdim\wdz@>\bigaw@\global\bigaw@\wdz@\fi
 \ifdim\wd@ne>\bigaw@\global\bigaw@\wd@ne\fi
 \ifCD@\enskip\fi
 \ifdim\wd\tw@>\z@
  \mathrel{\mathop{\hbox to\bigaw@{\leftarrowfill@\displaystyle}}%
       \limits^{#1}_{#2}}\else
  \mathrel{\mathop{\hbox to\bigaw@{\leftarrowfill@\displaystyle}}%
       \limits^{#1}}\fi
 \ifCD@\enskip\fi\ampersand@}
\begingroup
 \catcode`\~=\active \lccode`\~=`\@
 \lowercase{%
  \global\atdef@)#1)#2){~>#1>#2>}
  \global\atdef@(#1(#2({~<#1<#2<}}
\endgroup
\atdef@ A#1A#2A{\llap{$\m@th\vcenter{\hbox
 {$\ssize#1$}}$}\Big\uparrow\rlap{$\m@th\vcenter{\hbox{$\ssize#2$}}$}&&}
\atdef@ V#1V#2V{\llap{$\m@th\vcenter{\hbox
 {$\ssize#1$}}$}\Big\downarrow\rlap{$\m@th\vcenter{\hbox{$\ssize#2$}}$}&&}
\atdef@={&\enskip\mathrel
 {\vbox{\hrule width\minCDaw@\vskip3\ex@\hrule width
 \minCDaw@}}\enskip&}
\atdef@|{\Big\Vert&&}
\atdef@\vert{\Big\Vert&&}
\def\pretend#1\haswidth#2{\setboxz@h{$\m@th\scriptstyle{#2}$}\hbox
 to\wdz@{\hfill$\m@th\scriptstyle{#1}$\hfill}}
\message{poor man's bold,}
\def\pmb{\RIfM@\expandafter\mathpalette\expandafter\pmb@\else
 \expandafter\pmb@@\fi}
\def\pmb@@#1{\leavevmode\setboxz@h{#1}%
   \dimen@-\wdz@
   \kern-.5\ex@\copy\z@
   \kern\dimen@\kern.25\ex@\raise.4\ex@\copy\z@
   \kern\dimen@\kern.25\ex@\box\z@
}
\def\binrel@@#1{\ifdim\wd2<\z@\mathbin{#1}\else\ifdim\wd\tw@>\z@
 \mathrel{#1}\else{#1}\fi\fi}
\newdimen\pmbraise@
%      Note: because of the use of \mathpalette, if \pmb@ is
%      applied to a single math italic character (or a single
%      character from some other slanted math font), the italic
%      correction will be added.  This will cause subscripts
%      to fall too far away from the character in some
%      cases, e.g., $\pmb{T}_1$ or $\pmb{\Cal T}_1$.
\def\pmb@#1#2{\setbox\thr@@\hbox{$\m@th#1{#2}$}%
 \setbox4\hbox{$\m@th#1\mkern.5mu$}\pmbraise@\wd4\relax
 \binrel@{#2}%
 \dimen@-\wd\thr@@
   \binrel@@{%
   \mkern-.8mu\copy\thr@@
   \kern\dimen@\mkern.4mu\raise\pmbraise@\copy\thr@@
   \kern\dimen@\mkern.4mu\box\thr@@
}}
\def\documentstyle#1{\W@{}\input #1.sty\relax}
\message{syntax check,}
\font\dummyft@=dummy
\fontdimen1 \dummyft@=\z@
\fontdimen2 \dummyft@=\z@
\fontdimen3 \dummyft@=\z@
\fontdimen4 \dummyft@=\z@
\fontdimen5 \dummyft@=\z@
\fontdimen6 \dummyft@=\z@
\fontdimen7 \dummyft@=\z@
\fontdimen8 \dummyft@=\z@
\fontdimen9 \dummyft@=\z@
\fontdimen10 \dummyft@=\z@
\fontdimen11 \dummyft@=\z@
\fontdimen12 \dummyft@=\z@
\fontdimen13 \dummyft@=\z@
\fontdimen14 \dummyft@=\z@
\fontdimen15 \dummyft@=\z@
\fontdimen16 \dummyft@=\z@
\fontdimen17 \dummyft@=\z@
\fontdimen18 \dummyft@=\z@
\fontdimen19 \dummyft@=\z@
\fontdimen20 \dummyft@=\z@
\fontdimen21 \dummyft@=\z@
\fontdimen22 \dummyft@=\z@
\def\fontlist@{\\{\tenrm}\\{\sevenrm}\\{\fiverm}\\{\teni}\\{\seveni}%
 \\{\fivei}\\{\tensy}\\{\sevensy}\\{\fivesy}\\{\tenex}\\{\tenbf}\\{\sevenbf}%
 \\{\fivebf}\\{\tensl}\\{\tenit}}
\def\font@#1=#2 {\rightappend@#1\to\fontlist@\font#1=#2 }
\def\dodummy@{{\def\\##1{\global\let##1\dummyft@}\fontlist@}}
\def\nopages@{\output{\setbox\z@\box\@cclv \deadcycles\z@}%
 \alloc@5\toks\toksdef\@cclvi\output}
\let\galleys\nopages@
\newif\ifsyntax@
\newcount\countxviii@
\def\syntax{\syntax@true\dodummy@\countxviii@\count18
 \loop\ifnum\countxviii@>\m@ne\textfont\countxviii@=\dummyft@
 \scriptfont\countxviii@=\dummyft@\scriptscriptfont\countxviii@=\dummyft@
 \advance\countxviii@\m@ne\repeat                                           %1
 \dummyft@\tracinglostchars\z@\nopages@\frenchspacing\hbadness\@M}
\def\first@#1#2\end{#1}
\def\printoptions{\W@{Do you want S(yntax check),
  G(alleys) or P(ages)?}%
 \message{Type S, G or P, followed by <return>: }%
 \begingroup % to localize the following change to \endlinechar:
 \endlinechar\m@ne % to prevent a space or \par in \ans@ from ^^M
 \read\m@ne to\ans@
%  Define \ans@ to uppercase itself, and default to P if the user
%  just pressed <return>.
 \edef\ans@{\uppercase{\def\noexpand\ans@{%
   \expandafter\first@\ans@ P\end}}}%
%  Cast the new definition of \ans@ outside the group
 \expandafter\endgroup\ans@
 \if\ans@ P% fine, no action needs to be taken
 \else \if\ans@ S\syntax
 \else \if\ans@ G\galleys
 \else\message{? Unknown option: \ans@; using the `pages' option.}%
 \fi\fi\fi}
\def\alloc@#1#2#3#4#5{\global\advance\count1#1by\@ne
 \ch@ck#1#4#2\allocationnumber=\count1#1
 \global#3#5=\allocationnumber
 \ifalloc@\wlog{\string#5=\string#2\the\allocationnumber}\fi}
\def\document{\def\alloclist@{}\def\fontlist@{}}

\let\plainfootnote\footnote

\let\footnote\undefined
\let\=\undefined
\let\>\undefined

\catcode`\@=\active
\message{... finished}

\font\piccolo=cmr10   at 8.50 truept
\font\piccolot=cmr10   at 10.80 truept
 \def\d{\roman{d}}  
\font\corsivo=eusm10
 \def\connab{\hbox{\corsivo A}(A,B)}
\def\connaab{\hbox{\corsivo A}(A,\bar A, B)}

\def\chen{\int_{\text {\rm Chen}}}
\def\im{{\text {\rm Im}}}
\def\chena{\int_{\text {\rm Chen}(A)}}
\def\chenav#1{\int_{\text{\rm Chen}(#1)}}

\def\quot#1#2{\refe\xdef#1{[\the\norefe]}\immediate\write\fileref
 { \hangindent\parindent}
\ignorespaces\immediate\write\fileref{[\the\norefe]{#2}\par}\ignorespaces
[\the\norefe]\ignorespaces} 
  
\loadeufm 
 
\pageno=0

\def\claim#1#2{ \edef\aah{\string #1} 
\xdef#1{Claim
\the\sectnum.\the\thmnum }  \bigbreak\noindent{\bf Claim
\the\sectnum.\the\thmnum.} [\aah] {\it#2}\global\advance\thmnum by 1
}\def\Lra{\Longrightarrow} 
\def\T{\text {\rm T}}  \def\TM{{\text {\rm T}M}} \def\Hol{{\text
{\rm Hol}}} \def\TP{{\text {\rm T}P}} 
\def\V{{\text {V}}}  
\def\wtilde{\widetilde}  \def\lap{\Lscr_A P}
\def\loop{\Lscr M} \def\pap{\Pscr_A P} 
\def\path{\Pscr M} \def\pathv#1{\Pscr {#1}} 
\def\loopvv#1#2{\Lscr_{#1}{#2}} \def\loopv#1{\Lscr {#1}}
\def\papv#1{\Pscr_{#1} P} 
\def\Adhmt{{\text {Ad}_{h_{\gamma}^{-1}(t) } }}
\def\wsigma{\widehat\sigma} \def\AutP{{\text{\rm Aut}}\;P} 
\def\Aut{\text {\rm Aut}}  \def\aut{\text {\rm aut}} 
\def\autP{\text { \rm aut}\;P}   
\def\DiffX{\text{ \rm Diff}\;X}

\def\hliftomo{\frak L(\omega,\Gamma,  p_0)}

\def\hliftabor{\frak
L((A(\kappa),-F_{A(\kappa)}+ B(\lambda)),\Gamma_r,  p_0)}
\def\hliftaoo{\frak L((A,0),\Gamma,  p_0)}

\def\hliftraoo{\frak L((A,0),\Gamma_r,  p_0)}

\def\partrap{{\frak
L}(A,\gamma,p)}  
\def\ja{{J_A}}
\def\jsuba{{j}_A}
\def\lieg{{\frak g}}
  \def\partrak{{\frak
L}(\omega+\kappa\eta,\gamma, u)}   \def\partrao{{\frak
L}(\omega,\gamma,u)}   \def\loopabh{\frak
L(\connab ,\Gamma,  p)} 
\def\holp{\holvo{\connab }{\Gamma}{  p} }
\def\holvo#1#2#3{{\text{\rm Hol}}_{#1}({#2},{#3})} 
\def\curvaab{\hbox{\corsivo F}(A,\bar A, B)}  
\def\curvab{\hbox{\corsivo F}(A,B)}  
\def\as{\`a\ }

% DERIVATE
\def\ddr#1{ 
{d{#1}\over d r}  }
\def\dsopv#1{ \left. {\partial {#1}\over \partial s}\right\vert
_{s=0 }}  \def\dtv#1{\left. {d{#1}\over dt}\right.}
\def\pdsv#1{{\left.  {\partial{#1}\over \partial s}\right.  }}
\def\pds{{\left.  {\partial\over \partial s}\right.  }}
\def\pdt{{\left. {\partial\over \partial t}\right. }} 
\def\pdtv#1{{\left. {\partial{#1}\over \partial t}\right. }} 
\def\dsv#1{\left. {d{#1}\over ds}\right.}  \def\dkov#1{ \left.
{d{#1}\over d \kappa}\right\vert _{  \kappa=0}   } \def\dsov#1{
\left. {d{#1}\over ds}\right\vert _{ s=0} } 
{ \nopagenumbers     
\title{Loop and Path Spaces and Four-Dimensional $\bold
{BF}$  Theories:}\medskip
\title{Connections, Holonomies and Observables} \vskip0.4in \centerline{
{\bf Alberto S. Cattaneo},$^{(a)}$\plainfootnote{} {\piccolo This work has
been partly supported by  research grants of the Ministero dell'Universit\as e della Ricerca
Scientifica e Tecnologica (MURST). 
Part of this work was completed while A.S.C. was at
Harvard University supported by I.N.F.N. Grant No. 5565/95 and DOE Grant No.\ 
DE-FG02-94ER25228, Amendment No.\ A003. P.C.--R. developed some of the work related
to this paper while participating to the APCTP/PIms Summer Workshop at the University
of Vancouver, B.C.}}  \centerline{{\bf Paolo Cotta-Ramusino},$^{(a)}$
and {\bf Maurizio Rinaldi}$^{(b)}$ } \vskip12pt
\centerline{$^{(a)}$Dipartimento di Matematica}
\centerline{Universit\as di Milano} \centerline{Via Saldini 50
}\centerline{20133 Milano, Italy} \centerline{and}
\centerline{I.N.F.N., Sezione di Milano} \vskip12pt
\centerline{$^{(b)}$Dipartimento di Matematica} 
\centerline{Universit\as di Trieste} \centerline{Piazzale Europa 1}
\centerline{34127 Trieste, Italy} \vfill  \noindent
\plainfootnote{}{e-mail addresses:
{\piccolot cattaneo\@elanor.mat.unimi.it,
cotta\@mi.infn.it, 
rinaldi\@univ.trieste.it.}}\vskip-0.4in 
 
\vskip-0.25in 
\centerline{\bf Abstract}
 
{\piccolot We study  the differential geometry 
of principal $G$-bundles
whose base space is
the space of free paths (loops) on a  manifold $M$.
In particular we consider connections
defined in terms of pairs $(A,B)$, where $A$ is a connection
for a fixed principal bundle $P(M,G)$ and $B$ is a 2-form on $M$. 
The relevant curvatures, parallel transports and holonomies
are computed and their expressions in local coordinates are  exhibited.
When the 2-form $B$ is given by the
curvature of $A$,  then the so-called non-abelian 
Stokes formula follows.
\vskip0pt
For a generic 2-form $B$,
we distinguish the cases when the parallel transport depends on the whole path
of paths  and when it depends only on the spanned surface. 
In particular we discuss generalizations of 
the non-abelian Stokes formula.
We study also the invariance properties of the (trace of the) holonomy 
under suitable transformation groups acting on the pairs $(A,B)$. 
\vskip0pt
In this way we are able to define observables 
for both topological and non-topological
quantum field theories of the $BF$ type. In the non topological case,
the surface terms may be 
relevant for the understanding of the quark-confinement problem. In the 
topological case the (perturbative) four-dimensional 
quantum $BF$-theory is expected to yield 
invariants of imbedded (or immersed) 
surfaces in a 4-manifold $M$.}
\vfill\eject} 
\pageno=1
\sezione{Introduction}
In this paper we consider the spaces $\loop$
and $\path$ of free
loops and paths of a compact manifold $M$ and the principal 
$G$-bundles on $\loop$ and $\path$
obtained by pulling back, via  the evaluation map, a fixed 
principal bundle
$P(M,G)$. 
We are interested in the connections  on such bundles
that are 
determined by pairs $(A,B)$ where
$A$ is a connection on $P(M,G)$ and  $B$ is a 2-form
of the adjoint 
type on $P(M,G)$. We study the properties
of the curvature and of the holonomy of such connections.

The motivations for this study  are rooted
in the four-dimensional 
quantum field theories of the $BF$-type. One of our goals
is to
understand the relation between those QFT's and the (smooth)
invariants of four-manifolds
and of surfaces imbedded  (or immersed)  in four-manifolds.

Before discussing the differential geometrical
results, we   comment briefly 
on quantum $BF$-theories.

\sottosezione{Quantum Field Theory}
Four-dimensional $BF$-theories
may become increasingly
relevant both to the quantum-field 
theoretical description of  smooth four-manifold invariants and to the
understanding of quark-confinement problems.

What characterizes $BF$-theories, and distinguishes them from
ordinary gauge theories, is the fact that there are two fundamental
fields: a connection $A$ for some principal
$G$-bundle over a four-dimensional manifold $M$
and a 2-form field $B$ that
transforms under gauge transformations as the curvature of $A$.

Various actions (and observables) 
can be constructed with the two fields $A$ and
$B$ and, as
a result, $BF$-theories can be both topological 
and non-topological (see section {\bf 8}). 

One of the relevant non-topological $BF$ theories is the first-order 
formalism of the Yang-Mills theory introduced in
\quot\halpuno{M.B.\, Halpern, {\it
Field Strength Formulation of Quantum
Chromodynamics,} \pr{D 16}, 1798--1801 (1977);}, \quot\halpdue{M.B.\,Halpern {\it 
Gauge Invariant
Formulation of the Selfdual Sector,} \pr{D 16}, 3515--3519 (1977);}, \quot\rein{
H. Reinhardt, {\it Dual Description of QCD,}
hep-th/9608191;
}\;and then modified in
\quot\bfsette{A.S.\,Cattaneo, \hskip0.1in P.\,Cotta-Ramusino, F.\,Fucito, 
M.\,Martellini, M.\,Rinaldi, A.\,Tanzini,   
M.\,Zeni,
 {\it  Four-Dimensional Yang--Mills Theory   as a Deformation of
Topological $BF$ Theory},  
to be published in \cmp{} (1998);}\; by replacing $B$ with
$B-d_A\eta$, where the  extra field $\eta$ is a 1-form.  The 
resulting theory, which
turns out to be  a deformation
of a topological theory, has been shown in 
\bfsette\ to be equivalent to the Yang-Mills theory.  

For general topological field theories the reader is referred to
\quot\witten{E.\,Witten, {\it Topological Quantum Field Theory},
\cmp{117}, 353--386 (1988);}.

Topological $BF$-theories have been introduced in \quot\ass{A.S.\,Schwartz,
{\it The Partition Function of a Degenerate Quadratic Functional and Ray-Singer
Invariants}, Lett.\, Math.\, Phys. {\bf 2}, 247--252 (1978);}\,
(see also \quot\hor{G.T.\,Horowitz, {\it Exactly Soluble Diffeomorphism 
Invariant Theories}, \cmp{125}, 417--436 (1989);})
and reviewed in
\quot\blau{D.\,Birmingham, M.\,Blau, M.\,Rakowski, G.\,Thompson,
{\it Topological Field Theories}, Phys. Rep. {\bf 209}, 129-340 (1991); 
}. The inclusion of observables in four-dimensional topological
$BF$-theories  is due to \quot\cm{P.\,Cotta-Ramusino, M.\,Martellini, {\it BF
theories and 2-knots}, in ``Knots and Quantum Gravity'', 
edited by J.\,Baez (Oxford
University, Oxford), 169--189 (1994);}. 

Here we begin a study of the {\it geometry of four-dimensional 
$BF$-theories}. The starting point is the observation that the
spatial components of $B$ are
the conjugate momenta to the spatial components of $A$.
If we formally identify the tangent and cotangent
spaces, we may see $B$ as
an infinitesimal connection.

The natural way to interpret the 2-form $B$  as a tangent
vector to a space of connections is to consider a principal G-bundle 
which has both as the total space and as the 
base manifold a space of loops or paths. 
By integrating the 2-form $B$ over a path, one obtains a 1-form.
We show in this paper
that  pairs $(A,B)$ represent  connections on a $G$-bundle over
the path or loop space of our manifold $M$.

An ordinary connection $A$ yields  the parallel transport along a path
of geometrical objects associated
to points. Similarly a connection on the path or loop space, 
represented by a pair $(A,B),$
yields the parallel transport along a path of paths or a path of loops
of geometrical objects associated
to paths or loops.

Surfaces are spanned by paths of paths (though not in a unique way), so 
the first question one has to ask is when a Stokes-like formula holds
for parallel transport. It is easy to see that
when the field $B$ is the (opposite of the) curvature
$F_A$ of $A$, then the parallel transport along a path of paths 
spanning a surface $S$
is uniquely determined by the $A$-parallel 
transport along the boundary of $S$. This is a new version 
of
the well-known
non-abelian Stokes formula (see \quot\halptre{M.B.\, Halpern, {\it Field
Strength and dual variable formulation of Gauge theory},
\pr{D 19}, 517--530 (1979);},  
\quot\arafe{I.Ya.\,Aref'eva, {\it Non-Abelian Stokes
Formula}, Teor.\ Math.\ Fiz. {\bf 43}, 111--116 (1980);}).
 
If $B$ is a small
deformation of a (small) curvature $F_A$, 
then a surface term appears in the parallel
transport with respect to the connection $(A,B)$. For large deformations 
the parallel transport with respect to $(A,B)$ depends on the whole path of 
paths structure and not only on the spanned surface.

We now recall that the ordinary holonomy along a loop in space--time
has a physical interpretation: namely, it represents the contribution
of a pair quark--antiquark forced to move along the loop. 
With our construction we have at our disposal more general objects.

For example we may consider the parallel transport 
with respect to the pair $(A,B)$ along
a path of paths. 
This provides us with a simple generalization of the previous
situation: we might think of this case as of
a pair quark--antiquark with an interaction that is not given
just by the $A$-parallel transport. 

As a second example we may take the holonomy 
corresponding to the pair $(A,B)$ of a loop in the path
(or loop) space. This should represent
the contribution of a pair of open (or closed) strings in interaction.

Thus, the equivalence between the
Yang--Mills theory and a deformation of the $BF$-theory plus the presence
of surface terms associated to the parallel transport involving a non-trivial 
$B$-field may be relevant for the problem of quark
confinement as formulated by 
Wilson\quot\wilson{Wilson, {\it Confinement of Quarks}, Phys.\,Rev.\,D,
{\bf 10}, 2445--2459 (1974);}. 

A  r\^ole of the $BF$-theories in the understanding 
of quark confinement has been
discussed in \quot\fmz{F.\,Fucito, M.\,Martellini, M.\,Zeni, 
{\it The $BF$ Formalism 
for QCD and Quark Confinement}, 
Nucl.\ Phys.\ {\bf B  496},  259--284 (1997);
} \,(see also \quot\kondo{K.\,Kondo, {\it
Yang-Mills Theory as a Deformation of Topological Field Theory, 
Dimensional Reduction and Quark Confinement}, hep-th/9801024;
}).

As far as topological $BF$ theories are concerned, we notice that 
one can define, at least in principle, a four-dimensional  analogue of the
Witten--Chern--Simons theory, whose vacuum expectation values (v.e.v.) of
(products of) Wilson loops represent link-invariants.

In four-dimensional topological $BF$-theories,
one can compute v.e.v.'s  of traces of holonomies of imbedded 
(or immersed) loops of paths (or of loops).
Again one of the delicate points is to check when
these invariants can be considered 
invariants of the surfaces spanned by  loops of paths (or of
loops). This requirement is close to the
parameterization invariance for the surface
since (for imbeddings) the loop-of-paths structure yields coordinates
on the surface. 

In this paper we show that when the fields $A$ and $B$,
restricted to the surface spanned by a loop of  paths, 
take values in an abelian
subgroup $T$ of $G,$ 
then the trace of the $(A,B)$-holonomy depends only on the surface and not 
on the loop-of-paths structure.

This reducibility condition is 
related to the abelian projection considered 
in \quot\hoof{G.`t\,Hooft, {\it On the Phase Transition towards
Permanent Quark Confinement,} \np{B 138}, 1--25 (1978);},
\quot\hooft{G.`t\,Hooft,  
{\it A Property of Electric and Magnetic Flux in Nonabelian Gauge
Theories,} \np{B 153}, 141-160 (1979);}.

The actual computation of v.e.v.'s for the topological $BF$
theories will be carried out elsewhere. There is some indication that
these v.e.v.'s may be related to 
the invariants of imbedded (or immersed)
surfaces considered by Kronheimer and 
Mrowka \quot\km{P.B.\,Kronheimer, T.S.\,Mrowka,
{\it Gauge theory for embedded surfaces, I,II }, 
Topology, {\bf 32}, 4, 773--826,
(1993)  and {\bf 34}, 1, 37--97 (1995);}.

Even though the original motivations 
for this paper lies in the development of $BF$
quantum field theories, our main goal here is to  study 
connections, curvatures and holonomies  of principal
fiber bundles over path and loop spaces. Our language will be therefore
the language of differential geometry.

\sottosezione{Geometry}

We begin  Sect. {\bf 2} by considering a fixed 
principal fiber bundle $P(M,G)$
together with  a connection $A$. The space of $A$-horizontal paths,
denoted
by $\pap$, is 
a principal $G$-bundle itself
which has the space of paths on $M$ as base manifold.

We describe explicitly the tangent bundle of $\pap$ as a submanifold
of the path space of $\T P$.
We consider connections on $\pap$ and particularly those connections
that are defined in terms of a 2-form $B$. We call these connections
on $\pap$ 
{\it special connections.}

The curvatures, horizontal
distributions and parallel transports corresponding  to special connections
are  computed. The explicit expression for 
the curvature involves Chen integrals. 

In Sect. {\bf 3} we discuss the parallel transport of paths of paths
with respect to special connections.
We single out the case when $B$ is given by the (opposite) of the curvature
$F_A$ of $A$: it is the only case
where we have an abelian Stokes formula, namely a relation
between  the parallel
transport 
of a path of paths and the $A$-holonomy of the loop given
by the boundary of such path of paths.  

More generally we study  the possible conditions that force  the (trace of the)
parallel transport along
a path of paths (or of loops) and the 
(trace of the) holonomy of a loop of paths (or of loops) to depend only on
the spanned surface.  The first of these conditions
is a ``perturbative" one: namely we assume that 
$P$ admits a flat connection and we expand  both $B$ and $F_A$
around zero. Then, for imbedded paths of paths, the (trace of the) 
parallel transport with respect to a special
connection $(A,B)$ depends only on the spanned surface, up to second-order 
terms.

In Sect. {\bf 4} we compute the expressions of the special connections
and of the parallel transport of paths of paths in local coordinates.

After recalling the general transformation properties of the holonomy
considered as a function of the space of connections for a generic
principal bundle (Sect. {\bf 5}), we discuss in Sect. {\bf 6}
the ``non-perturbative" conditions that guarantee that the holonomy of
a loop of paths is independent of those automorphisms of $P$
that maps the spanned surface into itself.

Here we require both $B$ and $A$ to be reducible to an abelian
subgroup of the structure group $G$ once they 
are restricted to the image
of the spanned surface.

In Sect. {\bf 7} we study  
the action on the space of pairs $(A,B)$
of those transformations groups  that happen to be 
symmetries for the $BF$-theories. 

The invariance of the (trace 
of the) holonomy under those transformations  
is guaranteed provided that
we ask
again the reducibility conditions for both $A$ and $B$. 

It is worthwhile noticing that the
group of gauge
transformations on $\pap$, which  preserves
the trace of the holonomy,  
is not a symmetry group for the $BF$-theories. More precisely
the symmetries of the $BF$ theories 
 are ``close" to being gauge transformations on $\pap$, 
the missing terms being boundary terms and 
higher-order Chen integrals.

The full group of gauge transformations on $\pap$, the space of
all connections on $\pap$ and the relation between $\pap$ and the
free loop bundle $\Lscr P$ (whose structure group is
the loop group of $G$)
will be discussed in a forthcoming paper\quot\finaledue{P.\,Cotta-Ramusino,
M.\,Rinaldi, {\it in preparation};}.
  
In section {\bf 8} we describe the observables for $BF$ quantum field
theories, both  in the topological and in the 
non-topological case.

\sezione{Differential geometry of horizontal paths} We describe
here the general setting of this paper. We 
consider a smooth manifold $M$
that is assumed to be closed, compact, oriented and Riemannian, 
a compact Lie Group $G$ with an $\Ad$-invariant inner product on its Lie algebra ${\frak g}$
and a fixed principal $G$-bundle $P=P(M,G)$ over $M$.  The group of
gauge transformations of $P$ will be denoted by $\Gscr$, while the
space of connections on $P$ will be denoted by $\Ascr$. 

Also we denote by $\form{*}$ the graded Lie algebra of forms on $M$
with values in the adjoint bundle $\adP=P\times_{\Ad}\lieg.$
We will consistently consider the elements of $\form{*}$ also 
as forms on $P$ that
are both of the adjoint type and tensorial \quot\kobanomi
{S.\,Kobayashi, K.\,Nomizu, {\it Foundations of Differential Geometry}, Vol. I,
Interscience Publishers, New York (1963);}.

The group $\Gscr$ acts on $\Ascr,$ and this action is free provided
that we restrict $\Ascr$ to be the space of irreducible connections
and divide $\Gscr$ by its center. We denote this action
as follows: $$\Ascr\times \Gscr \ni(A,g)\rsa A^g \in \Ascr.$$
 
In the course of this paper we will have to consider other principal
$G$-bundles, say $P_X(X,G)$, over some  manifold $X$, possibly
infinite-dimensional.  We will then denote then 
$\Gscr(P_X)$ and by
$\Ascr(P_X)$ the relevant group of gauge transformations and the 
space
of connections.  
If no confusion arises, we use the symbol $\pi$ to denote the projection
of any fiber bundle. 

For any manifold $X$ we denote by $\pathv{X}$ the space of
smooth paths on $X$. The space of  smooth 
free loops on $X$ will be  denoted by the symbol $\loopv{X}$ and the
space of $x$-based loops ($x\in X$) by the symbol $\loopvv{x}{X}.$ 
With some extra work we could consider also
piecewise smooth paths and loops, but we do not wish
to discuss this problem here. 

We will  also be interested in the space of smooth maps assigning to
each point $x\in X$ a path or a loop with initial point $x$. We call
such maps {\it path-fields} and, respectively, 
 {\it loop-fields}.  If we denote by $Map(X)$ the space of smooth
maps of $X$ to itself, then  a path field and a loop field on $X$
are  represented respectively by a path or a loop on $Map(X)$ with
initial point the identity map.

Most of this paper deals with horizontal lifts of paths on $M$
with respect to a
given connection $A\in\Ascr$. We use the following notation for
horizontal lifts;  for any $\gamma\colon[0,1]\to M$ and for any
$p\in P$ with $\pi(p)=\gamma(0)$, the $A$-horizontal lift of
$\gamma$ with initial point $p$ is denoted by the  symbol $$\partrap.
 \ref\lgothic$$ 
Our first
task is to study the differential geometry of the space of
$A$-horizontal paths.

\sottosezione{The principal bundle of horizontal paths and its
tangent bundle}  Let $\pap$ denote the space of $A$-horizontal
paths in $P$. This is a principal $G$-bundle $$\pap\riga{}{20}\path$$ 
where the right $G$-action is given by the right $G$-action on the initial
points of the horizontal paths. 

If we consider two distinct connections $A,\bar A\in \Ascr$,  then
we have two distinct and isomorphic principal $G$-bundles $\pap$ and
$\papv{\bar A}$. They are isomorphic since, for any connection $A$,
the bundle  $\pap$ is isomorphic to the pulled-back bundle
$ev_0^*P$. By the symbol
 $ev$ we denote in 
general the evaluation map and, in this particular
case, the map $$ev\colon \path\times I\to M,\quad I=[0,1], \quad
ev_t\triangleq ev(\cdot,t).\ref\evalone$$  Let us call $J_A$
the isomorphism between $ev^*_0 P$ and  $\pap$ given by
$$J_A(\gamma,p)\triangleq\partrap,\quad \gamma\in\path,\, p\in
\pi^{-1} {\gamma(0)} .\ref\isobundles$$  We denote by $\jsuba$ the
evaluation map  $ev^*_0 P\times I\to P$ given by 
$$\jsuba((\gamma,p),t)\equiv  {\frak
L}(A,\gamma,p)(t).\ref\isobundlesone$$

We have the following bundle morphisms:
$$\pap\times I \riga{ev}{30} P;\quad \pap\riga{ev_t}{30} P.$$
 
As a particular case we can  consider the loop space $\loop$ on
$M$,  instead of $\path$ and the corresponding principal bundle
$\lap$ whose elements are the $A$-horizontal {\it paths on $P$ whose
projections are loops}.

We now study the properties of the tangent bundle of the bundle $\pap.$

Firstly we identify $\T\path$ (the tangent bundle of the path space)
with $\pathv{(\T M)}$ (the path space of the tangent bundle).

In other words,  given any path $\gamma\in \path$, a vector $X\in
\T_{\gamma}(\path)$ is given by the assignment
for each $t\in I$ of a vector  $X(t)\in \T_{\gamma(t)}M$. 
Equivalently the same vector $X$ can be
represented by a smooth map $\Gamma\colon(-\epsilon,\epsilon)\times
I \to M$, so that $$\Gamma (0,t)=\gamma (t),\quad
\Gamma'(0,t)\equiv \dsopv {\Gamma(s,t)}= X(t).$$  

For any horizontal path $q$ on $P$, let us consider
a tangent vector ${\frak
q}\in \T_q(\pathv{P})$, defined by a smooth map
$Q\colon(-\epsilon,\epsilon)\times I \to P$ 
satisfying the following conditions
$$Q(0,t)=q(t), \quad \left(Q_* \pdt\right)_{Q(0,t)}=\dot q(t),\quad
\left(Q_* \pds\right)_{Q(0,t)}= {\frak q}(t).\ref\Qts$$
The tangent vector ${\frak q}$
belongs
to $\T_q(\pap)$ if and only if the following extra
requirement is satisfied: 
\blank 
$$\displaystyle{ {\left(Q_*\pds\right){ A ( {\dot Q(s,t)}
)}=0},\quad }  s=0,\;\forall t\in I.\ref\tangpap$$ 
%$$\displaystyle{ {\dsopv{ A ( {\dot Q(s,t)}
%)}=0},\quad }\forall t\in I.\ref\tangpap$$ 

Here we used the dot to denote the derivative with respect to the variable
$t\in I$. We will use this notation also in 
the future. Moreover when
dealing with two variables $(s,t)\in I\times I$
we will use the  prime to denote the derivative with respect to the
variable $s\in I.$

Condition \tangpap\ is independent of the choice of the map $Q(s,t)$
representing ${\frak q}$. In fact any two such choices $Q(s,t)$ and
$\widetilde Q(s,t)$ would satisfy (in local coordinates) the
condition $Q(s,t)-\widetilde Q(s,t)=s g(s,t)$ for some map $g$ with
$g(0,t)=0$. 

By considering the Lie derivative $L$ and inner product $i$ operators,
condition \tangpap\ is written as
$$\dy{L_{\pds}i_{\pdt} Q^*A=0}, \quad s=0, \forall t\in I.\ref\tangpapi$$

An important consequence of \tangpapi\ which will be used 
several times in the rest of this paper is given by
the following
\theorem\qgothic{For any element ${\frak
q}\in \T_q( \pap)$, the following equations hold: $$\gathered \dtv
{A\left({\frak q} (t)\right)}+ F_A\left({\frak q}(t),\dot
q(t)\right)=0,\quad \forall t\in I,\\
  A\left({\frak q} (t)\right)-A({\frak q}(0))= -\int_0^t \d t_1
\,F_A\left({\frak q}(t_1),\dot q(t_1)\right),\quad \forall t\in
I,\endgathered\ref\relaf$$  where $F_A$ denotes the curvature of
$A$.} 

\proof\qgothic \blank

Condition \tangpapi\ and the commutation property
$\dy{\left[\pds,\pdt\right]=0}$
imply
$$(Q^*dA)\left(\pdt,\pds\right)=L_{\pdt}i_{\pds} Q^*A,\quad s=0,
t\in I.
\ref\tangpapbis$$ 
Since $q$ is horizontal we have also 
$$(Q^*[A,A])\left(\pdt,\pds\right)=0,\quad s=0, t\in I.\ref\commutlocal$$
Equations \tangpapbis\ and \commutlocal\ and the structure equation for the 
curvature imply
$$(Q^*F_A)\left(\pdt,\pds\right)= L_{\pdt}i_{\pds}Q^*A,
\quad s=0, t\in I.\ref\fa$$
If we recall \Qts\ then \fa\ becomes immediately
the first equation of \relaf, namely
$$F_A\left(\dot q(t),{\frak
q}(t)\right)=\dtv
{A\left({\frak q} (t)\right)},$$ while the second equation of \relaf\
is obtained by integrating the
first.\endproof 

To any path $(q,{\frak q})$ in $\T P$ and $t\in I$
we associate  the tangent
vectors $(\dot q(t),\dot {\frak q}(t))\in \T_{(q(t),{\frak q}(t))} \TP$. 
Altogether the quadruple
 $\left(q,{\frak q},\dot q, \dot {\frak q}\right)$ represents
a path in $\T \T P.$ 

To any connection $A\in\Ascr$ we can canonically associate
a connection $\dota$ on the $\T G$-bundle $\T P$
\quot\kob {S.\,Kobayashi, {\it Theory of Connections}, Ann. Mat. Pura Appl. {\bf
43},  119--194 (1967);},
\quot\ccr{A.S.\,Cattaneo, P.\,Cotta-Ramusino, 
M.\,Rinaldi, {\it BRST symmetries for the
tangent gauge group}, \jmp{  39} 1316--1339 (1998);}, 
called the
tangential connection.

We recall that the tangential connection $\dota$ applied to an
element $\T\TP$ represented by a smooth map $Q\colon (-\epsilon,
\epsilon) \times (t_0-\epsilon,t_0+\epsilon)\to P$ yields 

$$ \left(A ({\dot Q} )(0,t_0) , { \dsopv{ A({\dot Q}
)(s,t_0)}}\right)\in {\frak g}\times {\frak g}.\ref\tangconn$$ 
\blank
We have then the following \remark\tangential
{A path $(q, {\frak q})$
in $\T P$ represents an element of   $\T \pap$ 
if and only if 
it is a $\dota$-horizontal path in $\T P$, 
namely if we have $\dota\left(q,{\frak q},\dot q, 
\dot {\frak q}\right)=0.$}
\blank
In other words $(q, {\frak q})$ is the $\dota$-horizontal lift
of a path $(\gamma, \rho)$ in $\T M$ 
with initial point $(q(0),{\frak q}(0))\in \pi^{-1}(\gamma(0),\rho(0))$.
\blank
A {\it vertical} vector ${\frak q}\in \T_q(\pap)$ is required to
satisfy the extra condition $${\frak q}(t)\in \V_{q(t)}P, \quad
\forall t\in I,\ref\verttang$$ where $\V_pP$ denotes the vertical
subspace of $\T_pP.$ 

Finally we have the following
\cor\vertical{As a consequence of   \qgothic\ and condition
\verttang, vertical vectors in $\T_q(\pap)$ satisfy  the
equation $$\displaystyle{\dtv{ A({\frak q}(t))}=0, \; \forall t\in
I}.$$} 

\sottosezione{Connections and curvatures on the bundles of
horizontal paths}

We now consider connections on $\pap$. In particular 
we are interested here in those connections on $\pap$ that are
determined by 2-forms in $\form{2},$ as shown in the following 
\theorem\connessione{Let $B\in \form{2}$ and $A,\bar A$ be any pair
of connections on $P$. The form  $$ 
 ev^*_0\bar A  + \int_Iev^* B  \ref\conne$$ defines a
connection on $\pap.$} 

\proof\connessione \blank The ${\frak g}$-valued 1-form $ev^*_0\bar
A$ is a connection on $\pap$. Moreover the 1-form
$\displaystyle{\int_I ev^*B}$ is of the adjoint type and is tensorial, as
can be seen by inspecting the explicit expression: $$\left(\int_I ev^*
B\right)(\frak q) = \int_0^1 {\text d}t\; B\left({\frak q}(t),\dot
q(t)\right),\quad {\frak q}\in \T_{q}(\pap). \ref\bigamma$$ 
 \endproof  We will call any connection of the above form a {\it special 
connection} on $\pap$ and will  denote it by the triple $(A,\bar
A,B)$. The space of special connections on $\pap$ is an affine space
modeled on $\form{1}\oplus\form{2}$.\blank
The reason why we are particularly interested in special connections is that 
elements of $\form{2}$ and connections
are the essential ingredients of
four-dimensional $BF$-theories 
\blau, \cm,
\quot\bfquattro{A.S.\,Cattaneo, P.\,Cotta-Ramusino, A.\,Gamba, 
M.\,Martellini, {\it The Donaldson-Witten Invariants and Pure QCD 
with Order and Disorder
't Hooft-like Operators}, Phys.\ Lett.\ {\bf B  355},  245--254 (1995);}, \bfsette.

The space of special connections is a proper subspace
of the space of all smooth connections on $\pap$. A simple
example of a connection on $\pap$ that is not special is 
very easy to construct. Let $t\rsa\wtilde B_t$ be a path in 
$\form{2}$ and $A,\bar A$ be any
pair of connections on $P$.  Then  we have  a connection on
$\pap$, defined on  ${\frak q}\in \T_{q}(\pap)$ as $$\bar
A(\frak q(0))+\int_0^1 \d t \,\wtilde B_t\left({\frak q}(t),\dot
q(t)\right).\quad  \ref\bigammaone$$ 
\blank
Other examples of connections 
on $\pap$ that are not special will be 
discussed extensively in a subsequent paper \finaledue.

In \conne\ we  often choose $\bar A=A$ and denote the triple
$(A,A,B)$ simply by a pair $(A,B)$. Here $A$ is kept fixed, so
the space of special connections on $\pap$ that are represented by
pairs $(A,B)$ is an affine space modeled on $\form{2}.$

\blank A vector ${\frak q}\in \T_q(\pap)$ is, by definition, {\it
horizontal} with respect to the connection $(A,\bar A,B)$ if the
following condition is satisfied: 

$$\bar A\left({\frak q} (0)\right) + \int_0^1 \d t \left[
B\left({\frak q}(t), \dot q(t)\right)\right]=0.\ref\horzpap$$ \blank

We consider two particular connections on $\pap$: \blank
\item{{\bf 1}} The {\it trivial} connection $(A,0)$ with curvature
$ev^*_0F_A.$ Here condition \horzpap\ is equivalent to requiring
that ${\frak q}(0)$ is $A$-horizontal. \blank \item{{\bf 2}} The
{\it tautological} connection $(A,-F_A).$ \vskip 1truecm
As a consequence of \qgothic, we have the following
\cor\tauto{The tautological
connection is given by $ev^*_1 A$ and its curvature is given
by $ev^*_1 F_A$. Condition \horzpap\ for the tautological connection
is the requirement that ${\frak q}(1)$ is $A$-horizontal.} \blank 
Let us add that on $\lap$ the tautological and the trivial connections
are {\it gauge equivalent}, but we refer to
\finaledue\  for the study of the gauge group
of $\lap$.

The 
computation of the curvature for a generic 2-form $B$
involves Chen integrals. These integrals
are defined in\quot\Chen{K.\,Chen, {\it Iterated integrals of differential forms and
loop space homology}, Ann. of Math. {\bf 97}, 2, 217--246 (1973);}\,\, for scalar
forms, but their extensions to forms in
$\form{*}$ is relatively easy and will be 
discussed extensively in \finaledue.

Here it is enough to say that the Chen integral
of a form $w\in \form{\deg (w)}$ is just the form given 
by the ordinary integral$$
\int_{\text {\rm Chen}} w\triangleq\int_I ev^*w \in \Omega^{\deg(w)-1}(\path, \ad(\pap)),$$
while the Chen bracket of two forms $w_1,w_2 \in \form{*}$
is the form in $\Omega^*(\path, \ad(\pap))$ of degree $\deg(w_1)+\deg(w_2)-1$
defined as 
$$\int_{\text{\rm Chen}}\{w_1; w_2\}\triangleq
\int_{0<t_1<t_2<1} \!    \left[w_1(\cdots,\dot\gamma(t_1))\d t_1,
w_2(\cdots, \dot
\gamma(t_2)  )\right]\d t_2 =$$ $$ = (-1)^{\deg( w_2)-1}
\int_{0<t_1<t_2<1} \!   \left[ w_1(\cdots,\dot\gamma(t_1)), w_2(\cdots,
\dot
\gamma(t_2)  )\right]\d t_1\,\d t_2,$$
where for each value of $t_i$
$w_i(\cdots,\dot\gamma(t_i))$ are $(\deg(w_i)-1)$-forms on $P$
to
be evaluated 
at tangent vectors in $\T_{\gamma(t_i)}P$.
Notice that the Chen bracket is bilinear, but
neither skew-symmetric nor
graded-skew-symmetric.
 
We have  the following \theorem\curvature{
The curvature $F_{(A,B)}$ of $(A,B)$ is given by the following
2-form on $\pap:$ $$\gather{ ev^*_0F_A - ev^*_1B + ev^*_0B
+\int_Iev^* d_AB  + (1/2)\left[\int_Iev^* B,\int_I ev^*B\right] -
\int_I \left[ev^* A -ev^*_0A, ev^* B\right]\\=ev^*_0F_A - ev^*_1B +
ev^*_0B +\int_I ev^* d_AB   + \chen \left\{ B  +  F_A ;  
B\right\}.}\endgather $$}    
\proof\curvature \blank  The
curvature of the connection $(A,B)$ is given by: $$F_{(A,B)}=ev^*_0
F_A + d_{ev^*_0 A}\int_I ev^* B + (1/2)\left[\int_I ev^*B,\int_I
ev^*B\right].$$

We then recall that we have the following relation between exterior
derivatives $$d\big|_{ \pap\times I}= d\big|_{\pap} \pm d\big|_I,$$
where the sign is given by the parity of the order of the form on
$\pap$. 

Hence we have the following chain of identities: $$\gather{ \int_I
ev^* d_AB = \int_I d_{ev^*A}ev^*B =\int_I d\,ev^*B
+\left[ev^*A,ev^*B\right]\\= d\int_I ev^*B +\int_I
\left[ev^*A,ev^*B\right]+ ev^*_1B -ev^*_0B\\= d_{ev^*_0 A}\int_I
ev^*B  +\int_I\left[ev^* A- ev^*_0A,ev^* B\right]+ev^*_1B -ev^*_0B.
}\endgather $$ 

We also have   $$(1/2)\left[\int_I ev^*B,\int_I
ev^*B\right] =\chen \left\{B; B\right\}$$  and  by taking into
account \qgothic, we  conclude the proof. \endproof
\blank
It is now natural to look for flat connections on $\pap$. If we restrict to
special connections $(A,B)$, then $(A,0)$ is a flat connection if $A$ is flat.
In order to find other flat connections, we have to require
some reducibility conditions.
\blank

Let $T$ be an abelian subgroup of $G$. We use the following
\definizione\redB{We say that a form $\omega\in \form{*}$ is reducible to $T$ 
if there exists a
$T$-subbundle of $P$, such that $\omega$ restricted to it takes values in 
$Lie(T)$.}
\blank
When we
require the reducibility of the connection $A$ and of some
forms $\omega_i\in \form{*}$, 
it will be understood that there will exist a $T$-subbundle of $P$
where the above forms are reducible {\it simultaneously}.
When we restrict ourselves to considering the bundle $\lap$
of horizontal paths whose projections are loops, then
a sufficient condition for the flatness of $(A,B)$ is given by the following
\theorem\flatness{The curvature of 
the connection $(A,B)$ on $\lap$ is zero if the
following conditions are satisfied:
\item{1.} $F_A=0$,
\item{2.} $d_AB=0$,
\item{3.} $A$ and $B$ are reducible to $T$.}
\blank
To conclude this
section, we  recall  that the bundles $ev_0^*P$ and $\pap$ are
isomorphic via \isobundles. We denote by $\connaab$ and $\connab$
the connections on
$ev_0^*P$ induced respectively 
by the connection $(A,\bar A, B)$ and $(A,B)$ on $\pap.$

Namely we
set  
$$\connaab \triangleq ev_0^*{\bar A} +\int_I \jsuba^*B,\ref\connectaab$$ 
$$\connab \triangleq ev_0^*A +\int_I \jsuba^*B,\ref\connectab$$ 
where
$\jsuba$ has been defined in \isobundlesone.  The curvature of
$\connectaab$ and $\connectab$ will be denoted respectively
by the symbols $\curvaab$ and 
$\curvab$.

\sezione{Horizontal
lift of paths of paths and the non-abelian Stokes formula}  
In this section we consider the parallel transport of paths of paths
and 
in particular of imbedded paths of paths. We discuss when
the relevant parallel transport is invariant under isotopy. 

Any 
connection on $\pap$ defines  horizontal lifts of a
path of paths $\Gamma$ in $M$, namely, of a map
$$\Gamma\colon [0,1]\times [0,1]\to M.$$  Each of these horizontal lifts
depends on the choice of the {\it initial path} $q\in\pap$, with
$\pi(q(t))=\Gamma(0,t).$ In turn this initial path $q$, being
$A$-horizontal, depends only on the choice of a {\it
initial point} $q(0)\in \pi^{-1}\Gamma(0,0)$. So we will speak of
{\it horizontal lift of paths of paths with respect to an initial 
point} $p_0\in \pi^{-1} \Gamma(0,0) .$

The horizontal lift of a path of paths $\Gamma$ with respect
to a given connection on $\pap$ is, by definition, a
path of  $A$-horizontal paths.  So 
we have the following:
\theorem\lif{The horizontal lift of $\Gamma\colon I\times I\to M,$
with respect to a given connection on $\pap$ and an initial point
$p_0\in \pi^{-1}\Gamma(0,0),$
is  uniquely determined
by the lift of the path of the initial points of the given paths
$s\rsa\Gamma(s,0)\in M$.} 
\blank
Notice that the lift of the path of initial points 
considered in \lif\ coincides with the horizontal lift
with respect to a connection $\bar A\in \Ascr$ 
only if we choose the special connection $(A,\bar A,0)$
on $\pap$. For a general connection
on $\pap,$ the lift of the path of initial points is more general 
than the horizontal lift: it is still $G$-equivariant but  
depends on the whole
path of paths $\Gamma.$

It is therefore convenient to consider
the following general
definition of path-lifting for a principal bundle $P(M,G)$.

\definizione\lifting{A lift is a smooth map:
${\frak h}\colon ev_0^*P\to\pathv{P}$ satisfying   the
conditions
$${\frak h}(p,\gamma)(0)=p;
\quad \pi\left({\frak h}(p,\gamma)\right)=\gamma,$$
and the $G$-equivariance 
$${\frak h}(ph,\gamma)= [{\frak h}(p,\gamma)]h,\quad\forall h \in G.$$}

Our definition of lift is a smooth 
$G$-equivariant version of the definition
of a ``connection" for  the
{\it fibration} $\pi\colon P\to M$, 
as given, e.g., in \quot\white{G.W.\,Whitehead, {\it Elements
of Homotopy theory}, 
Springer Verlag, Berlin Heidelberg, New York (1979);}. 
But we will use the term ``lift" instead of ``connection"
in order 
to avoid confusion with the ordinary connections on $P(M,G)$.

We recall that a path-field is a smooth map $M\to \path$
that assigns to each $x\in M$ a path beginning at $x$.
Each path-field $Z$ composed with $\gamma$ yields an element  
of $\pathv{(\path)}$.
When we lift $Z\circ\gamma$
via a connection on $\pap$, we obtain a path of paths in $P$
whose initial points
are a lift of $\gamma$ in 
the sense
of \lifting. Hence \lif\ can be rephrased as follows: 
\theorem\pathfield{If we denote by ${\frak P}(M)$
the space of path-fields on $M$, by 
${\frak H}(P)$ the space of lifts as  
in \lifting\ and by $\Ascr(ev^*_0P)$ the 
space of connections on $ev^*_0P$, then we have a map:
$${\frak P}(M) \times \Ascr(ev^*_0P)\to {\frak H}(P).\ref\sgothic$$}

In particular standard horizontal lifts correspond either 
to the choice of a special connection of the type
$(A,\bar A,0)$ on $\pap$ together with an
arbitrary choice of a path-field or to an arbitrary choice
of a connection on $\pap$ together with the choice of the trivial 
path-field (i.e. the path field
assigning the constant path to every point of $M$). This shows in 
particular that \sgothic\ is far from being injective. 
\blankm
 
Following definition \lgothic\ we denote the horizontal lift with 
respect to the connection $\omega$ on $\pap$ 
by the symbol $\hliftomo$. 
The comparison 
of the lift of the initial
points of $\Gamma\in\pathv{(\path)}$ with respect
to the two  connections $\omega$  and  $(A,0)$
defines a path $k_{\Gamma,\omega} \in \pathv{G}$  such that 
$$\hliftomo(s,0)=\hliftaoo(s,0)\cdot k_{\Gamma,\omega}(s).\ref\liflift$$

Due to the 
$A$-horizontality of the
the lifted paths $t\rsa \hliftomo(s,t),$
equation \liflift\ holds also for a generic $t$, namely we have
$$\hliftomo(s,t)=\hliftaoo(s,t)\cdot k_{\Gamma,\omega}(s).\ref\kgammadef$$
 
When we work with a fixed path of paths $\Gamma$ and a fixed base
point $p_0\in\pi^{-1} \Gamma(0,0) $ we  use a
simplified notation, i.e., we  set  $${\frak L}_A(s,t)\equiv
\hliftaoo(s,t).
 \ref\pstbis $$
 
From \liflift\ we conclude that $k_{\Gamma,\omega}$
satisfies the following differential equation:
$${dk_{\Gamma,\omega}(s)\over
ds}k_{\Gamma,\omega}(s)^{-1}= -\omega ({\frak L}'_A(s,\bullet)),
\ref\psttter$$
where we have to keep in mind that for any $s\in I$,
the map 
$$ {\frak L}'_A(s,\bullet)(t)\equiv {\frak L}'_A(s,t)$$ represents 
a tangent vector in $\T\pap.$

We will use consistently in this paper the notation of 
\psttter. Namely for any function
$f$ of several variables $a,b,c,d,\cdots$
we denote by $f(\bullet,b,c,d,\ldots)$   the function of one
variable ($a$) obtained by evaluating $f$ at $b,c,d,\ldots$).

When we choose 
$\omega$ to be the special connection $(A,A+\eta,B),$
then equations \kgammadef\ and \horzpap\
imply the following
differential equation: $${dk_{\Gamma,A,\eta,B}(s)\over
ds}k_{\Gamma,A,\eta,B}(s)^{-1}= -\int_0^1 B({\frak L}'_A(s,t),\dot {\frak
L}_A(s,t))\d t -\eta \left({\frak L}'_A(s,0)\right).\ref\pstter$$ 
The solution is  given by a
path-ordered
exponential 
(in the variable $s$)  
$$k_{\Gamma,A,\eta,B}(s)=\Pscr\exp\left\{-\int_{[0,s]}
\d s_1\left[ \int_0^1 B({\frak L}'_A(s_1,t),\dot {\frak
L}_A(s_1,t))\d t +\eta \left({\frak L}'_A(s_1,0)\right)\right]
\right\}.
\ref\kgammasol$$ If $G$ is an
abelian group (e.g., $U(1)^n$), then path-ordering is
not needed.
\blank
Consider now the evaluation map $ev\colon \pathv{(\path)}\times I\times
I\to M$ and the pulled-back bundle $ev_{0,0}^*P$ whose elements
are represented precisely 
by pairs $(\Gamma,p_0)$ where $\Gamma$ is a
path of paths on $M$ and $p_0\in P$ is an element in the fiber over
$\Gamma(0,0).$

We have the following
\theorem\holopath{ Any connection $\omega$ on $\pap$ determines a map
$$\Hscr_{\omega}\colon ev_{0,0}^*P\to G$$ of the adjoint type,
i.e. satisfying the equation
$$\Hscr_{\omega}(\Gamma, pg)={\Ad}_{g^{-1}}
\left(\Hscr_{\omega}(\Gamma,p)\right),\quad \forall g\in G.
\ref\Hscruno$$ In
particular when $\omega$ is a special connection $(A,A+\eta,B)$, then the
map $\Hscr_{(A,A+\eta,B)}$ has the following properties
$$\Hscr_{(A^{\psi},A^{\psi}+{\Ad}_{\psi^{-1}}\eta,
{\Ad}_{\psi^{-1}}B)}(\Gamma,p)={\Ad}_{\psi^{-1}(p)}
\left(\Hscr_{(A,A+\eta,B)}(\Gamma,p)\right), \quad \forall
\psi\in\Gscr.\ref\Hscrdue$$} \proof\holopath \blank 
We define $\Hscr_{\omega}(\Gamma,p)\equiv k_{\Gamma,\omega}(1)$, where
the r.h.s. is in turn defined by \liflift. \Hscrdue\ is a consequence of
\kgammasol.
\endproof
\blank
\holopath\ summarizes the properties of the 
{\it horizontal lift of paths of paths}. 
Let us now consider the ``square" associated to (i.e. the image 
of) a path of paths   $I\times I$ to $M$.
If we compute
the map $\Hscr$
applied to two {\it different paths of paths} with the {\it same image} 
in $M$, is the result the same?   
The answer 
is in general no, but some special situations are worth of
consideration.

First we consider the case of the trivial connection on $\pap$.
It follows from the definition,
that for any $\Gamma$ and for
any $A\in \Ascr$, we have $\Hscr_{A,0}(\Gamma,p_0)=1$.

\blank

Next we consider 
the {\it tautological connection} $(A,-F_A)$. 
For this connection  the {\it non-abelian Stokes
formula} holds, namely we have the following
\theorem\arafeva{For any path of paths
$\Gamma\in\pathv{(\path)}$ and for any $p_0\in
\pi^{-1} \Gamma(0,0) $ we have
$$\Hscr_{(A,-F_A)}(\Gamma,p_0)=\Hol_A(\partial\Gamma,p_0).$$}

Here $\partial\Gamma\colon [0,1]\to M$ denotes the (smooth) \foot {We
assume that, when needed, all corners are properly smoothed.} loop
defined by
$$
\halign{
\hskip80pt$ # $&$#\hfil$&\hskip20pt$#$&$#$\cr 
 (\partial\Gamma)(\tau)& =\Gamma(0,4\tau),\quad && 0\leq \tau\leq
\frac 14\cr (\partial\Gamma)(\tau)&=\Gamma(4\tau-1,1),\quad &&\frac 14  
\leq \tau\leq  \frac 12\cr (\partial\Gamma)(\tau)&=\Gamma(1,3-4\tau),
\quad &&\frac
12
\leq \tau\leq
\frac 34\cr (\partial\Gamma)(\tau)& =\Gamma(4-4\tau,0).\quad &&\frac 34 \leq
\tau\leq 1.
\cr}  $$

\proof\arafeva
\blank
Here the connection on
the bundle of horizontal paths, is given by $ev^*_1A$. 
Hence the $(A,-F_A)$-horizontal lift at $p_0\in P$ of the
path of initial points $\Gamma(\bullet,0)$ is  obtained as
follows. 

We  first consider the $A$-horizontal lift of $\Gamma(0,\bullet)$
and its end point $p_1=\Gamma(0,1)$. Then we consider the
$A$-horizontal lift of the path of end-points $\Gamma(\bullet,1)$
beginning at $p_1$ and the $A$-horizontal lift of all the paths
$\Gamma(s,\bullet)$ for $s\in (0,1]$ with assigned end-point.
The resulting path of
initial points is the $(A,-F_A)$-horizontal lift of
$\Gamma(\bullet,0)$. The theorem follows  immediately from 
\holopath.\endproof
\blank
The non-abelian Stokes Formula
has a long history, starting from \halptre, \arafe.  
For some more recent papers see \quot\broda{B.\,Broda, {\it Non-Abelian Stokes
Theorem},  in ``Advanced Electromagnetism: Foundations, Theory and 
Application'' (T.\,Barrett, D.\,Grimes eds.) World Scientific, Singapore, 
 496--505 (1995);}, \quot\diak{D.\,Diakonov, V.\,Petrov, {\it Non-Abelian Stokes 
Theorem and Quark-Monopole Interaction}, hep-th 9606104;}, \quot\lun{F.A.\,Lunev,
{\it Pure Bosonic Worldline Path Integral Representation for Fermionic 
Determinants, Non-Abelian Stokes Theorem, and Quasiclassical Approximation
in QCD}, \np{B 494}, 433--470 (1997);}.  The treatment of the problem as a problem
of parallel transport in a
space of paths is new. 

\blank

We now consider different paths of paths with the same image in $M$ and 
see if their image with respect to the
map given in 
\holopath\
is the same.

{\it In this section, from now on, 
we limit ourselves to considering imbedded paths (or loops)
of paths}. 
We assume,  in particular, that we
 have an {\it isotopy} $\Gamma_r\colon I\times I\to M,\;r\in[0,1]$
satisfying the following
assumptions:
\blank
\item{\bf G.1} $\Gamma_r(0,0)=\Gamma_0(0,0),\quad \forall r\in[0,1]$
\item{\bf G.2} $\im(\Gamma_r)=\im(\Gamma_0),\quad \forall r\in[0,1]$
\item{\bf G.3} $\im(\Gamma_r(\bullet,0))=\im(\Gamma_0(\bullet,0)),
\quad \forall r\in[0,1].$\blank
\blank
By taking derivatives of $\Gamma_r$ with respect to the parameter
$r$, we define for
each $r$ a smooth map  $Z_r\colon I\times I\to \TM$,
with $Z_r(s,t)\in
\T_{\Gamma_r(s,t)}M.$  The  
  following conditions
for $Z_r$ are a consequence of the corresponding 
conditions for $\Gamma_r$:
\blank
\item{\bf Z.1} $Z_r(0,0)=0$ 
\item{\bf Z.2} $Z_r(s,t)$ is tangent to $\im(\Gamma_0)$ and 
the restriction of $Z_r$ to $\partial(I\times I)$ 
is tangent to  $\im(\partial\Gamma_0)$.
\item{\bf Z.3} $Z_r(s,0)$ is tangent to $\im(\Gamma_0(\bullet,0))$
and $Z_r(1,0)=0$
\blank

When do we have $\Hscr(\Gamma_r,p_0)=\Hscr
(\Gamma_0,p_0)$? 
We have  the following partial answer
\theorem\smoothomo{If $\Gamma_r\colon I\times I\to M$ is an isotopy 
satisfying the conditions G.1 and  G.2 above and if, 
moreover, $F_A=0$,
then for any $B\in \form{2}$ we have
$$\Hscr_{A,\lambda B}(\Gamma_r,p_0)=
\Hscr_{A,\lambda B}(\Gamma_0,p_0)+o(\lambda),\quad \forall
r\in [0,1].\ref\homouno $$
If, in addition, condition G.3 is satisfied,
then we have also 
$$\Hscr_{(A,A+\lambda \eta,\lambda B)}(\Gamma_r,p_0)=
\Hscr_{(A,A+\lambda \eta,\lambda B)}(\Gamma_0,p_0) +o(\lambda),
\quad \forall r\in [0,1].
\ref\homodue$$}
\blank
Here we have assumed that some representation of the group
$G$ has been
chosen so that the sum in \homouno\ and \homodue\ makes sense.
Even though we do not have in general   
a true horizontal lift of squares or of surfaces, the implications
of \smoothomo\ are that in some particular cases
such horizontal lifts do exist.
This is true if  we consider imbeddings or immersions as paths of paths, 
and small deviations from
a flat connection and from $B=0$.
\proof\smoothomo
\blank
Consider $\hliftraoo$, i.e., the $(A,0)$-horizontal lift of
$\Gamma_r$. Since $A$ is flat, the image under
$\hliftraoo$ of any curve in
$I\times I$ is $A$-horizontal in $P$.
We take derivatives with respect to $r$ in $\hliftraoo$
and obtain a map $\bar Z_r\colon I\times I\to \T P$. For each 
$(s,t)$, $\bar Z_r(s,t)$ is now a horizontal lift of $Z_r(s,t).$
The map $r\rsa \Hscr_{A,\lambda B}(\Gamma_r,p_0)$ defines a curve 
in $G$.
By taking the logarithmic derivative of
the above map, we obtain  
an element of $\lieg.$  If this element is zero,
up to terms of order $\lambda^2$ and  for any $r$, then the theorem
is proved. 

We use again a simplified notation by setting
$${\frak L}_{A,r}(s,t)\equiv
\hliftraoo(s,t).
 \ref \pster$$ 
We first consider equation \homouno.
By taking into account \kgammasol, we see
that the element in $\lieg$ we are looking for
is
$$-\int_{I\times I}\ddr{\left({\frak L}_{A,r}^*B\right)}+
o(\lambda).\ref\logzero$$
The integrand above coincides with  
$${\frak L}_{A,r}^*\left(L_{\bar Z_r}B\right).$$
But the Stokes theorem and property $Z_2$ imply that the integral of 
${\frak L}_{A,r}^*di_{\bar Z_r}B$ vanishes.
Moreover we have 
$$\int_{I\times I}{\frak L}_{A,r}^*i_{\bar Z_r}dB
=\int_{I\times I}{\frak L}_{A,r}^*i_{\bar Z_r}d_AB=0$$ since,
for any $X,Y\in \T_{s,t}
(I\times I),$ the three vectors $\left({\frak L}_{A,r}\right)_*X,$
$\left({\frak L}_{A,r}\right)_*Y$, $\bar Z_r(s,t)$ are $A-$horizontal
and linearly dependent.

As for \homodue, we
set
$${\frak L}_{A,r,0}(s)\equiv
\hliftraoo(s,0).
 \ref \pstermine$$ 
In order to prove \homodue\ we have
to show the vanishing of the term
$$\int_{I}\ddr{\left({\frak L}_{A,r,0}^*\eta\right)}=
\int_I {\frak L}_{A,r,0}^* \left(L_{\bar Z_r}\eta\right) 
.\ref\logzerozero$$
The integral \logzerozero\ vanishes since the Stokes theorem and 
conditions 
Z.1-Z.2 imply that $\int_I {\frak L}_{A,r}^*(di_{\bar Z_r}\eta)$
vanishes, while condition Z.3 implies $\int_I
{\frak L}_{A,r}^*(i_{\bar Z_r}d_A\eta)=0.$
\endproof

\remark\subman{
If the image of $\Gamma_r$ is contained in a submanifold
$i\colon N\incul M$, then in order for the conclusions of \smoothomo\ to remain true,
it is enough to require $i^*F_A=0$. Moreover in \homouno\ and
\homodue\ we may replace $\lambda B$ with any $B(\lambda)$ such that
$i^*B=o(\lambda)$.} 
\blank

Let us come back to the non-abelian Stokes formula.
This formula implies that
$\tr\Hscr_{(A,-F_A)}(\Gamma,p_0)$ 
coincides with the Wilson loop of the boundary $\partial \Gamma$.

We recall that the Wilson loop is defined precisely as
$\tr \Hol_A(\gamma,p_0)$ for  $\gamma\in
\loop$, $A\in \Ascr$ and $\pi(p_0)=\gamma(0).$

When we consider instead of the tautological connection
a generic special connection $(A,B)$, the corresponding
generalized Wilson loop $\tr\Hscr_{(A,B)}(\Gamma,p_0)$ depends on
the path of paths $\Gamma$ and not only on $\partial\Gamma$.

This may be relevant for the understanding of
the {\it quark-confinement
problem} in the framework of $BF$-theories.

In particular we are
interested in considering  
generalized Wilson loops represented by
deformations
of  the ordinary Wilson loop, where,
up to the second  order
in the perturbative expansion, 
$\tr\Hscr_{(A,B)}(\Gamma,p_0)$ depends only on
the surface $\im(\Gamma)$
and not on the particular path of paths $\Gamma$.

Accordingly  we consider a special connection
given by a perturbation series in a neighborhood of 
a {\it flat connection} $(A,0)$, where $F_A=0$. 

We may use two different
variables $\kappa$ and $\lambda$ to describe 
respectively the deformation of the connection $A$ and of the
2-form field $B$, i.e. we set:
$$A(\kappa)\equiv
A+\kappa\eta +o(\kappa),\quad 
B(\lambda)\equiv\lambda B + o(\lambda).\ref\deformab$$

We now choose a smooth isotopy of imbeddings $\Gamma_r$
(or smooth homotopy of immersions) as before 
and  set 
 
$$\Hscr(\kappa,\lambda,r,p_0)\equiv
\Hscr_{(A(\kappa),-F_{A(\kappa)} +B(\lambda))}(\Gamma_r,p_0).
\ref\abr$$
If we have only one parameter $\kappa,$ we set 
$$\Hscr(\kappa,r,p_0)\equiv
\Hscr_{(A(\kappa),-F_{A(\kappa)} +B(\kappa))}(\Gamma_r,p_0).
\ref\abrr$$ 
\blank
When $\lambda=0$,
then {\it the non-abelian Stokes formula
implies that $\Hscr(\kappa,\lambda=0,r,p_0)$
is independent of $r$.}

In the general case the power series expansions
of 
$\Hscr(\kappa,\lambda,
r,p_0)$  and $\Hscr(\kappa,r,p_0)$
depend on $r$ but satisfy the following:

\theorem\wilson{Let $P(M,G)$ be a principal $G$-bundle
admitting a flat connection $A$.  Let $\Gamma_r$ satisfy G.1 and G.2
and let $A(\kappa)$ and $B(\lambda)$ be defined as above. 
For $\Hscr(\kappa,\lambda,r,p_0)$
given by \abr\ we have the following equation:
$$\left.{\partial^2 \Hscr(\kappa,\lambda,r,p_0) 
\over \partial \lambda \partial r}\right\vert_{\kappa=\lambda=0}=0.
\ref\quark
$$
If we have only one parameter $\kappa=\lambda$, and we assume also
G.3,
then we have
$$\left.{\partial^2 \Hscr(\kappa,r,p_0)
\over \partial \kappa \partial r}\right\vert_{\kappa=0}=0,
\ref\superquark$$
where definition \abrr\ has been assumed.}
\blank
\wilson\ provides a {\it surface law for the generalized
Wilson loop} in $BF$-theories. The
main difference between \smoothomo\
and \wilson\ lies in the fact that in the latter the field $B$
deforms a (non-trivial) tautological connection (for which
the non-abelian Stokes formula holds) at any order in $\kappa.$

\proof\wilson
\blank
As in \liflift\ and \pstbis\ we set 
$$\hliftabor(s,t)\equiv
\hliftraoo(s,t)k_{(r,\kappa,\lambda)}(s),$$ 
with $k_{(r,\kappa,\lambda)}\in \pathv{G}$ 
and
$${\frak L}_{(r,\kappa)}(s,t)\equiv \hliftraoo(s,t).$$
Analogously to \pstter\ we have
$$\eqalign{{dk_{(r,\kappa,\lambda)}(s)\over
ds}k_{r,\kappa,\lambda}(s)^{-1} &= -\kappa
\eta({\frak L}'_{(r,\kappa)}(s,0))
+\cr
\quad &\int_0^1 [F_{A(\kappa)}- B(\lambda)]
({\frak L}'_{(r,\kappa)}(s,t),\dot {\frak
L}_{(r,\kappa)}(s,t)) \d t.\cr}\ref\pstquater$$
By taking the derivative of \pstquater\ with respect to
$\lambda$ at $\kappa=\lambda=0$
the r.h.s. of 
\pstquater\ becomes
$$-\int_0^1 B \left({\frak L}_{(A,r)}'(s,t), \dot{\frak L}_{(A,r)}(s,t)
\right)
\d t,$$ (see \pster\ for the notation) 
and the proof of \smoothomo\ applies verbatim to our case.

As for \superquark\ we notice that we have to replace 
$\lambda$ with $\kappa$ in \pstquater.
The derivative with respect to $\kappa$ at $\kappa=0$
of the r.h.s of \pstquater\ becomes
$$-\eta\left({\frak L}_{(A,r)}'(s,0)\right)+ \int_0^1 [-B+d_A\eta] 
\left({\frak L}_{(A,r)}'(s,t), \dot{\frak L}_{(A,r)}(s,t)
\right)
\d t.$$
We differentiate again with respect to $r$. Using G.3 and the 
same arguments as in \smoothomo, we obtain  
\superquark.
\endproof
We end this section by considering the special
case of $\Gamma$ being an imbedded {\it loop} of paths. 
The 
holonomy with respect to the connection $(A,B)$
is then given by
$$\Hol_{(A,B)}(\Gamma,p_0)=  \Hol_A(
\Gamma(\bullet,0),p_0)\Pscr_s\exp\left(-\int_{ {\frak L}_A([0,1]\times
[0,1])}B\right). \ref\holAB$$ 

It is clear that  $\Hol_{A,B}(\Gamma,p_0)$  is given by the group
element $g=g\left(\Gamma, (A,B), p_0\right)$ such that $p_0g$ is the
end-point of the $(A,B)$-horizontal lift of the loop of initial
points $\Gamma(\bullet,0)$. If $B=0$, then the above holonomy is nothing
else than the $A$-holonomy of the loop of initial points. 

\remark\hscrhol{If $\Gamma\in\loopv{(\path)}$ then we have:
$$\Hscr_{(A,B)}(\Gamma,p_0)=
\Hol_A^{-1}(\Gamma,p_0)\Hol_{(A,B)}(\Gamma,p_0)$$and if $\Gamma$ is a loop of
loops:
$$\Hscr_{(A,-F_A)}(\Gamma,p_0)= \Hol_A^{-1} (\Gamma(\bullet,0),p_0)
\Hol_A (\Gamma(0,\bullet),p_0) \Hol_A (\Gamma(\bullet,0),p_0)
\Hol_A^{-1} (\Gamma(0,\bullet),p_0).$$ }

\sezione{Local
coordinates}

Here we discuss the expressions of the special connections on
$\pap$ and of
the relevant
parallel transport in {\it local coordinates}. 

Let $U$ be the domain of a local chart in $M$. We denote by
$\pathv{U}$ the space of paths in $U$ and by  $\Pscr_UM$ the space
of paths in $M$ with initial point in $U$.
Any section $\sigma\colon U\subset M\lora P$ determines a section
$$\eqalign{ \wsigma\colon& \Pscr_UM\to \pap, \cr \gamma\rsa&
\wsigma(\gamma)\triangleq {\frak L}(A,\gamma,\sigma\left(\gamma
(0))\right), \cr}\ref\widesigma$$  where, as before,  $\frak L$ 
 denotes the horizontal lift. So the bundle of horizontal
paths is trivial if and only if  the bundle $P$ is trivial. 
\definizione\defh{For any  section $\sigma\colon U \to P$  we define 
$h\colon\pathv{U} \times I\to G$ by the equation $$\sigma(\gamma(t))
h(\gamma,t)\equiv [\wsigma (\gamma)](t).$$  }  The  map $h$ 
allows us to compare the $A$-parallel transport with the
image of a section $\sigma$ and is given by the standard path-ordered
exponential of the integral 
$$h(\gamma,t)= \Pscr \exp \int_{[0,t]}\gamma^*(-\sigma^*A).$$
\sottosezione{Connections}

Given any $X\in \T_{\gamma}(\path)$, we set $h_\gamma(t)\equiv
h(\gamma,t)$, $q=\wsigma(\gamma)$ and ${\frak q}\equiv \wsigma_*
X$.  We have $${\frak q}= \sigma_* X\,h_\gamma +
i\left((h^{-1}dh)X\right), \ref\frq$$ where the map $i\colon {\frak g}\to
{\frak X} (P)$ is, by definition, the map yielding fundamental
vector fields.  Moreover we have $$\dot q = (\sigma_* \dot \gamma)  
h_\gamma + (\sigma\circ\gamma) \dot h_\gamma.\ref\qudot$$ The second
terms of both \frq\ and \qudot\ are vertical vectors fields along
$q$, so we finally obtain $$\wsigma^*\left(\int_I q^*B\right) (X)=
\int_0^1 {\text d}t\;\;\Adhmt\big(\sigma^*B\left(X(t),
\dot\gamma(t)\right)\big ).\ref\coordconn$$ This is the {\it
expression in local coordinates of the difference between the
connection $(A,B)$ and the connection} $(A,0).$\foot{The expression
\coordconn\ was firstly considered in \cm\ where the
notation $\Hol_A(\gamma)_0^t$ for $h_{\gamma}(t)$ was 
employed.}

\sottosezione{Horizontal lift of paths of paths} 

Now we compute the $(A,B)$-horizontal lift
of  $\Gamma\in \pathv{\left(\path\right)}$ in local coordinates.  

Consider a map
$\Gamma\colon I\times I \to M,$ where the first variable ($s$)
describes the path of paths, while the second variable ($t$)
describes each individual path. We assume that the image of $\Gamma$
is all contained in the domain $U$ of a local section
$\sigma\colon U\subset M\to P$. 

We consider the section $\wsigma$ ( \widesigma) on $\pap$. We
have explicitly, for each {\it fixed} $s\in I$,
$$(\wsigma\Gamma)(s,t)=\sigma\left(\Gamma(s,t)\right)
h_{\Gamma(s,\bullet)}(t),$$ where $h$ is as in  \defh\ and
$\Gamma(s,\bullet)$ denotes the path in $\path$ given by $t\rsa
\Gamma(s,t)$. 
 
The $(A,0)$-horizontal lift of $\Gamma$ is given, in local
coordinates, by $$(s,t)\rsa
(\wsigma\Gamma)(s,t)h_{\Gamma(\bullet,0)}(s)=
\sigma\left(\Gamma(s,t)\right)
h_{\Gamma(s,\bullet)}(t)h_{\Gamma(\bullet,0)}(s).$$

When we consider
as in \liflift\ the path $k_{\Gamma,A,B}\colon I\to
G$, then the $(A,B)$-horizontal lift of $\Gamma$ is given, in local
coordinates, by $$(s,t)\rsa (\wsigma \Gamma)(s,t)
h_{\Gamma(\bullet,0)}(s)k_{\Gamma,A,B}(s).\ref\dueh$$ 

We now consider the following vectors in $\T_{(\wsigma\Gamma)(s,t)}P:$
$$\Gamma_1(s,t)\triangleq \pdsv{\wsigma\Gamma}(s,t),\quad
\Gamma_2(s,t)\triangleq \pdtv{\wsigma\Gamma}(s,t).$$  We also set
$$K_{\Gamma,A,B}(s) \triangleq
h_{\Gamma(\bullet,0)}(s)k_{\Gamma,A,B}(s).\ref\kdis$$  
The $(A,B)$-horizontality of \dueh\
translates into the following equations
$$\dsv{K_{\Gamma,A,B}(s)}K_{\Gamma,A,B}^{-1}(s) +
A\left(\Gamma_1(s,0)\right) + \int_0^1{ \d t\;\;
B\left(\Gamma_1(s,t),\Gamma_2(s,t)\right)}=0,\ref\Bholonomy$$ 
 $$\dsv{h_{\Gamma(\bullet,0)}(s)}h^{-1}_{\Gamma(\bullet,0)}(s) +
A(\Gamma_1(s,0))=0.\ref\hA$$ Hence we have 
$$\dsv{k_{\Gamma,A,B}(s)}k^{-1}_{\Gamma,B}(s) +
{\Ad}_{h^{-1}_{\Gamma(\bullet,0)}(s)}\int_0^1{ \d t\,\;
B\left(\Gamma_1(s,t),\Gamma_2(s,t)\right)}=0,\ref\equak$$ and
$$\int_0^1{\d t \;\; B\left(\Gamma_1(s,t),\Gamma_2(s,t)\right)}=
\int_0^1{ \d t\;\;{\Ad}_{h^{-1}_{\Gamma(s,\bullet)(t)}} (\sigma^*
B)\left({\Gamma'(s,t)},{\dot\Gamma(s,t)}\right)}.$$ 

The solution  of \equak\ is finally given by
$$k_{\Gamma,A,B}(s')=\Pscr_{s'}
\exp\left\{
-\int_0^{s'} \d s\; {\Ad}_{h^{-1}_{\Gamma(\bullet,0)}(s)} \int_0^1{
\d t\;{\Ad}_{h^{-1}_{\Gamma(s,\bullet)}(t)} (\sigma^* B)\left({
\Gamma'(s,t)},{ \dot\Gamma(s,t)}\right)}\right\}, \ref\holonomy$$
where $\Pscr_{s'}$ denotes path-ordering in the variable $s'$. \vskip
1truecm 

\sezione{Transformation properties of the holonomy as a function of the 
connection} 

In this section we consider a generic manifold $X$ and a principal
$G$-bundle $\pi\colon P_X\to X$ (typically we have in mind either 
$X=M$ or $X=\path$) and  we recall the main properties of
the parallel transport of paths and of the 
holonomy of loops, both seen as functions
on  the space of connections $\Ascr(P_X). $
\blank
Let $\omega\in \Ascr(P_X)$, $\eta\in \T \Ascr(P_X)$,
$\gamma\in \pathv X,$
and  $u\in\pi^{-1} \gamma(0) \subset P_X$.  We consider the
horizontal  lift ${\frak L}(\omega, \gamma, u)$ \lgothic.  To the
path of connections given by  $\kappa\rsa \omega+\kappa\eta,\kappa
\in (-\epsilon,\epsilon)$  we associate
the path in $\pathv{G}$, $\kappa\rsa
g_\kappa=g_\kappa(\gamma, \omega,\eta, u)\in \pathv{G}$ given by
the solution of the following
equation: $${\frak L}(\omega+ \kappa \eta,\gamma,u)(t)={\frak
L}(\omega,\gamma, u)(t) g_\kappa(t),\quad g_{\kappa}(0)=1.\ref\pathG$$ By definition
we have
$$ 
\left[(\omega+\kappa\eta)\left(\dtv{\partrak(t) }\right) \right]=0,
\forall \kappa\in (-\epsilon,\epsilon), \forall t.\ref\varaut$$ 

In this section
we use again a simplified notation for the horizontal lift by
setting $$ {\frak L}(t) \equiv \partrao(t), $$
and the relevant evaluation map
$$ev: I\times {\frak L}\to P.$$
The paths $g_\kappa(t)$ satisfy the following equation in the
variable $t$  $$ g_\kappa^{-1}{\dot g_\kappa} + \kappa{\text
Ad}_{g_\kappa}^{-1} \eta(\dot {\frak L}) 
 =0\ref\eqhkappa.$$
The solution is the
path-ordered exponential
$$g_\kappa(t)=\Pscr\exp\left(-\kappa\int_0^t \d \tau\,\eta (\dot
{\frak L}(\tau) )\right)= \Pscr\exp\left(-\kappa\int_{[0,t]}ev^*\eta
\right),\ref\formalvaria$$ We are
interested in $H(t)\triangleq\dy{\dkov{g_\kappa(t)}}$. 

By differentiating at $\kappa=0$   \eqhkappa, we get
$$H (t)=-\int_{[0,t]}\d t\, \eta (\dot {\frak L}(t) )=
-\int_{[0,t]} ev^*\eta .\ref\varholpath$$ Thus we have 
proved the following

\theorem\varolonomia{For any loop $\gamma\in\loopv X$,   the
logarithmic exterior derivative  of the holonomy,  seen as a
function of the connection $\omega\in\Ascr_X$, is given by
$$ \Hol_\omega^{-1}(\gamma, u)\delta \left(\Hol_\omega(\gamma,
u)\right) (\eta) =   -\int_{I}  ev^*\eta  .\ref\varholloop$$  } 
 
We denote by $\AutP_X$ the group of automorphims of $P_X$ and by
$\autP_X$ the Lie algebra of infinitesimal automorphisms of $P_X$. There is
an action of $\AutP_X$ on $\Ascr_X$ and  a projection (group
homomorphism) $$\rho\colon \AutP_X\to \DiffX\ref\jproj$$ whose kernel is  the
group of gauge transformations $\Gscr(P_X)$.

This projection allows us to define an action of $\AutP_X$ on
$\pathv X$. Hence any $\psi\in\AutP_X$ defines an isomorphism
of  bundles of horizontal paths
$$\psi\colon\papv\omega_X\to  \papv{\psi^*\omega}_X\ref\isopomega$$ 
We now
want 
to discuss the effect of this isomorphism on the parallel transport
and the holonomy.

The isomorphism \isopomega\ satisfies the following
equation $${\frak L}({\psi^*\omega},\rho(\psi^{-1})\circ 
\gamma,\psi^{-1}(u))=\psi^{-1}\left( \partrao\right), \quad
\gamma\in \pathv X.\ref\tripleaction$$  This implies that 
the infinitesimal action of $\autP_X$ 
on $\papv\omega_X$ is just the opposite of the corresponding action on
$\Ascr_X$. For any $Z\in\autP_X$  we compute the corresponding 
Lie derivative
$$L_Z\omega=d_\omega
i_Z\omega+i_ZF_\omega.$$  By setting $\eta=L_Z\omega$ in 
\varholloop, we get
$$\eqalign{  &\dsov{\,\Hol_{\omega+sL_Z\omega}
(\gamma,u)}  =-\Hol_\omega(\gamma,u) \int_{I} 
ev^*\left(i_ZF_\omega+ d_\omega i_Z\omega\right) =\cr &-
\Hol_\omega(\gamma,u)\,\left(\int_{I} ev^*  i_ZF_\omega\right)
-\Hol_\omega(\gamma,u)\left(   i_Z\omega ({\frak L}(1))-  
i_Z\omega({\frak L}(0))  \right).\cr
 }\ref\variationhol$$  If we consider 
the variation of {\it the trace} of
the holonomy (in any representation of $G$),
we have \theorem\vartrace{Let  $\omega\in
\Ascr_X$,  $\gamma\in\loopv X$,  $Z\in \aut P_X$, and  $u\in
\pi^{-1} \gamma(0) $. Then we have $$ \left(\delta
{\tr}\Hol_\omega(\gamma,u)\right)(L_Z\omega)     = -
\tr\left(\Hol_\omega(\gamma,u)\int_{{I}} ev^*  i_ZF_\omega\right) 
\ref\variationtrhol$$   }\proof\vartrace  \blank We have  $$
i_Z\omega ({\frak L}(1))= i_Z\omega ({\frak
L}(0)\Hol_\omega(\gamma,u)) =\Hol_\omega(\gamma,u)^{-1} \left(i_Z
\omega({\frak L}(0)\right) )\Hol_\omega(\gamma,u)  $$   and therefore
$$\tr\left(\Hol_\omega(\gamma,u)  i_Z\omega ({\frak L}(1))-  
\Hol_\omega(\gamma,u) i_Z\omega({\frak L}(0))\right)=0$$  The result
now follows from \variationhol.\endproof 

\cor\zinv{The variation \variationtrhol\ vanishes if the restriction
of $Z$ to the image of $\gamma$ is proportional to the
tangent vector $\dot \gamma$}
\blank
In particular if the loop is an imbedding, then 
the corresponding trace of
the holonomy is invariant under the action of  any
$\psi\in \AutP_X$ connected to the identity  for which $\rho(\psi)\in
\DiffX$  maps the image of the loop  into itself.

\sezione{Holonomy of cylinders and the group of automorphisms of $P$.}

In this and the following 
section we consider loops of paths  and loops of
loops and study the corresponding holonomies as
functions on the space of (special) connections.

We will use the name {\it cylinders} to denote the image
of loops of paths,
even though we are not assuming that such loops of paths
are necessarily imbeddings or 
immersions.

In this section we look for the conditions  which guarantee that
the (trace of the) holonomy of a loop of paths  is invariant
under those automorphisms  of $P$ which 
project onto diffeomorphisms connected to the
identity, that  map the corresponding image (cylinder) into itself.

Since we are considering the action of $\AutP$ on the space of connections $\Ascr$,
it is convenient to work primarily with the bundle $ev^*_0P$ instead
of $\pap,$ for which
the choice of a fixed connection $A$ is required. We will, 
though, make constantly 
use of the isomorphism $J_A\colon ev^*_0P\to\pap$ \isobundles.

Equation \tripleaction\ says 
that the group $\AutP$ of automorphisms of $P$ acts in a natural way on
the bundle $ev^*_0P$. In fact we have  $$P\times \path\ni
(p,\gamma)\rsa (\psi(p),\rho(\psi)(\gamma)),\quad \psi\in \AutP$$ 
with $p\in \pi^{-1} \gamma(0) \Lra \psi(p)\in
\pi^{-1} \rho(\psi)(\gamma(0)).$
\blank 
The group $\AutP$ can be identified with a
subgroup of $\Aut(ev^*_0P).$
The Lie algebra $\autP$ can be accordingly identified with a subalgebra
of $\aut(ev^*_0P).$    

Given now $Z\in\aut P$ and the corresponding element in
$\aut(ev^*_0P)$ which we denote by the same symbol, we want 
to describe $(J_A)_*Z\in\aut(\pap)$ explicitly.
Consider $q\in \pap.$
The path $$t\rsa (\pi q(t), \rho_*Z(\pi q(t)))\ref\autpath$$ is an element of 
$\T\path$. We now lift $\dota$-horizontally \autpath\
(see \tangential )
with initial point $(p, Z(p))\in \T P$. This
lifted path is $((J_A)_*Z)(q)$. For any $t$, $((J_A)_*Z)(q)(t)$ is a vector
in $\T_{q(t)}P$. 
Notice that  
in general $((J_A)_*Z)(q)(t)$ is different from $Z(q(t))$ unless
$t=0$.
\blank
The isomorphism \isobundles\ 
$\ja\colon ev_0^*P\to\pap$ and the corresponding evaluation map
\isobundlesone\ 
$\jsuba\colon ev_0^*P\times I \to P$ allow
us to transform forms on $\pap$
defined by Chen integrals into forms defined on
$ ev^*_0P$. 
The result of performing first Chen integrals and then pulling
back the forms to $ev^*_0P$ via $\ja$ will be represented by the
symbol
$\chena.$ 
In the special case of line-integrals,  we have for a
generic $k$-form $\phi$ on $P$  
$$\chena \phi =\int_I \jsuba^*\phi.$$ 

Then we have the following
\theorem\pullback{The pullback of the connection $\connab$
\connectab\ via $\psi\in \Aut P$ is given by $$\psi^*(\connab )=
ev^*_0\psi^*A + \int_I j_{\psi^*A}^* \psi^*B. \ref\connectab$$
At the infinitesimal level, for any $Z\in \autP$, we have
$$L_Z\connab = ev^*_0 L_Z A +\chena  L_ZB + \chena \{L_ZA;  B\}.
\ref\infinitesimal$$}
 
\proof \pullback
\blank
We have $$L_Z\connab = ev^*_0 L_Z A +\chena  L_ZB + \dkov{}
\chenav{A+\kappa L_ZA}  B\ref\infinitesimal$$ 
If we are given
$\eta\in\form{1}$ and $\zeta\in\form{*}$ we have $$\dkov{}
\chenav{A+\kappa \eta}  \zeta=\dkov{}\int_I
j^*_{A+\kappa\eta}\zeta=\dkov{}\int_I{\Ad}_{g^{-1}_\kappa}j_A^*\zeta$$
where $g_\kappa$ is defined as in 
in \pathG. Now the proof follows from  \varholpath.
\endproof
The curvature $\curvab$ of $\connab$ at $(q,p)$, is  given by
$$ev^*_0F_{A}-\jsuba(1)^*  B+ev^*_0B+  \chena
d_AB + \chena\left\{
B+ F_A; B\right\}.\ref\curvatureab$$

A direct consequence of \vartrace\ is

\theorem\holocyl{Let  $\Gamma\in\loopv{(\path)}$  and $Z\in \autP$.
The trace of the holonomy $\holp $ in $ev_0^*P$ with respect to $\connab $
transforms as follows $$ \delta  \tr \holp( Z)=- \tr\left(\holp
\int_{I} {ev^*i_{Z}{ \curvab}}\right),\ref\streptosyl $$  
with $ev\colon I\times \loopabh\to P$.} 
\blankm
We now compute explicitly \streptosyl. First
we  set
$${\frak L}_{A,B}(s,t)\equiv\loopabh(s)(t),\quad
Z(s,t)\equiv  ((J_A)_* Z)({\frak L}_{A,B}(s,\bullet))(t) \in
\T_{{\frak L}_{A,B}(s,t)}P.\ref\pst$$ 

We can write down \streptosyl\ as follows:
$$\displaylines{ \delta \tr\holp(Z)=    
 -\tr\holp\bigg\{\int_0^1\d s\,  F_A\left(Z(s,0), {{\frak
L}_{A,B}'(s,0)}\right)
   \cr
 + \int_0^1\d s\, \int_0^1 \d t \,{d_AB \left( Z(s,t),{{\frak
L}_{A,B}'(s,t) }  ,{\dot{\frak L}_{A,B}(s,t)} \right)}\cr
+\int_0^1\d s\, \int _0^1 \d t\int_0^t \d \tau\,\left[  (B+
F_A)\left(Z(s,\tau),{\dot{\frak L}_{A,B}(s,\tau)}\right),
 B\left({{\frak L}_{A,B}'(s,t)},{\dot{\frak
L}_{A,B}(s,t)}\right)\right] \cr
-\int_0^1\d s\, \int _0^1 \d t\int_0^t \d \tau\,\left[  (B+
F_A)\left({{\frak L}_{A,B} '(s,\tau)},\dot{\frak L}_{A,B}
(s,\tau)\right), B\left(Z(s,t),{\dot{\frak L}_{A,B}(s,t)}\right)
\right]\bigg\}. \cr}$$

In order to obtain the vanishing of the previous expression, we make
some extra assumptions on the vector field 
$Z$, namely: \blank
\item{\bf A)}$\pi_*Z(s,0)$ is proportional to the tangent
vector to the path of initial
points
$\Gamma(\bullet,0)$ 
  i.e., $ Z(s,0)$ is proportional to ${{\frak L}_{A,B}'(s,0)}$ up to
vertical vectors, \item{\bf B)}
 $\pi_* Z(s,t)$ is a  linear combination
 of $\Gamma'(s,t)$ and $\dot\Gamma(s,t),$ with coefficients that, in
general,  are
functions of $s$ and $t.$
\blank Moreover let $\Sigma$ be a submanifold of $M$ containing
 $\im(\Gamma)$. For the restriction of $P$ to 
 $\Sigma$, we make the following assumptions:
\blank \item{\bf C)} the connection $A$
restricted to the bundle $P_\Sigma$ is reducible to an abelian
subgroup $T$ of $G,$ \item{\bf D)} the form $B\in\form{2}$ restricted
to  $P_\Sigma$ is (simultaneously) reducible to $T$ \blank 
We have finally the following \theorem\vartracabb  {We have
$$\delta \tr\holp(Z)=0\ref\ultima$$ provided either that conditions {\bf A), B),
C), D)} are satisfied  or that conditions {\bf A), B)} are satisfied
together with the extra requirement that on $\Sigma$
we have  either $B=-F_A$ or $B=0$.}
\blank
\vfill\eject
\sezione{Invariance properties of the (trace of the)
$(A,B)$-holonomy}

The space of connections on $ev^*_0P$ of the type $\connab$ is isomorphic
to the affine
space $\Ascr\times\form{2}$ which is acted 
upon by some transformation groups, that arise in the framework of quantum 
field theories of the $BF$ type (see below). 

In this section we want to check under what conditions 
the trace of the $(A,B)$-holonomy is invariant under those 
transformation.
group.  

In quantum field theories
one considers  first of all the
gauge group $\Gscr$. If we divide $\Gscr$ by its center
and consider only irreducible connections,
then $\Gscr$  acts freely on
$\Ascr\times\form{2}$
$$(A,B)g=\left(A^g,{\Ad}_{g^{-1}}B\right).\ref\gaugetrasf$$
We have moreover  the group
$\Gscr_T$ given by the semidirect product
$\Gscr \ltimes \form{1}$, where $\Gscr$ acts on the abelian group
$\form{1}$ via the adjoint action. The group $\Gscr_T$ acts {\it 
non-freely} in two
ways on $\Ascr\times\form{2}$.
The first action is given by the transformation $$(A,B)\rsa (A^g +
\eta, {\Ad}_{g^{-1}}B -d_{A^g}\eta-\frac 12[\eta,\eta]), \quad (g,\eta)\in
\Gscr_T,\ref\actuno$$ while the second action is given by
$$(A,B)\rsa (A^g, {\Ad}_{g^{-1}}B -d_{A^g}\eta), \quad (g,\eta)\in
\Gscr_T.\ref\actdue$$

Before seeing how the above transformation groups act on the holonomy,
we compute
the derivative of the $(A,B)$-holonomy
as a function on $\Ascr\times\form{2}$
at $(\eta,\beta)\in \T\Ascr\times\T\form{2}$

Under the transformation $A\rsa A+\eta, B\rsa B+ \beta$
the connection $\connab=ev^*_0A +\int
\jsuba^*B$ on the bundle $ev^*_0P\to\path$ transforms into
$$ev_0^*(A+\kappa\eta)+\chenav{A+\kappa\eta} B+\kappa\beta .$$ The
corresponding derivative of the holonomy is given by:
$$\gathered\delta \tr \holp(\eta,\beta) =\\ -\tr\bigg(\holp\int_I  
\bigg(ev^*_0\eta+\chena\beta+  \chena\{\eta;
 B\}\bigg)\bigg),\endgathered \ref\varhol$$ 
  
The integral in the r.h.s. of \varhol\ can be written explicitly as
$$\gather\int_0^1  \d s\,\left(\eta({\frak L}_{A,B}'(s,0))+\int_0^1
\d t \,\beta({\frak L}_{A,B}'(s,t),\dot{\frak L}_{A,B}(s,t))
\right)+\\ \int_0^1  \d s\,\int_0^1 \d t\left[\int_0^t
\d\tau\,\eta(\dot{\frak L}_{A,B}(s,\tau)), B({\frak
L}_{A,B}'(s,t),\dot{\frak L}_{A,B}(s,t)) \right],\endgather $$ where
the prime denotes, as usual, the derivative with respect to the
variable $s$ and the dot the derivative
with respect to the  variable $t,$ and
${\frak L}_{A,B}(s,t)$ has been defined in \pst.

A direct consequence of \pullback\ is the following 
\theorem\gauge{Let $\Gamma\colon S^1\times I\to M$  be any loop  in $\path$,
let $p$ be such that $\pi(p)=\Gamma(0,0)$ and let $g\in \Gscr$ be
any gauge transformation.  The trace of the $(A,B)$-holonomy of
$\Gamma$ with initial point $p$ is invariant under  the
transformation $$(A,B)\rsa (A^g,{\Ad}_{g^{-1}}B).$$}

Now we can study the transformation properties of the $(A,B)$-holonomy
under \actuno\ and \actdue\ in the special case when 
we restrict the elements of $\Gscr$ to be the identity.

In this case  \actuno\
becomes the transformation
$$A\rsa A+\eta,\quad B\rsa B-d_A\eta-\frac 12 [\eta,\eta].$$
and we have the following:

\theorem\urka{When $\beta=-d_A\eta-\frac 12[\eta,\eta]$ then \varhol\ becomes
$$\gathered \delta \tr \holp(\eta,-d_A\eta) = \\ -\tr\bigg(\holp\int_I
\bigg(ev^*_1\eta-\frac 12\chena [\eta,\eta]- \chena\{B+F_A;\eta 
\}\bigg).\endgathered 
\ref\varholbis$$}

\proof\urka \blank We have $$\gather\int_0^1\int_0^1 \d s\d t
\left[A,\eta\right] \left(  {\frak L}'_{A,B}(s,t),\dot{\frak
L}_{A,B}(s,t)\right)=
  \int_0^1\int_0^1 \d s\d t \left[A({\frak
L}_{A,B}'(s,t)),\eta(\dot{\frak L}_{A,B}(s,t))\right]=\\
 \int_0^1\int_0^1 \d s\d t \left(\left[A({\frak L}_{A,B}'(s,0))
+\int_0^t \d\tau \, F_A(\dot{\frak L}_{A,B}(s,\tau),{\frak
L}_{A,B}'(s,\tau)),\eta(\dot{\frak L}_{A,B}(s,t)) \right]\right)=\\
 \iint_{I\times I} \d s\d t \left[\int_0^1\d \tau B(\dot{\frak
L}_{A,B}(s,\tau),{\frak L}_{A,B}'(s,  \tau)) +\int_0^t \d\tau
F_A(\dot{\frak L}_{A,B}(s,\tau),{\frak
L}_{A,B}'(s,\tau)),\eta(\dot{\frak
L}_{A,B}(s,t))\right],\endgather$$ where we have used \horzpap.
 Therefore $$\gather
-\int_I\left(\chena
 [A,\eta]-  \chena\{\eta; B\}\right)=\\=-\int_0^1\d s\int_0^1 \d t
\left[\int_0^t \d \tau\,(F_A+B)(\dot{\frak L}_{A,B}(s,\tau),{\frak
L}_{A,B}'(s,\tau)), \eta(\dot{\frak L}_{A,B}(s,t)
 )\right].\endgather$$  Notice also that $$\gather -\int_I\chena
d\eta= \int_0^1\d t\,\left( -\eta(\dot {\frak L}_{A,B}(1,t))+\eta(\dot
{\frak L}_{A,B}(0,t))  
 \right)-\\
 \int_0^1\d s\,\left( \eta(  {\frak L}'_{A,B}(s,0)) - 
 \eta(  {\frak L}'_{A,B}(s,1))\right)=\\ \int_0^1\d t\,\left( \eta(\dot {\frak
L}_{A,B}(0,t)) -
 \holp^{-1} \eta(\dot {\frak L}_{A,B}(0,t))\holp\right)  \\ -\int_0^1\d
s\,\left( \eta(  {\frak L}'_{A,B}(s,0)) -
 \eta(  {\frak L}'_{A,B}(s,1))\right) .\endgather$$ Therefore $$\gather
-\tr\holp\Bigg(\int_0^1  \d s\,\left(\eta({\frak
L}_{A,B}'(s,0))-\int_0^1 \d t \,d_A\eta({\frak
L}_{A,B}'(s,t),\dot{\frak L}_{A,B}(s,t)) \right)+\\ \int_0^1  \d
s\,\int_0^1 \d t\left[\int_0^t \d\tau\,\eta(\dot{\frak
L}_{A,B}(s,\tau)), B({\frak L}_{A,B}'(s,t),\dot{\frak L}_{A,B}(s,t))
\right] \Bigg)=\\ \int_Iev^*_1\eta-\chena\{F_A+B;\eta\}\endgather $$
\endproof

Finally we take into account \urka\ and consider the invariance
properties under \actdue\  of the trace of the holonomy
of a loop of paths $\Gamma\colon S^1\times I\to M.$ 
The previous discussion
yields the following
 
\theorem\urkauno{Corresponding to the action \actdue\ we have $$
\delta \tr \holp(0,d_A\eta) =$$
$$\tr\bigg[\holp\bigg(\int_I  
   \left(\chena \{F_A;\eta\}+\left[\int_I B,\int_I\eta\right] \right)
+$$ $$\int_I{\frak L}_{A,B}(s,1)^*\eta-\int_I{\frak L}
_{A,B}(s,0)^*\eta\bigg)\bigg].$$ }
We now consider {\it loops of loops}. In this case we have:
\cor\aburka{
Let $T$ be an abelian subgroup of $G$. If conditions ${\bold C)}$ and 
${\bold D)}$ of the
previous section are satisfied and if the restriction of
$\eta\in
\form{1}$ to $\Gamma\colon S^1\times S^1\to M$ is also reducible to
$T$, then the  trace of the $(A,B)$-holonomy of the loop
of loops $\Gamma$ is
invariant under \actdue. If, besides the above conditions, 
we have also  $$\int_{I}ev^*\eta=0,\quad 
ev\colon I\times \Gamma(\bullet,0)(p) 
\to P,$$
 then the trace of  the $(A,B)$-holonomy for loops of loops
 is also invariant under
\actuno.}
\blank

The conclusion of this section is that the symmetry \actdue, which arises
from the $BF$ (quantum) field theory, does not leave the trace of the holonomy
invariant, unless some reducibility constraints are imposed on the connection
$A$ and on the field $B$. 

In this sense the transformations \actdue\ represent {\it almost} a
good symmetry for the observable given
by the trace of the $(A,B)$-holonomy.

A good symmetry for the same observable would certainly be represented
by the group of gauge transformations for $\pap$. 
Unfortunately gauge transformations
for $\pap$ do not map special connections into special connections and hence
are not good symmetries for the $BF$ theories. 

In general gauge transformations for $\pap$ map
special connections into special connections plus some extra terms
given by Chen integrals and boundary terms. By neglecting these extra terms 
one obtains exactly \actdue. In this sense the 
transformations \actdue\ are
{\it almost} gauge transformations for $\pap$. 

\sezione{Observables, actions and quantum field theories} An
application of the ideas developed in this paper is the construction
of new {\it observables for quantum field theories} (QFT).

A QFT is described by an action functional, and by observable one
means another functional that is invariant under the same symmetries
that leave the action functional unchanged. A weaker requirement for
the observables is the invariance only {\it on shell} (i.e., upon using
the Euler--Lagrange equations);  in this case the quantization of the
theory requires the use of the Batalin--Vilkovisky formalism
\quot\bat{I.A.\,Batalin and G.A.\,Vilkovisky, {\it Relativistic
S-Matrix of Dynamical Systems with Boson and Fermion Constraints},
\pl{69 B}, 309--312 (1977);}, \quot\frad{E.S.\,Fradkin and 
T.E.\,Fradkina, {\it Quantization of Relativistic 
Systems with Boson and Fermion First- and Second-Class Constraints},
\pl{72 B}, 343--348 (1978);}, 
but we
will discuss this elsewhere.
Throughout this section we will restrict ourselves to considering a 
four-dimensional manifold $M$.

\vskip 12pt 
\sottosezione{Non-topological QFT's} 

The first QFT we consider is the
Yang--Mills theory described by the action functional $$ S_{\text
YM}[A] = ||F_A||^2 = -\int_M\tr(F_A\wedge *F_A), $$ where $*$ is the
Hodge dual with respect to the Riemannian metric on $M$. The
invariance group of the Yang--Mills action functional is the group of gauge
transformations
$A\to A^g$. In this framework we have two
natural elements of $\Omega^2(M,\adP)$ at our disposal, viz., $F_A$
and $*F_A$. Therefore, we may consider the following family of
observables: $$ \Oscr_{\alpha\beta}(\Gamma) =  \tr\Hscr_{(A,\alpha
F_A +\beta *F_A)}(\Gamma,p), \ref\oscr$$ where $\Gamma$ is a path of paths or
loops. \holopath\ guarantees that this is indeed an observable.

Notice that for $\alpha=\beta=0$ the observable reduces to  the
trace of the identity, while for $\alpha=-1$ and $\beta=0$ it yields
the trace of the $A$-holonomy along the boundary $\partial\Gamma$. Taking
$\alpha=0$ and $\beta=1$ ($\beta=-1$) is an interesting choice if
the background connection---i.e., the solution of the 
Euler--Lagrange equations $d_A^*F_A=0$ around which we are 
working---is anti-self-dual (self-dual); in this case, on shell the
observable is  the $A$-holonomy along the boundary of $\Gamma$ but
off shell it depends on $\Gamma$ (see \wilson).

Another family of observables can be obtained by replacing $\Hscr$
by $\Hol$ in the above formula.

As discussed in the Introduction, a 
 physical interpretation of these observables may be the following:
as the Wilson loop---i.e., the trace of the $A$-holonomy---describes 
the displacement of a point-like charge, so the
observable $\Oscr$ describes the displacement of a path-like (or
loop-like) charge, namely of an open  or closed string.

Notice that, even if  the image of
$\Gamma$ represents a smooth surface, the
observable $\Oscr$ depends in general
on its underlying path-of-paths structure.
If, however, we impose assumption {\bf C} of section {\bf 6}
as a {\it boundary condition} for $A$, then $\Oscr$ will depend
only on the surface represented by $\Gamma$ and on the loop of
initial points.

There are other theories that are equivalent to the Yang--Mills
theory,  like the first order Yang--Mills theory 
\halpuno\halpdue,
$$ S_{\text
YM'} = {1\over 4}\,||B||^2 +i \int_M\tr(B\wedge F_A).\ref\bfym $$ In this
case, however, we have at our disposal a bigger family of observables
than those given by \oscr\. In fact
as our form in $\Omega^2(M,\adP)$, we can take a generic
linear combination $$ \alpha F_A +\beta *F_A + \gamma B + \delta *B.
$$ 
Another version of \bfym\ is the so-called $BF$-Yang--Mills theory \bfquattro,
\bfsette, 
where $B$ is replaced by $B-d_A\eta$, $\eta\in\Omega^1(M,\adP)$, in
the above action and, consequently, in the observable. 

The $BF$-Yang--Mills
theory has been extensively studied in \bfsette\ where it has been shown to
be equivalent to the Yang--Mills theory.  
This equivalence makes more interesting the appearance 
of a surface term for Wilson loops.

\sottosezione{Topological QFT's}

Topological Quantum Field Theories are QFT's whose action functional
does not depend on the Riemannian structure of $M$ and so it is expected
to yield topological or smooth invariants as its vacuum expectation values.

We consider  the following TQFT's: \blank \item {1)} the
topological Yang--Mills theory $$ S_{\text tYM} = 
\int_M\tr(F_A\wedge F_A), $$ \item{2)} the $BF$ theory with a
cosmological term $$ S_{BF-BB} =  \int_M\tr(B\wedge F_A) +
{1\over 2}\int_M\tr(B\wedge B), $$ and\blank \item {3)} the pure $BF$ theory $$
S_{BF} =  \int_M\tr(B\wedge F_A). $$

We do not have a non trivial loop-of-loops observable for the topological
Yang--Mills theory. 

As for the $BF$ theory with a cosmological term, we notice that the
symmetries read, at the infinitesimal level, $$ \delta A = d_A\xi
+ \eta, \qquad \delta B = [B,\xi] - d_A\eta, $$ with
$\xi\in\Omega^0(M,\adP)$ and $\eta\in\Omega^1(M,\adP)$. These 
transformations correspond to  \actuno. 

Since the
Euler--Lagrange equations are $B+F_A=0$, then the trace of $\Hol_{(A,B)}$
is {\it almost} invariant on shell. The problem is the presence of
boundary terms in $\eta$, see \varholbis. To get a good on-shell
observable for {\it loops of paths} $\Gamma$, 
we have to eliminate these boundary terms; so we 
may consider $$ \Oscr(\Gamma) =
\tr\left[\Hol_{(A,-F_A)}(\Gamma,p)^{-1}\, \Hol_{(A,B)}(\Gamma,p)
\right]. $$ Notice that on shell this observable is trivial. Off
shell one must add Batalin--Vilkovisky corrections.
Alternatively one can assume conditions {\bf C} and {\bf D} of 
section {\bf 6} as boundary conditions. In this case the above observable
is automatically invariant both on shell and off shell.

In the case of the pure $BF$ theory, the infinitesimal symmetries
are \actdue, i.e. $$ \delta A 
= d_A\xi, \qquad \delta B = [B,\xi] - d_A\eta.
$$ The Euler--Lagrange equations read $F_A=0$, $d_AB=0$.
These 
conditions correspond 
{\it almost} to the flatness of the connection for loops of
paths, the missing requirement being the reducibility of $B$, see
\flatness. 

We have
then a first observable for pure $BF$ theory, namely,
$$ \Oscr(\Gamma) = \tr \Hol_{(A,B)}(\Gamma,p). $$
In fact, by \urkauno\ we get $$ \delta\Oscr(\Gamma) = 0,$$
provided that we assume conditions {\bf C} and {\bf D} of 
section {\bf 6} as boundary conditions. In this case the above observable
is invariant both on shell and off shell.

Another possible choice for pure $BF$ theory is 
given by the observable
$$
\widetilde\Oscr(\Gamma) = \tr\exp\left[
{d\over dt}\Big|_{t=0}\,\Hol_{(A,tB)}(\Gamma,p)
\right].\ref\bftop
$$
On shell (i.e. when $F_A=0$) \smoothomo\ guarantees that \bftop\ is an 
observable that can be rightly associated to the surface spanned by a 
loop of paths.
To compute the transformation properties of this observable, we must consider
the transformation of the holonomy and not of its trace but only up
to the first order in $t$. So $\widetilde\Oscr$ turns out to
be invariant on shell if 
one requires $\eta$ to vanish on the restriction of $P$ over a submanifold
$\Sigma$ containing the image
of $\partial\Gamma$.
To get a good
observable also off shell, i.e., also in the case when $A$ is not flat,
one must add Batalin--Vilkovisky corrections. Notice that \bftop\ is the exact
counterpart of the observable for 3-dimensional $BF$-theory considered in
\quot\CCM{A.S.\,Cattaneo, P.\,Cotta-Ramusino and M.\,Martellini, 
{\it Three-Dimensional $BF$ Theories and the 
Alexander--Conway Invariant of Knots}, 
\np{B 346}, 355--382 (1995);
},\quot\ccfm{{A.S.\,Cattaneo, P.\,Cotta-Ramusino, J.\,Fr\"ohlich,
M.\,Martellini, {\it Topological BF theories in 3 and 4 dimensions},
\jmp {   36},
  6137--6160 (1995); }},
\quot\cjmp{A.S.\,Cattaneo, {\it Cabled Wilson loops in BF theories},
\jmp  { 37}, 3684--3703 (1996).}.

Since the $BF$ theories are {\it topological}---i.e., do not depend 
on the choice of the  Riemannian metric on $M$---one expects that the vacuum
expectation values of the above metric-independent observables
will yield smooth invariants of the image of an imbedded (immersed) 
loop of paths (of loops).

When $M$ is a four-dimensional simply connected manifold,
we conjecture that these invariants are related to the Kronheimer--Mrowka
invariants \km\, of imbedded (immersed) surfaces. Both in 
their theory and in our framework, a special r\^ole is played by connections
that are reducible when restricted to the given surface. Moreover both
in \km\ and in the preliminary perturbative 
calculations of the four-dimensional quantum
$BF$ theory (see \cm, \bfquattro), the reducible connections
(``monopoles on the surface" ) yield loops and surfaces
that are non-trivially linked.

\sezione{Conclusions}
The natural geometrical setting for field theories of the $BF$ type
is a principal bundle on the space of paths (``open strings") or loops 
(``closed strings") of a (four-dimensional) manifold $M$. The fields $A$ and $B$
of the $BF$ theory describe collectively a connection on such principal bundles.

Out of the trace of the corresponding  holonomy one can  define
observables associated to paths (loops) of paths (of loops). 
These  can be seen as
associated to imbedded (or immersed) surfaces only if some extra
conditions are met and if those extra conditions are taken into account in the
calculations of Feynman integrals.
 
The geometrical analysis of $BF$ theories suggests two physically relevant
considerations:
\blank
\item{1.} In those $BF$ theories that are 
related (equivalent) to the Yang--Mills
theory, one can consider  $B$-dependent observables associated to
paths of paths which, when $B$ is a
deformation of the curvature, are a  deformation of the Wilson 
loop along the boundary of the surface spanned by
the path of paths. In other words a deviation from
the non-abelian Stokes formula appears and 
this may be relevant for a correct understanding of the 
problem of quark-confinement.
\item{2.} Four-dimensional topological $BF$ theories yield  invariants of 
the four-manifold. When no $B$-dependent observable  is included, the invariants
to be considered should be related to the
Donaldson invariants. When $B$-dependent 
observables are considered, one expects the corresponding quantum field theory
to yield invariants of imbedded (or immersed)
surfaces (like the Kronheimer--Mrowka invariants).
\blank
Four-dimensional $BF$ theories can then be  considered as 
a sort of gauge theories for loops and paths.
The main difference is the fact
that the action functional is not integrated over the whole space of
paths (loops) but over the original four-manifold $M$. 
As a consequence, the action functional is not invariant under the full 
gauge group of the principal bundle over the path space but is 
only approximatively
invariant 
(i.e. when one neglects boundary terms and higher-order Chen integrals).

The full structure of the gauge group, of the space of connections and
of the space of gauge orbits for paths and loops as well as the
relation with Hochschild (and cyclic) (co)homology, 
will be discussed elsewhere.
\ack{
We thank  A.Belli, L.Bonora, J.D.S Jones, M.Martellini
for useful discussions. P.C.-R. thanks G.
Semenoff for inviting him
to Vancouver B.C. (July 1997,
APCTP/PIms Summer Workshop).}
\immediate\closeout\fileack
                \par
                \null\blankm
                \centerline{\bf Acknowledgments}
                \blankm
                \input ref.tmp2\vfill\eject
\sezione{Appendix: Iterated loop spaces} Most of the construction
in ths paper can be easily iterated, namely we can consider 
principal $G$-bundles  on iterated free path and loop spaces. 
Let us denote 
those by the symbols $\Pscr^nM$ and $\Lscr^nM$.

We describe here the special connections and the relevant
curvatures for iterated path spaces (in the case $n=2$). 

If we are given a connection 
$(A,B)$ on $\pap$, then we can consider connections
on the $G$-principal bundle of $(A,B)$-horizontal paths of paths
$$\Pscr_{(A,B)}^2P\riga{\pi}{20}\Pscr^2M.$$ 
We have the following diagram 
$$\matrix  I\times I\times
\Pscr^2_{(A,B)}P&\riga{ev^{13}}{30}&I\times \pap&\riga{ev}{20}&P\\
\mapdown{id\times
id\times\pi}&&\mapdown{id\times\pi}&&\mapdown{\pi}\\ I\times I\times
\Pscr^2M&\riga{ev^{13}}{20}&I\times \path&
\riga{ev}{20}&M.\endmatrix$$ Here $ev^{13}$ acts on the first and the
third element of the product. 
 
Elements of $\Pscr^2_{A.B}P$ are maps $$\eqalign {Q\colon I\times I &
\to P\cr (s,t)&\rsa Q(s,t)\cr}$$ satisfying  $$
\gather A\left( {\dot Q(s,t)}\right)=0 ,\quad \forall s,t\in I\\
\dy{A\left(  {Q'(s,0)}\right)=\int_0^1 \d t\; B\left(\dot Q(s,t),
{Q'(s,t)}\right), \quad \forall s\in I.}\endgather$$ 

Vectors tangent to $\Pscr^2_{(A.B)}P$ are maps from
$I\times I$ to $\TP,$ which are in turn defined by maps $$ R\colon
(-\epsilon,\epsilon)\times I \times I\to P$$ so that
$$\left.{\partial \over {\partial r}}\right |_{r=0}  {A\left(
{R'(r,s,0)}\right)-\left.{\partial \over {\partial r}}\right |_{r=0}
\int_0^1 dt\; B\left( {\dot R(r,s,t)}, {R'(r,s,t)}\right)=0, \quad
\forall s\in I.}$$ 
 
Following the definitions of sect. {\bf 2}, in order
to define a {\it special
connection} on $\Pscr^2_{(A,B)}$ we need another connection $(\bar A, \bar B)$
and a form
$C\in \Omega^3(M,\adP)$. Here we choose 
$\bar A= A, \bar B=B.$ By considering the double
evaluation map $Ev\colon I\times I \times \Pscr^2_{(A,B)}P \to P$ the
special connection $(A,B,C)$ is explicitly given by: $$Ev^*_{(0,0)} A + \int_I
Ev^*_{(0,\cdot)} B +\int_{I\times I} Ev^* C.$$ 

The space  of special connections considered above
is an affine space modeled on $\Omega^3(M,\adP)$. 

We have: $$\int_{I\times I} Ev^* C=\int_I
ev^{*}_{13}\left(\int_I ev^*C\right).$$ Any tangent vector
$X\in \T_Q\Pscr^2_{A,B}P$ is a map $I\times I\ni (s,t)\rsa \T_{Q(s,t)}P.$ So
we get $$\int_{[0,1]\times[0,1]}\d s \,\d t\, C\left(X(s,t),
{Q'(s,t)}, \dot Q(s,t)\right)= \int_{[0,1]}\d s\left(\int_I
ev^*C\right)\left(X(s,t),Q'(s,t) \right) .$$ 
 
The curvature of a special connection
$(A,B,C)$ is obtained directly
from \curvature\ via  the following replacements $$
\cases A\to ev^*_{0}A+\int_I ev^*_{(\cdot)}B\\ B\to \int_I
ev^*C.\endcases$$ 
\immediate\closeout\fileref
                \par
                \null\blankm
                \centerline{\bf References}
                \blankm
                \input ref.tmp1\vfill\eject
\end